\documentstyle[amstex,amssymb]{amsart}
\newcommand{\f}{\cplx}
\newcommand{\real}{{\bold R}}
\newcommand{\cplx}{{\bold C}}
\newcommand{\zint}{{\bold Z}}
\newcommand{\bs}{{\bold s}}
\newcommand{\bb}{b}
\newcommand{\mm}{{\bold m}}

\newcommand{\ra}{{\bold Q}}
\newcommand{\nat}{{\bold N}}

\newcommand{\cplxN}{\f^N}
\newcommand{\F}{F_{N,n}}
\newcommand{\pn}{^{\otimes n}}
 
\newcommand{\cz}{\f[z_1^{\pm 1},\dots,z_n^{\pm 1}]}
\newcommand{\cx}{\f[x_1^{\pm 1},\dots,x_n^{\pm 1}]}

\newcommand{\ssl}{ {{\frak s} {\frak l}}}
\newcommand{\gl}{ {{\frak g} {\frak l}}}
\newcommand{\glN}{ {\frak g}{\frak l}_N}
\newcommand{\sll}{ {\hat{{\frak s} {\frak l}} }_N }

\newcommand{\wt}{ {\mathrm {wt}}}
\newcommand{\Ker}{ {\mathrm {Ker}}}
\newcommand{\End}{ {\mathrm {End}}}

\newcommand{\vac}{{\mathrm v}{\mathrm a}{\mathrm c}}

\newcommand{\p}{\partial} 

\newcommand{\sprod}[2]{( \; #1 \;, \;#2 \; )}
\newcommand{\sprodd}[2]{\langle \; #1 \;, \;#2 \;   \rangle}
\newcommand{\sproddd}[2]{\langle\langle \; #1 \;, \;#2 \;  \rangle\rangle}
\newcommand{\ep}{\epsilon}
\newcommand{\s}{\sigma}
\renewcommand{\l}{\lambda}
\renewcommand{\o}{\omega}
\newcommand{\sgn}{{\mathrm {sign}}}
\newcommand{\ov}[1]{\overline{#1}}
\newcommand{\un}[1]{\underline{#1}}

\newcommand{\LC}{{\cal L}}

\newcommand{\setN}{\{1,\dots,N\}}

\newcommand{\halmos}{\rule{5pt}{5pt}}
\newcommand{\ba}{\begin{array}}
\newcommand{\ea}{\end{array}}
\newcommand{\bqq}{\begin{equation}}
\newcommand{\eqq}{\end{equation}}
\newcommand{\bqa}{\begin{eqnarray}}
\newcommand{\eqa}{\end{eqnarray}}
\newcommand{\bqas}{\begin{eqnarray*}}
\newcommand{\eqas}{\end{eqnarray*}}

\numberwithin{equation}{section}

\newtheorem{prop}{\bf Proposition}
\newtheorem{df}[prop]{\bf Definition}

\newtheorem{thm}[prop]{\bf Theorem}

\newtheorem{lemma}[prop]{\bf Lemma}
\newtheorem{cor}[prop]{\bf Corollary}
\newenvironment{rmk}{\noindent{\bf Remark}\hskip 5pt}{}

\newenvironment{eg}{\noindent\addtocounter{prop}{1}{\bf Example
\theprop}\hskip 5pt}{}
\numberwithin{prop}{section}
\begin{document}

\title[Symmetric functions and the Yangian decomposition ]
 {Symmetric functions and the Yangian decomposition of   
 the Fock and Basic modules of the affine Lie algebra $\sll$ }
\author[Symmetric functions and the Yangian decomposition ]{Denis Uglov}
\address[ ]{Research Institute for Mathematical Sciences,
Kyoto University, Kyoto 606, Japan   }
\email{ duglov@@kurims.kyoto-u.ac.jp }  


\begin{abstract}
The decompositions of the Fock and Basic modules of the affine Lie algebra  $\sll$ into irreducible submodules of the Yangian algebra $Y(\gl_N)$ are  constructed. Each of the irreducible submodules admits the  unique up to normalization eigenbasis of the maximal commutative subalgebra of the Yangian. The elements of this eigenbasis are identified with specializations of Macdonald symmetric functions where both parameters of these functions  approach an $N$th primitive root of unity.  
\end{abstract}

\maketitle


\section{Introduction}

In 1992 Ha, Haldane, Talstra, Bernard and Pasquier \cite{HHTBP} discovered that each of the level-1 irreducible highest weight modules of the affine Lie algebra $\hat{\ssl}_2$ admits an action of the Yangian $Y(\gl_2).$ Since that time there have been a number of further works devoted to this subject. Notably in \cite{BPS} and \cite{BLS1} the decomposition of the irreducible highest weight modules of  $\hat{\ssl}_2$  under the action of the Yangian was constructed and used to derive a new class of character formulas for these modules. These results were generalized to the case of $\hat{\ssl}_N $ in the works \cite{Schoutens},\cite{BS}.  

In the present article we consider one particular  aspect of this intriguing and still not completely understood subject. Namely, we give the explicit construction of the  so called Gelfand-Zetlin bases in the irreducible submodules of the Yangian action on level-1 integrable highest weight modules of  $\hat{\ssl}_N .$

The Yangian algebra $Y(\gl_N)$ comes equipped with a distinguished maximal commutative subalgebra, $A(\gl_N),$ which is sometimes called the {\em Gelfand-Zetlin algebra} of $Y(\gl_N).$ A finite-dimensional irreducible Yangian module is called {\em tame} if $A(\gl_N)$ acts on this module semisimply. The study of tame Yangian modules was initiated by Cherednik in \cite{Cherednik} and continued by Nazarov and Tarasov \cite{NT} who, in particular, gave a classification of tame modules and an explicit construction for any given tame module of $Y(\gl_N).$ By definition a tame $Y(\gl_N)$-module admits an eigenbasis of the commutative algebra $A(\gl_N).$ Moreover, one of the results of \cite{NT} is that such an eigenbasis is unique up to normalization of its elements. In this way each tame $Y(\gl_N)$-module is equipped with  a canonical basis, called the {\em Yangian Gelfand-Zetlin basis} of this module.

First, we consider the Fock space module of the affine Lie algebra $\hat{\ssl}_N$ \cite{JM,KR}. We show, that the Fock space admits a one-parameter family of $Y(\gl_N)$-actions. At generic values of the parameter the decomposition of the Fock space into Yangian submodules is irreducible, and each of the components of this decomposition is tame. Moreover, each of these components  is isomorphic to a tensor product of the so-called {\em fundamental} $Y(\gl_N)$-modules. The union of the Yangian Gelfand-Zetlin bases of the irreducible Yangian submodules gives a certain basis of the Fock space, and the main problem which we deal with in this paper is to give a detailed  description of this basis. 

It is well-known \cite{JM,KR}, that the Fock space module of $\hat{\ssl}_N$ can be realized as the linear space of symmetic functions \cite{MacBook}. This realization is often referred to as the principal bosonization of the Fock space.  The union of the Yangian Gelfand-Zetlin bases gives rise to a basis in the space of symmetric functions, and the natural   question  we ask is which  symmetric functions are elements of this basis. 

The answer is provided in terms of the Macdonald symmetric functions \cite{MacBook}. A Macdonald symmetric function depends on two parameters $q$ and $t.$ When both $q$ and $t$ approach an $N$th primitive root of unity in a controlled way, the Macdonald symmetric function degenerates into what we call the Jack($\gl_N$) symmetric function. This symmetic function depends on one parameter and may be regarded as an analogue of the Jack symmetric function to which it reduces when $N=1.$ 

The basis of the Jack($\gl_N$) symmetric functions labelled by all partitions is precisely the union of the Yangian Gelfand-Zetlin bases associated with the decomposition of the Fock space into $Y(\gl_N)$-submodules.
Each partition is uniquely determined by the {\em coordinate configuration} and {\em spin configuration} associated with this partition. The set of all coordinate configurations labels irreducible components of the Yangian action, and the set of all partitions with a fixed coordinate configuration labels the Yangian Gelfand-Zetlin basis of the Yangian submodule which corresponds to this coordinate configuration.

Next, we consider the basic module of the  affine Lie algebra $\hat{\ssl}_N.$ This module is realized as a certain subspace of the  Fock space \cite{JM,KR}. When the parameter of the Yangian action on the Fock space is fixed to an appropriate non-generic value, this action can be restricted onto the basic module. We would like to emphasize that it is not clear to us whether this Yangian action coincides with that of \cite{Schoutens} apart from the case $N=2$.

  The decomposition of the basic module relative to the Yangian action is again irreducible and each of the irreducible components is tame. The union of the Yangian  Gelfand-Zetlin bases is now identified with the set of all Jack($\gl_N$) symmetric functions  labelled by the {\em $N$-regular partitions.} Moreover the parameter in these symmetric functions is fixed in a certain way.

The irreducible components of the Yangian action on the basic module can be labelled by the ribbon skew Young diagrams, this labelling was proposed in the paper \cite{KKN}. We point out which $N$-regular partitions parameterize the Yangian Gelfand-Zetlin basis in a Yangian submodule that  corresponds to a given ribbon diagram and compute the associated Drinfeld polynomials. The obtained decomposition of the basic module into irreducible Yangian submodules agrees with the conjecture made in the work \cite{KKN}.

\vspace{1cm}
\begin{large} {\bf Acknowledgments } \end{large}
The present article is an extensively expanded version of the series of lectures which I gave at RIMS, Kyoto University in February-March 1996. I am grateful to Professors M.Kashiwara and T.Miwa who gave me the opportunity to deliver these lectures and gave me ample time to incorporate much new material afterwards. Special thanks are due to Kouichi Takemura who took the notes of the lectures and with whom I subsequently collaborated on several topics related to the subject of these notes. A number of results  contained in the present article were first obtained in our joint works \cite{TU1},\cite{TU},\cite{STU}.

\section{The Spin Sutherland Model  }
We consider the spin generalization of the Sutherland Model \cite{Sutherland} which was proposed in \cite{Cherednikadv} and \cite{BGHP}. This Model describes $n$ quantum particles with coordinates $ y_1,\dots,y_n $ moving along a circle of length $L$ $( 0 \leq y_i \leq L)$. Each particle carries a spin with $N$ possible values, and the dynamics of the Model are governed by the Hamiltonian   
\begin{equation}
\ov{H}_{\beta,N} = -\frac{1}{2}\sum_{i=1}^n \frac{\p^2}{\p y_i^2} + \frac{\pi^2}{2L^2}\sum_{1\leq i\neq j\leq n}\frac{\beta(\beta + P_{ij})}{\sin^2\frac{\pi}{L}(y_i - y_j)}. \label{eq:hamiltonian}
\end{equation}
In this Hamiltonian the symbol $P_{ij}$ stands for the spin exchange operator for particles $i$ and $j,$ and the $\beta$ is the coupling constant. We will take the $\beta$ to be a positive real number. 

The expression (\ref{eq:hamiltonian}) is only a formal one. To define it as an operator acting on a Hilbert space it is convenient to make a gauge transformation of (\ref{eq:hamiltonian}) by taking   $W = \prod_{1\leq i<j\leq n}\sin\frac{\pi}{L}(y_i - y_j)$ and defining the gauge-transformed Hamiltonian $H_{\beta,N}$ by 
\begin{equation}
H_{\beta,N} = \frac{L^2}{2\pi^2} W^{-\beta} \ov{H}_{\beta,N} W^{\beta}.  \label{eq:H0}
\end{equation}
If we  set $z_j = \exp(\frac{2\pi i}{L}y_j),$ then $H_{\beta,N}$ assumes the form 
\begin{gather}
H_{\beta,N} = \frac{L^2}{2\pi^2} W^{-\beta} \ov{H}_{\beta,N} W^{\beta} = \label{eq:H}\\ 
= \sum_{i=1}^n D_i^2 + \beta\sum_{i=1}^n(2i-n-1)D_i + 2\beta\sum_{1\leq i<j\leq n}\theta_{ij}\left(D_i - D_j + \theta_{ji}(P_{ij} + 1)\right) + \frac{\beta^2n(n^2-1)}{12}, \nonumber 
\end{gather}
where we defined: $D_i = z_i{\p}/{\p z_i}$ and $\theta_{ij}=z_i/(z_i-z_j)$.

The gauge-transformed Hamiltonian $H_{\beta,N}$ acts on the linear space 
\begin{equation}
\F = \left( \cz \otimes  (\cplxN )\pn   \right)_{antisymm} \label{eq:Fsp}
\end{equation}
where $antisymm$ stands for the total antisymmetrization. In what follows we will always be working with the gauge-transformed Hamiltonian rather than with the Hamiltonian (\ref{eq:hamiltonian}), and will regard the $\F$ as the space of quantum states of the Spin Sutherland Model. 

Now let us give a  complete description of the space of states $\F$, and introduce on this space a scalar product.

\subsection{Preliminary remarks and notations} \label{sec:preliminary}
Let $N$ be a positive integer, it will have  the meaning of the  number of spin degrees of freedom of each particle in the Spin Sutherland Model. For any integer $k$ define the unique $\un{k}$ $\in$ $\setN$ and the unique  $\ov{k}$ $\in$ $\zint$ by setting $k = \un{k} - N \ov{k}$. And for a $k$ $=$ $(k_1,k_2,\dots,k_n)$ $\in$ $\zint^n$ set $\un{k}$ $=$ $(\un{k_1},\un{k_2},\dots,\un{k_n})$, $\ov{k}$ $=$ $(\ov{k_1},\ov{k_2},\dots,\ov{k_n}).$ 

 For any sequence $k$ $=$ $(k_1,k_2,\dots,k_n)$ $\in$ $\zint^n$ let $|k|$ be the  weight: $|k|$ $=$ $k_1 + k_2 + \cdots + k_n$, and define the partial ordering ( the {\em natural} or the {\em dominance } ordering \cite{MacBook} ) on $\zint^n$ by setting for any two distinct $k,l$ $\in$ $\zint^n$:   
\begin{gather} k > l  \\ 
\text{ iff } \quad |k| = |l|, \quad \text{ and} \quad  k_1+\cdots+k_i \geq l_1+\cdots+l_i\quad \text{ for all} \quad i=1,2,\dots,n. \nonumber 
\end{gather}
For $r\in \nat$ let $\LC^{(r)}_n$ be a subset of $\zint^n$ defined as
\begin{equation}
\LC^{(r)}_n = \{ k = (k_1,k_2,\dots,k_n) \in \zint^n \; | \; k_i \geq k_{i+1} \quad \text{and} \quad \forall s \in \zint \; \#\{k_i \: | \: k_i = s\} \leq r \}. \label{eq:LC}
\end{equation}
In particular the $\LC^{(n)}_n$ is the set of non-increasing sequences of $n$ integers and   the $\LC^{(1)}_n$ is the set of {\em strictly decreasing } sequences, i.e. such $k$ $ = $ $ (k_1,k_2,\dots,k_n) \in \zint^n$ that $ k_i > k_{i+1}$.    

Let $V$ = $\cplxN$ with the basis $\{ v_1,v_2,\dots,v_N\},$ and for a formal variable $z$ let $V(z)$ $=$ $ \cplx[z^{\pm 1}]\otimes V $ with the basis $\{ u_k\:|\: k\in \zint\}$ where $   u_k = z^{\ov{k}}\otimes v_{\un{k}}.$ For monomials in the vector spaces $\cz$, $ V\pn$ and $ V(z)\pn$ $=$ $\cz\otimes(V\pn)$ we will use the convention of multi-indices:  
\begin{eqnarray}
&z^t  =  z_1^{t_1}z_2^{t_2}\cdots z_n^{t_n}, \qquad & t=(t_1,t_2,\dots,t_n) \in \zint^n; \\
&v(a) =  v_{a_1}\otimes v_{a_2}\otimes\cdots \otimes v_{a_n},\qquad &  a=(a_1,a_2,\dots,a_n) \in \setN^n;  \\ 
&u_k  =   u_{k_1}\otimes u_{k_2}\otimes \cdots \otimes u_{k_n}, \qquad & k = (k_1,k_2,\dots,k_n) \in \zint^n. 
\end{eqnarray}

\subsection{The space of states $\F$ and its wedge basis} 
Let $K_{ij}$ be the permutation operator for variables $z_i$ and $z_j$ in $\cz$ ( the operator of coordinate permutation ), and let $P_{ij}$ be the operator exchanging $i$th and $j$th factors in the tensor product $V\pn$ ( the operator of spin permutation ).

Let $A_n$ be the antisymmetrization operator in $V(z)\pn$: 
\begin{equation}
A_n (u_{k_1}\otimes u_{k_2}\otimes \cdots \otimes u_{k_n}) = \sum_{w\in S_n}\sgn(w) u_{k_{w(1)}}\otimes u_{k_{w(2)}}\otimes \cdots \otimes u_{k_{w(n)}}, \label{eq:antisymm}
\end{equation}
where $S_n$ is the symmetric group of order $n$. We will  use the notation $\hat{u}_k$ $=$ $u_{k_1}\wedge u_{k_2}\wedge \cdots \wedge u_{k_n} $ for a vector of the form (\ref{eq:antisymm}), and will  call such a vector {\em a wedge}. A wedge $\hat{u}_k$ is {\em normally ordered } if and only if $k \in \LC^{(1)}_n$, that is $k_1>k_2>\dots >k_n$. Let $F_{N,n}$ be the image of the operator $A_n$ in $V(z)\pn$. Then $F_{N,n}$ is spanned by  wedges and the normally ordered wedges form a basis in $F_{N,n}$. Equivalently the vector space $\F$ is defined as the linear span of all vectors $f$ $\in$ $\cz \otimes (\otimes^n V) $ such that for all $1\leq i\neq j \leq n$ one has   
\begin{equation}
                   K_{ij}f = -P_{ij}f.  \label{eq:K=-P}
\end{equation}
This is the meaning of the total antisymmetrization in the equation (\ref{eq:Fsp}).

\subsection{Scalar product} \label{sec:scalar}
Here we define a scalar product on the space of states $\F.$ Our definition has three  steps. First we define scalar products on the vector spaces $V\pn$ and $\cz$ separately. Then we define a scalar product on the tensor product $V(z)\pn$ $=$ $\cz \otimes ( V\pn).$ Finally we define a scalar product on $\F$ considered as a subspace of $V(z)\pn.$  

On $V\pn$ define a sesquilinear, i.e., $\f$-anti-linear in the first argument and $\f$-linear in the second argument,  scalar product ${\sprod{\cdot}{\cdot}}_N$ by requiring that monomials in $ V\pn$ be orthonormal: 
\begin{equation}
 {\sprod{v(a)}{v(b)}}_N = \delta_{ab}, \qquad a,b \in \setN^n. \label{eq:sps}
\end{equation}

For a Laurent polynomial $f(z)\in \cz$ let $[f(z)]_1$ denote the constant term in $f(z),$ and let $\ov{f(z)}$ be the Laurent polynomial with complex conjugated coefficients. For a non-negative real number  $\delta$ and $ w_1,\dots,w_n \in \f$ define the weight function $\Delta(w;\delta)$ as follows 
\begin{equation}
\Delta(w;\delta) = \prod_{1\leq i\neq j \leq n}(1 - w_i w_j^{-1})^{\delta}.
\end{equation}
This weight function is obviously a symmetric Laurent polynomial when $\delta$ is a non-negative integer.

Now for all $f(z),g(z)$ $\in$  $\cz$ define a sesquilinear scalar product ${\sprod{\cdot}{\cdot}}_{\delta}'$ by setting 
\begin{equation}
{\sprod{f(z)}{g(z)}}_{\delta}' = \frac{1}{n!}\prod_{j=1}^n\int\frac{dw_j}{2\pi i w_j}   \Delta(w;\delta) \ov{f(w^{-1})}g(w)  \label{eq:spc}
\end{equation}
where the integrations are taken over the unit circle in the complex plane. 

On the linear space $V(z)\pn$ $=$ $\cz \otimes (V\pn)$ we define a scalar product ${\sprod{\cdot}{\cdot}}_{\delta,N}'$ as the composition of the scalar products (\ref{eq:spc}) and (\ref{eq:sps}), i.e. for  $f(z),g(z)$ $\in$ $\cz ;$ $u,v$ $\in$ $\otimes^n V$ we set:  
\begin{equation}
{\sprod{f(z)\otimes u}{g(z)\otimes v}}_{\delta,N}' ={\sprod{f(z)}{g(z)}}_{\delta}' {\sprod{u}{ v}}_{N},\label{eq:spcv}
\end{equation}
and extend the definition on all vectors by requiring that ${\sprod{\cdot}{\cdot}}_{\delta,N}'$ be sesquilinear.

Finally, on the subspace $\F \subset V(z)\pn$ a scalar product  ${\sprod{\cdot}{\cdot}}_{\delta,N}$ is defined as the restriction of  the scalar product ${\sprod{\cdot}{\cdot}}_{\delta,N}'$. Note, that the normally ordered wedges form an orthonormal basis of $\F$ relative to this scalar product when $\delta =0$, i.e. we have 
\begin{equation}
{\sprod{ \hat{u}_k }{\hat{u}_l}}_{0,N} = \delta_{kl} \qquad (k,l \in \LC^{(1)}_n).
\end{equation}

\section{The spectrum and an eigenbasis of the Hamiltonian $H_{\beta,N}$} \label{sec:eigenbasis1}

In this section we compute the eigenvalue spectrum of the Hamiltonian (\ref{eq:H}) and describe the corresponding eigenbasis of $\F.$ First of all we observe, that on the space $\F$ we may, as implied by (\ref{eq:K=-P}), replace the spin permutation operators $P_{ij}$ that appear in (\ref{eq:H}) with the operators $-K_{ij}.$ This gives 
\begin{gather}
H_{\beta,N} = H_{\beta,N}^K :=  \label{eq:HK}\\ = \sum_{i=1}^n D_i^2 + \beta\sum_{i=1}^n(2i-n-1)D_i + 2\beta\sum_{1\leq i<j\leq n}\theta_{ij}\left(D_i - D_j - \theta_{ji}(K_{ij} - 1)\right) + \frac{\beta^2n(n^2-1)}{12}\nonumber
\end{gather}
where the equality is understood as that of operators acting on $\F.$ Notice that the right-hand side of this equality commutes with all coordinate permutation operators $K_{ij}.$

Let $k=(k_1,k_2,\dots,k_n)\in \LC^{(1)}_n$. Note that $k_1 > k_2 > \dots > k_n$ implies, in particular, that $\ov{k_1}\leq \ov{k_2}\leq \dots \leq \ov{k_n}$. We wish to compute the action of the  Hamiltonian $H_{\beta,N}$ on the normally  ordered wedge $\hat{u}_k$ $=$ $u_{k_1}\wedge u_{k_2} \wedge \dots \wedge u_{k_n}$ $=$ $A_n (u_k)$. Since $K_{ij}$ commute with the $H_{\beta,N}^K,$ the last operator commutes with the $A_n,$ and we have 
\begin{equation}
H_{\beta,N}.\hat{u}_k  = H_{\beta,N}^K .A_n(u_k) = A_n(H_{\beta,N}^K .u_k ).
\end{equation}
The action of the $H_{\beta,N}^K$ on the $u_k$ is easy to compute, for the $H_{\beta,N}^K$ acts non-trivially only on the first factor in $u_k = z^{\ov{k}}\otimes v(\un{k}).$ The computation  gives 
\begin{equation}
H_{\beta,N}\hat{u}_k = E(k)_{\beta,N}\hat{u}_k + 2\beta\sum_{1\leq i<j\leq n}h_{ij}\hat{u}_k, \label{eq:Hw=}
\end{equation}
where 
\begin{gather}
E(k)_{\beta,N}=\sum_{i=1}^n \ov{k_i}^2 + \beta\sum_{i=1}^n(2i-n-1)\ov{k_i} + \frac{\beta^2n(n^2-1)}{12}, \label{eq:Heigenvalue} \\
\text{and} \quad h_{ij}(u_{k_1}\wedge\dots\wedge u_{k_i}\wedge \dots\wedge u_{k_j} \wedge \dots \wedge u_{k_n}) =\label{eq:hij} \\ = \sum_{r=1}^{\ov{k_j}-\ov{k_i}-1}(\ov{k_j}-\ov{k_i}-r)(u_{k_1}\wedge\dots\wedge u_{k_i-Nr}\wedge \dots\wedge u_{k_j+Nr} \wedge \dots \wedge u_{k_n}). \nonumber 
\end{gather}
Normally ordering the wedges in the right-hand side of (\ref{eq:Hw=}) we find from (\ref{eq:hij}) that 
\begin{equation}
H_{\beta,N}\hat{u}_k = E(k)_{\beta,N}\hat{u}_k + \sum\begin{Sb} l \in \LC^{(1)}_n,\; 
\ov{l} > \ov{k}, \; l < k \end{Sb} h_{kl}^{(\beta,N)} \hat{u}_l, \label{eq:triang1}
\end{equation}
with certain real coefficients  $h_{kl}^{(\beta,N)}$.

Now recall from \cite{MacBook} that for a positive $\beta$ we have: $ E(k)_{\beta,N}\neq E(l)_{\beta,N}$ when $\ov{l} > \ov{k}$. Then (\ref{eq:triang1}) leads to
\begin{prop} \label{prop1}
{\em (i)} For any $k\in\LC^{(1)}_n$ there is a unique eigenvector $\Psi_k^{(\beta,N)}$ of $H_{\beta,N}$ with the following expansion in the basis of normally ordered wedges  
\begin{equation}
 \Psi_k^{(\beta,N)} =\hat{u}_k + \sum\begin{Sb} l \in \LC^{(1)}_n,\; 
\ov{l} > \ov{k} \end{Sb} \psi_{kl}^{(\beta,N)} \hat{u}_l  \qquad ( \psi_{kl}^{(\beta,N)}  \in \real ). \label{p1i}  
\end{equation}

{\em (ii)} The coefficient  $\psi_{kl}^{(\beta,N)}$ vanishes unless $l < k.$

{\em (iii)} The eigenvalue of $H_{\beta,N}$ on  this eigenvector is $E(k)_{\beta,N}.$ 
\end{prop}

Note that for any integer $M$ the wedge
\begin{equation}
\vac(M) = u_M\wedge u_{M-1}\wedge \cdots \wedge u_{M-n+1} \label{eq:vac}
\end{equation}
is an eigenvector of the Hamiltonian $H_{\beta,N}$ as implied by either (\ref{p1i}) or (\ref{eq:Hw=}, \ref{eq:hij}). We will call any vector of the form (\ref{eq:vac}) a {\em vacuum } vector. 

From the above proposition we see that $\{ \Psi_k^{(\beta,N)}\: | \:k \in \LC^{(1)}_n\}$ is a basis of the space of states $\F,$ since the $\{ \hat{u}_k \:| \:k \in \LC^{(1)}_n\}$ is. The spectrum of the Hamiltonian is highly degenerate because the eigenvalue  $E(k)_{\beta,N}$ depends only on $\ov{k}.$ This means that there are other eigenbases of the $H_{\beta,N}.$ The eigenbasis  $\{ \Psi_k^{(\beta,N)} \}$ is not orthogonal relative to the scalar product ${\sprod{\cdot}{\cdot}}_{\beta,N},$ and one of our  tasks will be to construct an orthogonal eigenbasis by using the Yangian symmetry of the Spin Sutherland Model. In this construction the basis  $\{ \Psi_k^{(\beta,N)} \}$ will appear in a supplementary role, and the triangularity with respect to the dominance order of its expansion in terms of the normally ordered wedges expressed by Proposition \ref{prop1} (i),(ii)  will be important.

\section{The Cherednik-Dunkl operators} \label{sec:gauge}

There is an intimate relationship between the gauge-transformed Hamiltonian $H_{\beta,N}$ and the representation of the degenerate affine Hecke algebra given by the Cherednik-Dunkl operators and the operators of coordinate permutation.  

The degenerate affine Hecke algebra of the type $\gl_n,$ ${\bold H'},$ is the unital associative algebra over $\f$ with generators 
\begin{equation}
s_i, s_i^{-1}\quad (i=1,\dots,n-1);\qquad \delta_j \quad (j=1,\dots,n);
\end{equation}
and the following relations 

\begin{gather}
 s_i s_i^{-1} = s_i^{-1} s_i = 1,\qquad  s_{i}^2 = 1,\qquad    s_{i}s_{j}=s_{j}s_{i} \quad \text{ if $ |i-j| > 1$}, \label{eq:dah1}\\
  s_{i}s_{i+1}s_{i} = s_{i+1}s_{i}s_{i+1}, \label{eq:dah2}\\ 
 \delta_j \delta_k = \delta_k \delta_j  , \qquad  
s_{i}\delta_i - \delta_{i+1}s_{i} = 1,\qquad  s_{i}\delta_j = \delta_j s_{i} \quad \text{ if $j\neq i,i+1 $}. \label{eq:dah3}
\end{gather}
The subalgebra of  ${\bold H'}$ generated by the elements $s_i, s_i^{-1}$ is isomorphic to the symmetric group $S_n.$ Note also that the map
\begin{equation}
s_i \mapsto -s_i , \qquad \delta_j \mapsto - \delta_j  \label{eq:d->-d}
\end{equation}
extends to an automorphism of the algebra ${\bold H'}.$

Let us recall, that the Cherednik-Dunkl operators \cite{Cherednikinv} are defined as follows  
\begin{equation}
d_i(\beta)=\beta^{-1}D_i +n - i + \sum_{i<j\leq n}\theta_{ji}(K_{ij}-1)-\sum_{i>j\geq 1}\theta_{ij}(K_{ij}-1), \qquad (i=1,2,\dots,n). \label{eq:CD}
\end{equation}
These operators act on $\cz$ and together with the operators of the coordinate permutations $K_{ij}$ give a representation of  the degenerate affine Hecke algebra. This representation is  defined by the assignment
\begin{equation}
\delta_j \mapsto d_j(\beta), \qquad s_i \mapsto K_{ii+1}. 
\end{equation}
Note that the space $\cz$ with the above action is simultaneosly left and right module of  ${\bold H'}$ as implied by the defining relations (\ref{eq:dah1})-(\ref{eq:dah3}).

The relation between the gauge-transformed Hamiltonian and the Cherednik-Dunkl operators is expressed by (\ref{eq:HK}) and the equality \cite{Cherednikinv}
\begin{gather}
H_{\beta,N}^K = \beta^2\sum_{i=1}^n \left( d_i(\beta) + \frac{1-n}{2}\right)^2.  \label{eq:d=H}
\end{gather}

The following proposition is obtained by a straightforward computation.
\begin{prop} The Cherednik-Dunkl operators are self-adjoint relative to the scalar product ${\sprod{\cdot}{\cdot}}_{\beta}'.$ 
\end{prop}

\begin{cor} The Hamiltonian $H_{\beta,N}$ is self-adjoint relative to the scalar product ${\sprod{\cdot}{\cdot}}_{\beta ,N}.$ 
\end{cor}

Let us make explicit in our notations the dependence of the Cherednik-Dunkl operators on the number of variables $z_1,\dots,z_n.$ We will write $d_i(\beta)^{(n)} $ for the operator $d_i(\beta)$ (\ref{eq:CD}). For $m<n$ there is a natural action of the $d_i(\beta)^{(m)} $ on the space $\cz$ which is non-trivial only with respect to the first $m$ variables.

\begin{lemma} \label{l:dm=dn}
Let the sequence $t=(t_1,\dots,t_n)\in \zint^n$ be such that $t_1,t_2,\dots,t_m< t_{m+1}=t_{m+2}=\cdots = t_n$ for some $1\leq m < n.$ With these notations let $L(m,n,t_n)$ be the linear span of monomials $z^r=z_1^{r_1}\cdots z_n^{r_n}$ such that $r_i \leq t_n$ for all $i=1,\dots,n,$ and $\#\{r_i | r_i = t_n\}$ $<$ $n-m.$   

Then $d_i(\beta)^{(n)}$ leave $L(m,n,t_n)$ invariant and, moreover, the following identities hold
\begin{align}
& d_i(\beta)^{(n)}.z^t \equiv d_i(\beta)^{(m)}.z^t \bmod L(m,n,t_n) \qquad( i=1,\dots,m), \label{eq:dm=dn1} \\
& d_i(\beta)^{(n)}.z^t \equiv (\beta^{-1}t_n - i +n + m )z^t\bmod L(m,n,t_n) \qquad( i=m+1,\dots,n).\label{eq:dm=dn2} 
\end{align}
\end{lemma}

\subsection{Non-symmetric Jack polynomials}

For any $t \in \zint^n$ let $t^+$ denote the element of the set $\LC_n^{(n)}$ (\ref{eq:LC}) obtained by arranging parts of $t$ in non-increasing order. The following proposition summarizes results of the papers \cite{Opdam} and  \cite{Macdonald1} which we use in this article.
\begin{prop} \mbox{} 

{\em (i)} In the space $\cz$ there is the unique up to normalization common eigenbasis $\{E_t^{(\beta)}(z)| t\in \zint^n\}$ of the Cherednik-Dunkl operators.

{\em (ii)} The Laurent polynomials $ E_t^{(\beta)}(z)$ have a triangular expansion in the monomial basis of $\cz.$ Which means that one has 
\begin{equation}
 E_t^{(\beta)}(z) = z^t + \sum_{r \prec t}e^{(\beta)}_{tr} z^r, \label{eq:nsJtri}
\end{equation}
where $e^{(\beta)}_{tr}$ are real  coefficients and
\begin{equation}
r \prec t \quad \Leftrightarrow \quad \begin{cases} t^+ > r^+  &  \text{or} \\ 
   t^+ = r^+  & \text{and the last non-zero difference $t_i - r_i$ is negative.}\end{cases} \nonumber 
\end{equation}

{\em (iii)} The Laurent polynomials $ E_t^{(\beta)}(z)$  are eigenvectors of the Cherednik-Dunkl operators: 
\begin{align}
& d_i(\beta)E_t^{(\beta)}(z) = f_i(t;\beta)E_t^{(\beta)}(z), \qquad (i=1,2,\dots,n),  \\ 
& \text{where} \quad f_i(t;\beta) = \beta^{-1}t_i +n - \rho_i(t), \label{eq:d-eigenvalue}\\ 
& \text{and} \quad \rho_i(t) = \#\{ j \leq i \: | \: t_j \geq t_i \} + \#\{ j > i \: | \: t_j > t_i \}.
\end{align}

{\em (iv)} The action of the operators of coordinate permutation in the basis  $\{ E_t^{(\beta)}(z)\}$ is given by the following formulas
\begin{align}
& K_{ii+1}E_t^{(\beta)}(z) = {\cal A}_i(t)E_t^{(\beta)}(z) + {\cal B}_i(t)E_{(ii+1)t}^{(\beta)}(z),   \label{eq:KactionPhi1}\\ 
\intertext{where for a given $t=(t_1,\dots,t_n)$ the $(i,i+1)t$ stands for the $t$ where the elements $t_i$ and $t_{i+1}$ are exchanged, and }
& {\cal A}_i(t) = \frac{1}{ f_i(t;\beta) - f_{i+1}(t;\beta) }, \quad  {\cal B}_i(t) = \begin{cases} \frac{\left(f_i(t;\beta) - f_{i+1}(t;\beta)\right)^2 - 1}{\left(f_i(t;\beta) - f_{i+1}(t;\beta)\right)^2} & (t_i > t_{i+1}) , \\ \qquad  0  &  (t_i = t_{i+1}), \\
\qquad 1 & (t_i < t_{i+1}). \end{cases} \label{eq:KactionPhi2}
\end{align}

\end{prop}

\section{The Yangian}
With  the paper \cite{NT} as our  primary source we review several known facts about the Yangian of the Lie algebra $\gl_N.$ This is  a complex associative unital algebra $Y(\gl_N)$ with the set of generators $T_{ab}^{(s)}$ where $s=1,2,\dots\;$ and $a,b$ $=$ $1,\dots,N.$ The defining relations of the algebra  $Y(\gl_N)$ are usually written in terms of the formal power series in $u^{-1}$   
\begin{equation}
T_{ab}(u) = \delta_{ab}\cdot 1 + T_{ab}^{(1)}u^{-1} + T_{ab}^{(2)}u^{-1} + \dots ,
\end{equation}
and are 
\begin{equation}
(u-v)[T_{ab}(u),T_{cd}(v)]= T_{cb}(v)T_{ad}(u)- T_{cb}(u)T_{ad}(v).\label{eq:Ydef}
\end{equation}
Let $E_{ab}$ $\in$ $\End(\f^N)$ be the standard matrix units. Sometimes it is convenient to combine all the series $T_{ab}(u)$ into the single element
\begin{equation}
T(u) = \sum_{a,b=1}^N E_{ab}\otimes T_{ab}(u) \in \End(\f^N)\otimes Y(\gl_N)[[u^{-1}]]. \nonumber
\end{equation}
This element is invertible in the algebra $\End(\f^N)\otimes Y(\gl_N)[[u^{-1}]];$ denote 
\begin{equation}
\tilde{T}(u) = T(u)^{-1} = \sum_{a,b=1}^N E_{ab}\otimes \tilde{T}_{ab}(u).
\end{equation}
Let $h \in \f,$ and let $f(u)$ be a formal series from  $1+u^{-1}\f[[u^{-1}]].$\begin{prop}
Each of the maps 
\begin{eqnarray}
&\o_{f}: &T_{ab}(u) \mapsto f(u) T_{ab}(u),  \label{eq:omf}\\ 
&\xi(h): &T_{ab}(u) \mapsto T_{ab}(u + h),  \label{eq:xi}\\ 
&\s_{N}: &T_{ab}(u) \mapsto \tilde{T}_{ab}(-u), \label{eq:sigma}
\end{eqnarray}
determines an  automorphism of the algebra $Y(\gl_N).$
\end{prop}
Due to the defining relations (\ref{eq:Ydef}) every sequence of distinct indices ${\bold l}=(l_1,\dots,l_M)$ where $1 \leq l_a\leq M+N,$ determines an embedding of algebras 
\begin{equation}
\phi_{{\bold l}}: Y(\gl_M) \rightarrow Y(\gl_{M+N}) : T_{ab}(u) \mapsto T_{l_a l_b}(u).  \label{eq:embeddphi}
\end{equation}
The embedding $\phi_{{\bold l}}$ with ${\bold l} = (1,\dots,M)$ is called standard.

For $m=1,\dots,N$ define the series 
\begin{equation}
A_m(u) = \sum_{w\in S_m}\sgn(w) T_{1 w(1)}(u)T_{2 w(2)}(u-1)\dots T_{m w(m)}(u-m+1),
\end{equation}
and set $A_0(u)=1.$ The series $A_N(u)$ is called the {\em quantum determinant} of the Yangian. The following proposition is well-known. A detailed proof can be found in \cite{MNO}. 
\begin{prop} \label{p:quantumdet}
The coefficients at $u^{-1}, u^{-2},\dots $ of the series $A_N(u)$ are free generators for the centre of the algebra $Y(\gl_N).$ 
\end{prop}
Consider the ascending chain of algebras 
\begin{equation}
Y(\gl_1)\subset Y(\gl_2) \subset \dots \subset Y(\gl_N) \label{eq:chain}
\end{equation}
determined by the standard embeddings. Denote by $A(\gl_N)$ the subalgebra of $Y(\gl_N)$ generated by the centers of all algebras in  (\ref{eq:chain}); this subalgebra is clearly commutative. By Proposition \ref{p:quantumdet} the coefficients of the series $A_1(u),\dots,A_N(u)$ generate the subalgebra $A(\gl_N).$ It is proven by I.Cherednik \cite{Cherednik}, that the subalgebra  $A(\gl_N)$  is a maximal commutative. However we will not use this fact in what follows.

\subsection{Finite-dimensional irreducible modules of $Y(\gl_N)$}

Let $W$ be a finite-dimensional irreducible module of  $Y(\gl_N).$ A non-zero vector $w \in W$ is called {\em singular } provided it is annihilated by all coefficients of the series $T_{ab}(u)$ with $a>b.$ The following classification theorem is due to Drinfeld \cite{Drinfeld}.  
\begin{thm}
Every finite-dimensional irreducible $Y(\gl_N)$-module contains a unique up to a multiplier singular vector, say $w.$ The $w$ is an eigenvector of all coefficients of the series $A_1(u),\dots,A_N(u).$ Furthermore, one has   
\begin{equation}
\frac{A_{m+1}(u)A_{m-1}(u-1)}{A_{m}(u)A_{m}(u-1)}\cdot w = \frac{P_m(u-1)}{P_m(u)}\cdot w;  \quad m=1,\dots,N-1 \label{eq:DP1}
\end{equation}
for certain monic polynomials $P_1(u),\dots,P_{N-1}(u)$ with coefficients in $\f.$ Every collection of $N-1$ monic polynomials arises in this way. Up to an isomorphism of $Y(\gl_N)$-modules, modules with the same polynomials  $P_1(u),\dots,P_{N-1}(u)$ may differ only by an automorphism of the algebra $Y(\gl_N)$ of the form $\o_f.$
\end{thm}
The $N-1$ polynomials that characterize an irreducible  $Y(\gl_N)$-module are called the {\em Drinfeld polynomials} of this module. From (\ref{eq:DP1}) it follows that if $w$ is a singular vector, then we have   
\begin{equation}
\frac{T_{m+1,m+1}(u-m)}{T_{m,m}(u-m)}\cdot w = \frac{P_m(u-1)}{P_m(u)}\cdot w;  \quad m=1,\dots,N-1. \label{eq:DP2}
\end{equation}
Denote by ${\mbox{}}^{\s}W$ the irreducible $Y(\gl_N)$-module obtained from $W$ by the pullback through the automorphism $\s_N.$ Let ${\mbox{}}^{\s}P_1(u),\dots,{\mbox{}}^{\s}P_{N-1}(u)$ be the Drinfeld polynomials of the module ${\mbox{}}^{\s}W.$ 
\begin{prop} \label{p:sigmaDP}
For each $m=1,\dots,N-1$ we have ${\mbox{}}^{\s}P_m(u) = (-1)^{\deg P_m}P_m(-u+m-1).$
\end{prop}
\subsection{Realization of a certain class of  $Y(\gl_N)$-modules} \label{sec:tame}
An explicit construction of an arbitrary finite-dimensional irreducible $Y(\gl_N)$-module is an open problem. A realization of any {\em tame} module, which by definition is a  finite-dimensional irreducible $Y(\gl_N)$-module with a diagonalizable action of the subalgebra $A(\gl_N),$ was proposed in \cite{Cherednik} and \cite{NT}. Here we recall some elements of this construction.  

Let $E_{ab}$ with $a,b=1,\dots,N$ be generators of the Lie algebra $\gl_N.$ The universal enveloping algebra $U(\gl_N)$ is embedded into $Y(\gl_N)$ as a subalgebra by the map $E_{ba} \mapsto T_{ab}^{(1)}.$ Moreover there is a homomorphism   
\begin{equation}
\pi_N : Y(\gl_N) \rightarrow U(\gl_N) : T_{ab}(u) \mapsto \delta_{ab} + E_{ba}u^{-1}. \label{eq:evaluation}
\end{equation}
This homomorphism is called the {\em evaluation} homomorphism. It is an important tool in construction of Yangian modules. 

Let ${\bold l} = (M+1,\dots,M+N),$ and let $\psi= \s_{M+N}\phi_{{\bold l}} \s_{N}$ be a twisted version of the embedding (\ref{eq:embeddphi}). Thus the $\psi$ is a homomorphism from $Y(\gl_N)$ into $Y(\gl_{M+N}).$ The algebra $U(\gl_M)$ is realized inside $U(\gl_{M+N})$ by the standard embedding.

\begin{prop}[\cite{NT}] \label{p:pipsi}
The image of the homomorphism
\begin{equation}
\pi_{M+N}\psi : Y(\gl_N) \rightarrow  U(\gl_{M+N})
\end{equation}
commutes with the subalgebra $U(\gl_M)$ in $U(\gl_{M+N}).$
\end{prop}
Let $M\geq 0.$ Consider a pair of partitions $\l \supset \mu$ such that $\l = (\l_1,\dots,\l_{M+N})$ and $\mu = (\mu_1,\dots,\mu_M).$ Let $V_{\l}$ be the irreducible $\gl_{M+N}$-module with the highest weight $\l,$ and let $V_{\l,\mu}$ be the subspace of $V_{\l}$ spanned by all highest weight vectors of $\gl_M$ with the weight $\mu,$ where $U(\gl_M)$ is realized inside $U(\gl_{M+N})$ by the standard embedding. Due to Proposition \ref{p:pipsi} the  $V_{\l,\mu}$ is a $Y(\gl_N)$-module. Nazarov and Tarasov \cite{NT} proved that this module is tame, and determined its Drinfeld polynomials. To describe these Drinfeld polynomials it is convenient to represent the pair of partitions $\l$, $\mu$ by the skew Young diagram of the shape $\l/\mu.$ This diagram is the set 
\begin{equation}
\{ (i,j) \in \zint^2 | i \geq 1, \l_i \geq j > \mu_i \}.
\end{equation}
In this definition we do not distinguish between partitions that differ by zeroes. We employ the usual graphic representation of a skew Young diagram: a point $(i,j) \in \zint^2$ is represented by a unit square on the plane $\real^2$ with the centre $(i,j)$, the coordinates $i$ and $j$ increase downward and from left to right respectively. The {\em content} of a square $(i,j)$ is the difference $j-i.$

\begin{eg}
The diagram $\l / \mu$ with $\l = (5,5,3,3,1)$ and $\mu = (3,3,2,2)$ is 

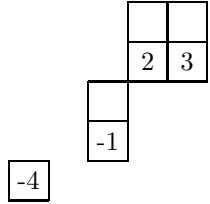
\begin{center}
\unitlength=0.75pt

\begin{picture}(120,110)
\put(60,100){\line(1,0){40}}
\put(60,80){\line(1,0){40}}
\put(40,60){\line(1,0){60}}
\put(40,40){\line(1,0){20}}
\put(40,20){\line(1,0){20}}
\put(0,20){\line(1,0){20}}
\put(0,0){\line(1,0){20}}
\put(0,0){\line(0,1){20}}
\put(20,0){\line(0,1){20}}
\put(40,20){\line(0,1){40}}
\put(60,20){\line(0,1){80}}
\put(80,60){\line(0,1){40}}
\put(100,60){\line(0,1){40}}
\put(80,60){\makebox(20,20){3}}
\put(60,60){\makebox(20,20){2}}
\put(40,20){\makebox(20,20){-1}}
\put(0,0){\makebox(20,20){-4}}
\end{picture}

\end{center}
where a number represents the content of the square in the bottom of each column. 
\end{eg}
\begin{prop}[\cite{NT}]\label{p:tameDP}
Let $\l/\mu$ be the skew Young diagram associated with the module $V_{\l,\mu}.$ Then the height of each column of the $\l/\mu$ is less or equal to $N.$ The Drinfeld polynomials of the  $V_{\l,\mu}$ are 
\begin{equation}
P_m(u) = \prod_{c}(u + c); \qquad m=1,\dots,N-1;
\end{equation}
where the product is taken over contents of bottom squares of  columns of the height $m$ in the skew Young diagram $\l/\mu.$
\end{prop}

For an irreducible $Y(\gl_N)$-module $W,$ and $h\in \f$ denote by $W(h)$ the pullback of $W$ trough the automorphism $\xi(h).$ Clearly the module $W(h)$ is irreducible and its Drinfeld polynomials are obtained from the Drinfeld polynomials of the $W$ by shifting all roots by $-h.$   

The Yangian is a Hopf algebra \cite{Drinfeld}. In particular it has the coproduct $\Delta: Y(\gl_N) \rightarrow Y(\gl_N)^{\otimes 2}$ defined by the map 
\begin{equation}
 T_{ab}(u) \mapsto \sum_{c=1}^N T_{ac}(u)\otimes T_{cb}(u),
\end{equation}
where the tensor product is taken over  $\f [[u^{-1}]].$

\begin{df}
 The irreducible modules $V_{(1)},V_{(1^2)},\dots,V_{(1^N)}$ are  called the  {\em fundamental } $Y(\gl_N)$-modules. The Drinfeld polynomials of the fundamental module $V_{(1^p)}$ are  
\begin{equation}
 P_m(u) = \begin{cases} 1 & \text{if $ m \neq p,$} \\  u + 1 - m &  \text{if $ m = p.$} \end{cases}
\end{equation}
\end{df}
\begin{rmk}
By abuse of terminology we will sometimes call the module $V_{(1^p)}(h),$ $(h\in\f)$ a fundamental module as well.
\end{rmk}

To finish our review of the Yangian algebra, we will quote the theorem, due to Nazarov and Tarasov, which gives classification of the tame modules, and an explicit construction for every tame module of $Y(\gl_N).$
\begin{thm}
Let $W$ be any finite-dimensional irreducible module of $Y(\gl_N).$ Let $P_m(u) = \prod_{i=1}^{\deg P_m}( u + z_{mi}),$ $(m=1,\dots,N-1)$ be the corresponding Drinfeld polynomials. 

The following three conditions are equivalent.

{\em (i)} For all $m\geq l$ one has $z_{li}-z_{mj} \neq 0,\dots,m-l$ unless $ (m,j) = (l,i).$

{\em (ii)} Up to some  automorphism $\o_f$ of $Y(\gl_N)$ the module $W$ is isomorphic to a
 tensor product 
\begin{equation}
   V_{\l^{(1)}/\mu^{(1)}}(h^{(1)})\otimes V_{\l^{(2)}/\mu^{(2)}}(h^{(2)})\otimes \dots \otimes V_{\l^{(n)}/\mu^{(n)}}(h^{(n)})
\end{equation}

for certain $ \l^{(s)} \supset \mu^{(s)},$ $h^{(s)} \in \f,$ $(s=1,\dots,n);$ such that $h^{(r)}-h^{(s)}\not\in \zint $ for all $r\neq s.$ 

{\em (iii)} The action of the subalgebra $A(\gl_N)$ on the module $W$ is diagonalizable. 
\end{thm}

\section{Eigenbases of the algebra $A(\gl_N)$ in tensor products of fundamental modules}\label{sec:fun}

In this section we will consider tensor products of fundamental modules of the form 
\begin{equation}
 V_{(1^{p_1})}(h^{(1)})\otimes V_{(1^{p_2})}(h^{(2)})\otimes \dots \otimes V_{(1^{p_r})}(h^{(r)}), \label{eq:fun}
\end{equation}
where $ p_s \in \setN$ and $h^{(s)} \in \real.$ We will see, that if the numbers $h^{(1)},\dots,h^{(r)}$ satisfy certain conditions, this module admits a unique up to normalization eigenbasis of the subalgebra $A(\gl_N),$ and the eigenvalue spectrum of this subalgebra is simple. Our main goal will be to deduce a certain triangularity property for the elements of this eigenbasis. 

Let $V=\f^N$ be the vector module of $\gl_N.$ Let $E_{ab} \in \End(V)$ be the matrix unit. It acts on the basis $v_1,\dots,v_N$ by $ E_{ab}v_c = \delta_{bc}v_a.$ Let $f\in \real$ and define the Yangian homomorphism 
\begin{equation}
\pi(f): Y(\gl_N) \rightarrow \End(V) : T_{ab}(u) \mapsto \delta_{ab} + E_{ba} (u + f)^{-1}. \label{eq:pihom}
\end{equation}
For a sequence of real numbers $f=(f_1,\dots,f_n)$ we denote by $\pi(f)$ the tensor product of the homomorphisms (\ref{eq:pihom}): $ \pi(f) = \pi(f_1)\otimes \dots \otimes \pi(f_n): Y(\gl_N)\pn \rightarrow \End(V\pn).$ The assignment $ T_{ab}(u) \mapsto \pi(f)\Delta^{(n)}T_{ab}(u)$ defines a $Y(\gl_N)$-module structure on $V\pn .$

With the numbers $p_1,\dots,p_r$ introduced in (\ref{eq:fun}) set $q_0 = 0$ and for each $s=1,\dots,r$ define the unique $q_s$ by $ p_s = q_s - q_{s-1}.$   

For a pair $(i,j)$ such that $0\leq i < j \leq n$ define the partial antisymmetrization operator $A_{(i,j) } \in \End(V\pn)$ by setting for each $a = (a_1,\dots,a_n)\in \setN^n$: 
\begin{multline}
A_{(i,j)}(v_{a_1}\otimes \cdots\otimes v_{a_n}) = \\ \sum_{w\in S_{j-i}}\sgn(w) v_{a_1}\otimes \cdots\otimes v_{a_i}\otimes v_{a_{i+w(1)}}\otimes v_{a_{i+w(2)}}\otimes \cdots \otimes v_{a_{i+w(j-i)}}\otimes v_{a_{j+1}}\otimes \cdots   \otimes v_{a_n}. 
\end{multline} 
Denote the sequence $(p_1,\dots,p_r)$ by $p$ and let $n=p_1+\dots+p_r.$ Let the subspace $(V\pn)_p $ be the image in  $V\pn$ of the operator 
\begin{equation}
 A_p = \prod_{s=1}^r A_{(q_{s-1},q_s)}.  
\end{equation}
The subspace $(V\pn)_p$ is obviously isomorphic to the tensor product (\ref{eq:fun}) as a $\gl_N$-module. 

Furthermore, suppose the sequence $f=(f_1,\dots,f_n)$ satisfies the condition
\begin{equation}
\text{ ${f}_i = {f}_{i+1}+1 $ when $ q_{s-1} < i < q_s $ for each $  s=1,2,\dots, r. $  }  \tag{C1}\label{eq:fbar1}
\end{equation}
For each $s=1,\dots,r$ define $h^{(s)} = f_{1+q_{s-1}}.$

\begin{prop} \label{p:Vpisomfun}
The coefficients of all the series $\pi(f)\Delta^{(n)}T_{ab}(u)$ leave the subspace $(V\pn)_p$ invariant. The $Y(\gl_N)$-modules $(V\pn)_p$ and {\em (\ref{eq:fun})} are isomorphic. 
\end{prop}
\begin{pf} The standard fusion procedure gives the required statement. 
\end{pf}
Define the  set ${\cal T}_p$ labelled by the sequence $p$ as follows    
\begin{multline}
{\cal T}_p =  \{ a = (a_1,\dots,a_n ) \in \setN^n \; | \; a_i < a_{i+1} \; \text{when}\; q_{s-1} < i < q_s \; \text{ for  } \; s=1,\dots, r \} \end{multline}
then the set 
\begin{equation}
\{ \varphi(a) = A_p v(a) \; | \; a \in {\cal T}_p \}  \label{eq:phi-def}
\end{equation}
is a basis of $ ( V\pn)_{p}$.

Let now the sequence $f$ satisfy the condition (\ref{eq:fbar1}) and the condition
\begin{equation}
 \qquad {f}_{q_s} > {f}_{q_s+1} + 1   \qquad \text{for} \quad   s=1,2,\dots, M-1.  \tag{C2} \label{eq:fbar2}
\end{equation}
The conditions (\ref{eq:fbar1}) and (\ref{eq:fbar2}) imply
\begin{equation}
h^{(s)}- h^{(s-1)} > p_s \quad \text{for each $s=1,\dots,r-1.$} \label{eq:h>h}
\end{equation}
A proof of the following proposition is contained in \cite{U2} Appendix A.
\begin{prop} \label{prop2} \mbox{}

{\em (i)} In the module $(V\pn)$ there is a unique up to normalization of eigenvectors eigenbasis of the commutative algebra $\pi(f)\Delta^{(n)}(A(\glN)).$ Denote this eigenbasis by $\{ \chi(a) \; | \; a \in {\cal T}_p \}.$

{\em (ii)} The expansion of this eigenbasis in the basis $\{ \varphi(a) \}$ is lower unitriangular. Which means that we have 
\begin{equation}
 \chi(a) = \varphi(a) + \sum_{b > a} c(a,b) \varphi(b)
\end{equation}
with certain real coefficients $c(a,b).$

{\em (iii)} The eigenvalues of the  algebra $\pi(f)\Delta^{(n)}A(\glN)$ are given by the formula
\begin{align}
& \pi(f) \Delta^{(n)}( A_m(u) ) \chi(a) = A_m(u;f;a) \chi(a), \qquad m=1,2,\dots,N; \label{p2iv}\\  
& \text{where} \quad  A_m(u;f;a)= \prod_{i=1}^n \frac{u+f_i + \delta(a_i \leq m) }{u+ f_i}. \nonumber  
\end{align}

{\em (iv)} The $N$-tuples of rational functions in the variable $u$: $  A_1(u;f;a), \dots,  A_N(u;f;a) $ are distinct for distinct $a \in {\cal T}_p$. In other words the spectrum of the  $\pi(f) \Delta^{(n)}( A_m(u) )$ on $(V\pn)_p$ is simple.
\end{prop}
In the expression for the eigenvalue $A_m(u;f;a)$ above we used the convention that for a statement $P,$ $\delta(P) = 1$ if $P$ is true, and $\delta(P) = 0$ otherwise. 

\section{The action of the Yangian on the space of states of the Spin Sutherland Model} \label{sec:YangactFn}
In this section we will define an action of the Yangian $Y(\gl_N)$ on the space of states $\F$ of the Spin Sutherland Model. This action is constructed by using the Drinfeld correspondence \cite{DrinfeldDual} between right modules of the degenerate affine Hecke algebra ${\bold H'}$ and left modules of the algebra   $Y(\gl_N).$ We will also find the decompositon of the   space of states $\F$ into invariant subspaces of this action.

\subsection{The Drinfeld correspondence}
Let $M$ be a right module of the degenerate affine Hecke algebra ${\bold H'}.$ For $a,b=1,\dots,N$ and $i=1,\dots,n$ define  $E_{ab}^{(i)} \in \End(V\pn)$ by  $E_{ab}^{(i)} = 1^{\otimes i-1}\otimes E_{ab} \otimes 1^{\otimes n-i}.$ Here the $E_{ab}$ is the matrix unit in $\End(V).$ Let $u$ be a formal variable, and define the $L$-operator
\begin{equation}
L_{ab}^{(i)}(u) = \delta_{ab}1 + ( u - \delta_i)^{-1}\otimes E_{ba}^{(i)} \in \End\left(M\otimes V\pn \right)[[u^{-1}]],
\end{equation}
where the denominator $( u - \delta_i)^{-1}$ is to be expanded into a series in $u^{-1}.$

The symmetric group $S_n$ acts from the left on the tensor product $V\pn$ by the permutation of  factors. The  $S_n$ also acts from the right on the module $M$ as the subalgebra of  ${\bold H'}.$

The following theorem is due to Drinfeld \cite{DrinfeldDual}.
\begin{thm}
The coefficients of the following series leave the space $ M\otimes_{S_n} V\pn $ invariant  
\begin{equation}
\widehat{T}_{ab}(u) = \sum_{c_1 = 1}^N \dots \sum_{c_{n-1} = 1}^N L_{a c_1}^{(1)}(u)L_{c_1 c_2}^{(2)}(u) \cdots L_{c_{n-1} b}^{(n)}(u).
\end{equation}
The assignment $ T_{ab}(u) \mapsto \widehat{T}_{ab}(u)$ defines a left $Y(\gl_N)$-module structure on the space  $ M\otimes_{S_n} V\pn .$  \label{t:DDual}
\end{thm}

Let us apply the Drinfeld correspondence to the ${\bold H'}$-module $\cz$ where the action of the degenerate affine Hecke algebra is given by the assignments (cf. (\ref{eq:d->-d})) 
\begin{equation}
\delta_i \mapsto -d_i(\beta), \qquad s_i \mapsto -K_{ii+1}.
\end{equation}
Note that we may regard the above action as either left or right due to the special structure of the defining relations of the algebra ${\bold H'}.$ Due to the definition (\ref{eq:K=-P}) of the space $\F$ we have $\cz\otimes_{S_n} V\pn  = \F.$ In view of  Theorem \ref{t:DDual} the left action of the Yangian $Y(\gl_N)$ on the space $\F$ is defined. We will denote this action by $Y(\gl_N;\beta),$ and the corresponding matrix of the generators by $T_{ab}(u;\beta).$

\subsection{Some properties of the Yangian action  $Y(\gl_N;\beta)$}
Let us denote by $A(\gl_N;\beta)$ the commutative subalgebra in $\End(\F)$ defined by the maximal commutative subalgebra $A(\gl_N)\subset Y(\gl_N)$ and the action $ Y(\gl_N) \mapsto Y(\gl_N;\beta).$ In particular we have 
\begin{equation}
A_N(u;\beta) = \prod_{i=1}^n\frac{u + 1 + d_i(\beta)}{u +  d_i(\beta)} \label{eq:qdet}
\end{equation}
and hence the Hamiltonian $H_{\beta,N}$ (\ref{eq:H}) is an element in $A(\gl_N;\beta),$ and commutes with $Y(\gl_N;\beta).$ In other words the  $Y(\gl_N;\beta)$ is the Yangian symmetry of the Spin Sutherland Model \cite{BGHP}.
 
For an operator $B$ acting on $\F$ let $B^*$ be its adjoint relative to the scalar product ${\sprod{\cdot}{\cdot} }_{\beta,N}$ and for a series $B(u)$ with operator-valued coefficients we will use $B(u)^*$ to denote the series whose coefficients are adjoints of coefficients of the series $B(u)$. 

The following proposition is proven in  \cite{TU}.
\begin{prop} \label{p:sa}
 For all $ a ,b = 1,\dots,N$ we have $T_{ab}(u;\beta)^* = T_{ba}(u;\beta).$ Moreover for all $m=1,\dots,N$ we have  $A_m(u;\beta)^* = A_m(u;\beta).$ 
\end{prop}

\subsection{The decomposition of the space of states into irreducible Yangian submodules} 
Our aim now is to describe the decomposition of the space of states $\F$ into irreducible Yangian submodules with respect to the action $Y(\gl_N;\beta).$ Let $\bs=(s_1,\dots,s_n)$ be an element of the set $\LC_n^{(N)}$ (\ref{eq:LC}), and let $r$ be the number of distinct elements in sequence $\bs.$ With this  $\bs$ we associate: the sequence of integers $q_0=0,q_1,\dots,q_{r-1},q_r=n$ uniquely defined by the condition $s_{q_j} >  s_{q_j+1}$ $( j=1,2,\dots,r-1);$ and the sequence of integers  $p(\bs)=(p_1,\dots,p_r)$  defined by $ p_j = q_j - q_{j-1}$ $( j=1,2,\dots,r).$ Compare these definitions with those made in Section \ref{sec:fun}. Clearly we have $ 1\leq p_j \leq N$ and $p_1+\dots+p_r =n.$

Further, let  $ f(\bs) = (f_1,f_2,\dots,f_n)$ be defined by $ f_i = f_i(\bs;\beta),$ where the $f_i(\bs;\beta)$ is the eigenvalue of the Cherednik-Dunkl operator (\ref{eq:d-eigenvalue}) on the non-symmetric Jack polynomial $E_{\bs}^{(\beta)}.$ Note that the sequence  $ f(\bs)$ satisfies the condition (\ref{eq:fbar1}) of Section  \ref{sec:fun}. In view of this, and Proposition \ref{p:Vpisomfun}, the space $(V\pn)_{p(\bs)}$ is the $Y(\gl_N)$-module with the action given by $\pi(f(\bs))\Delta^{(n)}T_{ab}(u).$ The following theorem is proven in \cite{TU}

\begin{thm} \label{t:isom}
Consider the space of states $\F$ as the Yangian module with the action $Y(\gl_N;\beta).$ 
We have \begin{equation}\F = \oplus_{\bs \in \LC_n^{(N)}} F_{\bs},\label{eq:Yangdecomp} \end{equation} where for each $\bs  \in \LC_n^{(N)}$ the subspace $F_{\bs}$ is invariant with respect to $Y(\gl_N;\beta)$ and is isomorphic to $(V\pn)_{p(\bs)}$ as a $Y(\gl_N)$-module.
\end{thm}

The  isomorphism between $(V\pn)_{p(\bs)}$ and $F_{\bs}$ is explicitly given by the operator $U(\bs;\beta):$ $(V\pn)_{p(\bs)}$ $\rightarrow$ $ F_{\bs}$ which is defined for any $v\in ( V\pn)_{p(\bs)}$ as follows   
\begin{equation}
U(\bs;\beta)v = \sum_{t \sim s} E_t^{(\beta)}(z)\otimes R_t^{(\beta)}v,   \label{eq:Uop}
\end{equation}
where the sum is taken over all distinct rearrangements $t$ of $s$, and $ R_t^{(\beta)}$ $ \in$ $\End(V\pn)$ is defined by the recursive realtions:  
\begin{gather}
 R_{\bs}^{(\beta)} = 1, \\ 
 R_{(i,i+1)t}^{(\beta)}= -\check{R}_{ii+1}( f_i(t;\beta) - f_{i+1}(t;\beta))R_{t}^{(\beta)} \quad \text{for} \quad t_i > t_{i+1}. 
\end{gather}
Here $\check{R}_{ii+1}(u)$ $=$ $u^{-1} + P_{ii+1}$ and $(i,i+1)t$ denotes the element of $\zint^n$ obtained by interchanging $t_i$ and $t_{i+1}$ in $t$ $=$ $(t_1,\dots,t_n)$.  

With notations of Section \ref{sec:fun}, for each  $a \in {\cal T}_{p(\bs)}$ define 
\begin{eqnarray} 
\Phi^{(\beta)}(\bs,a) & = & U(\bs;\beta)\varphi(a), \\ 
 X^{(\beta)}(\bs,a) & = & U(\bs;\beta)\chi(a).
\end{eqnarray}
Observe now that, since  $\beta $ is a positive real number, the condition (\ref{eq:fbar2}) of Section \ref{sec:fun} is satisfied by the sequence $ f(\bs)$. Proposition \ref{prop2} and Theorem \ref{t:isom} imply the following

\begin{prop} \label{prop3} \mbox{}

{\em (i)} $\{ X^{(\beta)}(\bs,a) \; | \; a \in {\cal T}_{p(\bs)} \}$ is the unique up to normalization of eigenvectors eigenbasis of the commutative algebra $ A(\glN;\beta)$ in the module   $F_{\bs}.$

{\em (ii)} The expansion of this eigenbasis in the basis $\{ \Phi^{(\beta)}(\bs,a) \}$ is lower unitriangular. Which means that we have 
\begin{equation}
 X^{(\beta)}(\bs,a) = \Phi^{(\beta)}(\bs,a) + \sum_{b > a} c(a,b) \Phi^{(\beta)}(\bs,b)\label{p3ii}
\end{equation}
with certain real coefficients $c(a,b).$

{\em (iii)} The eigenvalues of the  algebra $A(\glN;\beta)$ are given by the formula
\begin{align}
& A_m(u;\beta) X^{(\beta)}(\bs,a) = A_m(u;f(\bs);a)X^{(\beta)}(\bs,a), \qquad m=1,2,\dots,N; \label{p3iii}\\  
& \text{where} \quad  A_m(u;f(\bs);a)= \prod_{i=1}^n \frac{u+f_i(\bs;\beta) + \delta(a_i \leq m) }{u+ f_i(\bs;\beta)}. \nonumber  
\end{align}

{\em (iv)} The $N$-tuples of rational functions in the variable $u$: $  A_1(u;f(\bs);a), \dots,  A_N(u;f(\bs);a) $ are distinct for distinct $a \in {\cal T}_{p(\bs)}$. 
\end{prop}

\subsection{The eigenbasis of the commtative algebra $A(\gl_N;\beta)$}
We will now describe  properties of the eigenbasis $\{ X^{(\beta)}(\bs,a) \; | \;\bs \in \LC_n^{(N)}, a \in {\cal T}_{p(\bs)} \}$ in some detail. Our main goal is to establish the unitriangularity of the expansion of the elements of this basis in the basis of normally ordered wedges. 

For any pair $(\bs = (s_1,\dots,s_n) \in \LC_n^{(N)}, a=(a_1,\dots,a_n)\in {\cal T}_{p(\bs)})$ and each $i=1,\dots,n$ define $k_i = a_{n-i+1} - N s_{n-i+1}.$ It is easy to see that the sequence $k=(k_1,\dots,k_n)$ is an element of the set $\LC_n^{(1)}$ (cf. (\ref{eq:LC})). Moreover we have (cf. Section \ref{sec:preliminary}) $\un{k} = (a_{n},a_{n-1},\dots,a_1)$ and $\ov{k}=(s_{n},s_{n-1},\dots,s_1).$ The described correspondence between the set of all pairs $(\bs  \in \LC_n^{(N)}, a\in {\cal T}_{p(\bs)})$ and the set  $\LC_n^{(1)}$ is one-to-one, therefore we may label the elements of the bases $\{ \Phi^{(\beta)}(\bs,a)\}$ and  $\{ X^{(\beta)}(\bs,a)\}$ by sequences from $\LC_n^{(1)}.$ Let us define 
\begin{equation} 
X_k^{(\beta,N)} = (-1)^{\frac{n(n-1)}{2}}X^{(\beta)}(\bs,a).
\end{equation}
Recall that in Proposition \ref{prop1} we  introduced the  eigenbasis $\{\Psi_k^{(\beta,N)}| k \in\LC_n^{(1)}\} $ of the Hamiltonian $H_{\beta,N}.$ 
\begin{lemma} \label{l:Phi=Psi}
We have the equality $(-1)^{\frac{n(n-1)}{2}}\Phi^{(\beta)}(\bs,a) = \Psi_k^{(\beta,N)}.$
\end{lemma}
\begin{pf}
Using the expression (\ref{eq:Uop}) for the operator $U(\bs;\beta)$ we have
\begin{equation}
 \Phi^{(\beta)}(\bs,a) = \sum_{t \sim \bs} E_t^{(\beta)}(z)\otimes R_t^{(\beta)}\varphi(a), \label{eq:l20}
\end{equation}
where $t \sim \bs$ means that $t$ belongs to the set of all distinct rearrangements of the sequence $\bs.$
In view of the triangularity (\ref{eq:nsJtri}) of the non-symmetric Jack polynomial  $E_t^{(\beta)}(z)$ we may split the $E_t^{(\beta)}(z)$ with $t \sim \bs$ as follows
\begin{gather}
E_t^{(\beta)}(z) = E_t^{(\beta)}(z)' + E_t^{(\beta)}(z)'', \\
\text{where} \quad  E_t^{(\beta)}(z)' = \sum_{ r \preceq t, \; r \sim \bs} e^{(\beta)}_{tr} z^r , \quad \text{ and } \quad  E_t^{(\beta)}(z)'' = \sum_{ r^+ < \bs } e^{(\beta)}_{tr} z^r. \label{eq:l21}
\end{gather}
The vector $\Phi^{(\beta)}(\bs,a)$ belongs to the subspace  $F_{N,n}$ $\subset$  $\cz\otimes(V\pn).$ Therefore both of the vectors
\begin{equation}
\sum_{t \sim \bs} E_t^{(\beta)}(z)'\otimes R_t^{(\beta)}\varphi(a) \quad \text{and}\quad  \sum_{t \sim \bs} E_t^{(\beta)}(z)''\otimes R_t^{(\beta)}\varphi(a)
\end{equation}
also belong to this  subspace since  monomials $z^r$ which appear in the decomposition of the $E_t^{(\beta)}(z)'$ as a vector in $\cz\otimes(V\pn)$ are distinct from monomials $z^r$ which appear in the decomposition of the $E_t^{(\beta)}(z)''$  as implied by  (\ref{eq:l21}).

Taking into account the triangularity of the non-symmetric Jack polynomial  $E_t^{(\beta)}(z),$ and the equality  $ R_{\bs}^{(\beta)} = 1,$ we may write
\begin{equation}
\sum_{t \sim \bs} E_t^{(\beta)}(z)'\otimes R_t^{(\beta)}\varphi(a) = z^{\bs} \otimes \varphi(a) + \sum_{t \sim \bs,\; t \prec \bs} z^t \otimes \bar{\varphi}_t(a), \quad ( \bar{\varphi}_t(a) \in  V\pn).
\end{equation}
Since the expression above is a vector in $F_{N,n},$  by using the definition of the vector $\varphi(a)$ given in (\ref{eq:phi-def}) we obtain 
\begin{equation}
\sum_{t \sim s} E_t^{(\beta)}(z)'\otimes R_t^{(\beta)}\varphi(a) = A_n( z^s \otimes v(a) )  = (-1)^{\frac{n(n-1)}{2}} \hat{u}_k  \label{eq:e1prime}
\end{equation}
where $A_n$ (\ref{eq:antisymm}) is the operator of the total antisymmetrization in the space $\cz\otimes ( V\pn).$ 

The expansion of  the vector 
\begin{equation}
\sum_{t \sim s} E_t^{(\beta)}(z)''\otimes R_t^{(\beta)}\varphi(a) \label{eq:e2prime}
\end{equation}
in $\cz\otimes(V\pn)$ contains only monomials $z^r$ such that $ r^{+} < \bs.$ Therefore the expansion of the (\ref{eq:e2prime}) in the basis of normally ordered wedges $\{ \hat{u}_l \: | \: l \in \LC_n^{(1)} \}$ contains only  $\hat{u}_l$ such that $ \ov{l} > \ov{k} .$ Taking this, and (\ref{eq:e1prime}) into account we have
\begin{equation}
(-1)^{\frac{n(n-1)}{2}}\Phi^{(\beta)}(\bs,a) = \hat{u}_k  + \sum_{l \in \LC_n^{(1)},\; \ov{l} > \ov{k}} \varphi_{kl}^{(\beta)} \hat{u}_l,  \qquad ( \varphi_{kl}^{(\beta)} \in \ra). \label{eq:l22}
\end{equation}
The vector $\Phi^{(\beta)}(\bs,a)$ is an eigenvector of the quantum determinant $A_N(u;\beta)$ as implied by the equations (\ref{eq:d-eigenvalue}), (\ref{eq:qdet}) and (\ref{eq:l20}). Hence it is an eigenvector of the Hamiltonian $H_{\beta,N}.$ 

However according to Proposition \ref{prop1}(i) an eigenvector of the Hamiltonian $H_{\beta,N}$  with the expansion (\ref{eq:l22}) in the basis $\{ \hat{u}_l \: | \: l \in \LC_n^{(1)} \}$  is unique and equals to $\Psi_k^{(\beta,N)}.$ This proves the lemma.

\end{pf}

Theorem \ref{t:isom} and Proposition \ref{prop3} imply that the set $\{ X_k^{(\beta,N)}\: | \: k \in \LC_n^{(1)}\}$ is the unique up to normalization eigenbasis of the commutative algebra $A(\gl_N;\beta)$ in the space $\F.$ 

\begin{thm} \label{t:Xbasis} \mbox{} 

{ \em (i)} The transition matrix between the bases $\{ X_k^{(\beta,N)}\: | \: k \in \LC_n^{(1)}\}$ and  $\{ \hat{u}_k \: | \: k \in \LC_n^{(1)}\}$ is upper unitriangular. Which means that we have
\begin{equation}
 X_k^{(\beta,N)} = \hat{u}_k  + \sum_{ l < k } x^{(\beta,N)}_{kl} \;\hat{u}_l \label{p4i}  
\end{equation}
with certain real coefficients $x^{(\beta,N)}_{kl}.$

{\em (ii)} For each $m=1,\dots,N$ we have 
\begin{align}
&  A_m(u;\beta)  X_k^{(\beta,N)} = A_m(u;\beta;k) X_k^{(\beta,N)}   \label{p4ii}\\  
& \text{where} \quad A_m(u;\beta;k) = \prod_{i=1}^n \frac{u+\beta^{-1}\ov{k_i}+i-1 + \delta(\un{k_i} \leq m) }{u+\beta^{-1}\ov{k_i}+i-1}. \nonumber  
\end{align}

{\em (iii)} The $N$-tuples of rational functions in $u$:  $  A_1(u;\beta;k), \dots,  A_N(u;\beta;k) $  are distinct for distinct $k \in \LC_n^{(1)}$. In other words the spectrum of the commutative algebra $A(\gl_N;\beta)$ is simple.

{\em (iv)} $ {\sprod{X_k^{(\beta,N)}}{X_l^{(\beta,N)}}}_{\beta,N} = 0 $ if $  k \neq l. $
\end{thm}

\begin{pf}
(i) Lemma \ref{l:Phi=Psi} and Proposition \ref{prop3}(ii) give 
\begin{equation}
X_k^{(\beta,N)} = \Psi_{k}^{(\beta,N)} +   \sum_{l \in \LC_n^{(1)},\; \ov{l} = \ov{k},\; \un{l} < \un{k} } c_{kl}\: \Psi_{l}^{(\beta,N)} \qquad ( c_{kl} \in \real ).  \label{eq:X-Psi}
\end{equation}
Observe that $ \ov{l} = \ov{k} $ and $\un{l} < \un{k}$ imply  ${l} < {k}.$ The transition matrix between the bases $\{ \Psi_{k}^{(\beta,N)} \: | \: k \in \LC_n^{(1)}\}$ and $\{ \hat{u}_k \: | \: k \in \LC_n^{(1)}\}$ is upper unitriangular by Proposition \ref{prop1}(ii). Hence (\ref{eq:X-Psi}) leads to (i).

(ii) Follows immediately from Proposition \ref{prop3}(iii) and the explicit expressions for the eigenvalues of the Cherednik-Dunkl operators given by (\ref{eq:d-eigenvalue}).

(iii) Proposition \ref{prop3}(iv) shows that the sets of the eigenvalues  $ A_1(u;\beta;k), \dots,  A_N(u;\beta;k) $ are distinct for distinct $k$ which have equal $\ov{k}.$ In other words eigenvalues of the commutative algebra $A(\gl_N;\beta)$ separate eigenvectors that belong to the same irreducible component of the Yangian action $Y(\gl_N;\beta).$ We will prove that eigenvalues of the quantum determinant $A_N(u;\beta)$ separate between these  irreducible components. Eigenvalue $A_N(u;\beta;k)$ depends only on $\ov{k}$ and is a rational function of the form $P(u+1)/P(u)$ where $P(u)$ is a monic polynomial. Roots of this polynomial are $ 1- i - \beta^{-1}\ov{k_i},$ they are pairwise distinct due to the assumption that $\beta$ is a positive real number. Hence the polynomial $P(u)$ determines the sequence $(\ov{k_1},\dots,\ov{k_n})$ uniquely. This proves (iii).

(iv) Follows from (iii) and Proposition \ref{p:sa}.
\end{pf}

The $Y(\gl_N)$-action $Y(\gl_N;\beta)$  gives rise to an  action of the Lie algebra $\gl_N \subset Y(\gl_N)$ on the space of states $\F.$ The generators of this action are $T_{ab}^{(1)}(\beta)$ $(a,b =1,\dots,N).$  
As a corollary to Theorem \ref{t:Xbasis} we obtain 
\begin{prop} \label{p:wtn}\mbox{} \\
{\em (i)} Elements of the basis  $\{ X_k^{(\beta,N)} \; | \; k \in \LC_n^{(1)} \}$ are weight vectors of $\gl_N.$ For each $a=1,\dots,N$ we have  
\begin{equation}
T_{aa}^{(1)}(\beta)X_k^{(\beta,N)} = \sum_{i=1}^n \delta(\un{k_i} = a) X_k^{(\beta,N)}. \nonumber 
\end{equation}

\noindent {\em (ii) }Elements of the basis  $\{ X_k^{(\beta,N)} \; | \; k \in \LC_n^{(1)} \}$ are eigenvectors of the degree operator $D = z_1\frac{\p}{\p z_1} + \cdots + z_n\frac{\p}{\p z_n}.$ We have 
\begin{equation}
D X_k^{(\beta,N)} = |\ov{k}|  X_k^{(\beta,N)}. \nonumber 
\end{equation}
\end{prop}

Let $\mm = (m_1,\dots,m_n)$ be a non-decreasing sequence of integers such that each element $m_i$ appears in $\mm$ with multiplicity less or equal to $N.$ We denote the set of all such sequences by $M_n^{(N)},$ and represent an element $\mm$   of $M_n^{(N)}$ as $((r_1)^{p_1}(r_2)^{p_2}\dots (r_l)^{p_l})$ where $r_1 < r_2 < \dots < r_l$ and $p_i$ $( 1\leq p_i \leq N)$ denotes  the  multiplicity of $r_i$ in $\mm.$ If $k$ is an element of $\LC_n^{(1)}$ then $\ov{k}$ is an element of $M_n^{(N)}.$ Conversely for each $\mm \in M_n^{(N)}$ there is at least one $k \in \LC_n^{(1)}$ such that $ \ov{k} = \mm.$ 

\begin{prop} \label{p:decn}
For each $\mm = ((r_1)^{p_1} (r_2)^{p_2} \dots (r_l)^{p_l})\in M_n^{(N)}$ the linear space $ F_{N,n}(\mm) = \oplus_{\{ k | \ov{k} = \mm \}} \f X_k^{(\beta,N)} $ is invariant and irreducible with respect to the Yangian action $Y(\gl_N;\beta).$ Moreover  $ F_{N,n}(\mm)$  is isomorphic as a $Y(\gl_N)$-module to the tensor product 
\begin{equation}
V_{(1^{p_1})}( a_{1})\otimes V_{(1^{p_2})}( a_{2}) \otimes \cdots \otimes V_{(1^{p_l})}( a_{l})
\end{equation}
where $ a_{s} = \beta^{-1}r_s -1 + p_1 + \dots + p_s$ $(s=1,\dots,l).$
\end{prop}
\begin{pf}
In notations of Theorem \ref{t:isom} $F_{N,n}(\mm) = F_{\bs}$ where $\bs = ((r_l)^{p_l} (r_{l-1})^{p_{l-1}}\dots  (r_1)^{p_1}).$ Theorem \ref{t:isom} and Proposition \ref{p:Vpisomfun} imply the required statements. 

\end{pf}


\section{An isomorphism between the space of states of the Spin Sutherland Model and the space of symmetric Laurent polynomials}
In this section we will introduce an  isomorphism between the space of states of the Spin Sutherland Model and the space of symmetric Laurent polynomials. This isomorphism may be regarded as a finite-particle version of the well-known fermion-boson correspondence in the representation theory of infinite-dimensional Kac-Moody algebras \cite{KR}. We will determine the image of the eigenbasis of the commutative algebra $A(\gl_N;\beta)$ under this isomorphism. This image is a basis of the space of symmetric Laurent polynomials. A subset of this basis gives a basis of the space of symmetric polynomials which coincides with a degeneration of the basis of Macdonald polynomials \cite{MacBook}. This degeneration is described as follows. A Macdonald polynomial depends on two parameters $q$ and $t.$ Let $\o_N$ be an $N$th primitive root of unity, and let $p$ be a parameter. Set $q=\o_N p, \; t=\o_N p^{N\beta + 1}$ and take the limit $ p \rightarrow 1.$ In this limit a Macdonald polynomia!
!
l degenerates into a polynomial 
which we call a Jack($\gl_N$) polynomial. When  $N=1$ this polynomial is just the usual Jack polynomial. The Jack($\gl_N$) polynomials were recently used to compute certain dynamical correlation functions in the Spin Sutherland Model \cite{U2}. We, however, will not consider this subject in the present article.  
 
\subsection{Scalar products on the space of Laurent polynomials} Let $L_n =\cx$ be the algebra of Laurent polynomials in variables $x_1,\dots,x_n.$ For each  $f$ $=$ $f(x_1,\dots,x_n)$ $\in$ $L_n$ let $ \overline{f(x_1,\dots,x_n)}$ be the Laurent polynomial with complex conjugated coefficients. And let $f^* $ $ =$ $  \overline{f(x_1^{-1},\dots,x_n^{-1})}.$ As in Section \ref{sec:preliminary} we use  $[f]_1$ to denote the constant term in $f$. Throughout this section we let $b$ denote an arbitrary positive integer number.

Now let us introduce 
\begin{align}
& {\nabla}(\bb,N) = \prod_{1\leq i \neq j \leq n} ( 1 - x_i^N x_j^{-N})^{\bb} \label{eq:nabla}\\
\intertext{and}
& {\Delta}(\bb,N) = \prod_{1\leq i \neq j \leq n} ( 1 - x_i^N x_j^{-N})^{\bb}( 1 - x_i x_j^{-1} ) = {\nabla}(\bb,N)\prod_{1\leq i \neq j \leq n}( 1 - x_i x_j^{-1} ).
\end{align}
Both ${\nabla}(\bb,N)$ and ${\Delta}(\bb,N)$ are symmetric Laurent polynomials since $\bb$ is a positive integer.

We define two scalar products on the linear space $L_n$ by  
\begin{align}
{\sproddd{f}{g}}_{\bb,N}& = \frac{1}{n!}[{\nabla}(\bb,N)f^* g]_1 \label{eq:scAL} \\
\intertext{and}
{\sprodd{f}{g}}_{\bb,N} & = \frac{1}{n!}[{\Delta}(\bb,N)f^* g]_1. \label{eq:spbN}
\end{align}
Using ${\sproddd{f}{g}}_{0}$ as a shorthand notation for ${\sproddd{f}{g}}_{0,N}$ we have 
\begin{equation}
{\sproddd{f}{g}}_{\bb,N} = {\sproddd{f}{{\nabla}(\bb,N) g}}_0, \label{eq:e1}
\end{equation}
where we regard ${\nabla}(\bb,N)$ as a multiplication operator on $L_n.$

Let $A_n^{\pm}$ be the subspace of skew-symmetric Laurent polynomials in $L_n$. For each $l =(l_1,\dots,l_n)$ $\in$ $\zint^n$ define the antisymmetric monomial $\hat{a}_l$  as follows: 
\begin{equation}
\hat{a}_l = x^{l_1}\wedge x^{l_2} \wedge \cdots \wedge x^{l_n} = \sum_{w \in S_n}\sgn(w)  x^{l_1}_{w(1)} x^{l_2}_{w(2)}\cdots  x^{l_n}_{w(n)}.
\end{equation}
The set  $\{\hat{a}_l\:|\: l\in \LC_n^{(1)} \}$ is a basis of  $A_n^{\pm}.$

Let $\delta = (n-1,\dots,0).$ The monomial $\hat{a}_{\delta}$ is equal to the Vandermonde determinant $ \prod_{i<j}(x_i - x_j).$ Let $\Lambda_n^{\pm}$ be the subspace of symmetric Laurent polynomials in $L_n$. For each $l =(l_1,\dots,l_n)$ $\in$ $\zint^n$ define the symmetric Laurent polynomial $s_l$ by 
\begin{equation}
s_l = \frac{\hat{a}_{l + \delta}}{\hat{a}_{\delta}}.
\end{equation}
The set  $\{s_l\:|\: l\in \LC_n^{(n)} \}$ is a basis of  $\Lambda_n^{\pm}.$ A sequence $l\in 
\LC_n^{(n)} $ which contains only non-negative  elements is obviously identified with a partition $\l$ of length less or equal to $n.$ For any such $l$ the $s_{l}$ is a symmetric polynomial equal to the Schur polynomial $s_{\l}.$ Schur polynomials $s_{\l}$ where $\l$ runs through all partitions of length less or equal to  $n$ form a basis in the space of symmetric polynomials $\Lambda_n.$ 

\subsection{An isomorphism between the space of states of the Spin Sutherland Model and the space of symmetric Laurent polynomials}\label{sec:isomOmega} Now we  define an isomorphism which maps the space of states $\F$ into the space of symmetric Laurent polynomials  $\Lambda_n^{\pm}.$ Let us fix the sequence $o = (0,-1,\dots,-n+1).$ We will call this $ o$ the {\em vacuum sequence }. Recall the definition of the vacuum wedges in (\ref{eq:vac}). We have the equalities 
\begin{equation}
\hat{u}_{o} =  X_{o}^{(\beta,N)} = \vac(0).
\end{equation}

For any  $k \in \LC_n^{(1)}$ the difference $ k - {o}$ is an element of  $\LC_n^{(n)}.$ Define the isomorphism of linear spaces  
\begin{equation}
\Omega : \F \rightarrow \Lambda_n^{\pm} :  \hat{u}_k \mapsto s_{ k - {o}}. \label{eq:Om}
\end{equation}
The isomorphism $\Omega$ is a composition of the two isomorphisms $\o_*$ and $\o_{**}$ defined by   
\begin{equation}
\o_*:  \F\rightarrow A_n^{\pm} : \hat{u}_k \mapsto \hat{a}_k;\quad \o_{**}:   A_n^{\pm}\rightarrow \Lambda_n^{\pm} : \hat{a}_k \mapsto s_{ k - {o}}; 
\end{equation}
so that  $ \Omega =\o_{**} \o_*.$ Note that for any $\hat{a} \in A_n^{\pm}$ we have 
\begin{equation}
\o_{**}(\hat{a}) = (x_1 x_2 \cdots x_n)^{n-1} \frac{\hat{a}}{ \hat{a}_{\delta}}.
\end{equation}

\begin{lemma} \label{l:Pf=}
Let $P = P(z_1,\dots,z_n)$ be a symmetric Laurent polynomial in the variables $z_1,\dots,z_n.$ Regard $P$ as a multiplicative operator on the space $\F.$ For any $f \in \F$ we have  
\begin{equation}
\o_*\left( P f \right) = P(x_1^{-N},\dots,x_n^{-N})\o_*\left(  f \right). \label{eq:Pf=}
\end{equation}
\end{lemma}
\begin{pf}
The algebra of symmetric  Laurent polynomials in the variables $z_1,\dots,z_n$ is generated by the elements\begin{equation}
p_r = p_r(z_1,\dots,z_n) = z_1^r + z_2^r + \cdots + z_n^r\qquad (r=0,\pm 1,\pm 2, \dots \;).
\end{equation}
Therefore it is sufficient to prove the equality (\ref{eq:Pf=}) for  $P =p_r.$ Since wedges $\hat{u}_k$ span the space $\F$ we may assume $f = \hat{u}_k.$ We have 
\begin{gather}
 \o_*(p_r \hat{u}_k ) = \o_*\left( \sum_{i=1}^n u_{k_1}\wedge \cdots \wedge u_{k_i - Nr}\wedge \cdots \wedge u_{k_n} \right) = \\ 
 = \sum_{i=1}^n x^{k_1}\wedge \cdots \wedge x^{k_i - Nr}\wedge \cdots \wedge x^{k_n} = p_r(x_1^{-N},\dots,x_n^{-N}) \hat{a}_k. \nonumber 
\end{gather}
This proves the required statement.
\end{pf}
Recall the definition of the scalar product ${\sprod{\cdot}{\cdot}}_{\beta,N}$ given in Section \ref{sec:scalar}. This definition implies that for any $f,g \in \F$ and $\bb\in \nat$ we have the identity 
\begin{equation}
{\sprod{f}{g}}_{\bb,N} = {\sprod{f}{{\Delta}(z;\bb)g}}_{0,N}, \label{eq:e2}
\end{equation}
where we regard the symmetric Laurent polynomial 
\begin{equation}
{\Delta}(z;\bb) = \prod_{1\leq i \neq j \leq n} ( 1 - z_i z_j^{-1})^{\bb}
\end{equation}
as a multiplication operator on the space $\F.$

The basis of normally ordered wedges $\{\hat{u}_k\:|\:k\in\LC_n^{(1)}\}$ is orthonormal relative to 
the scalar product ${\sprod{\cdot}{\cdot}}_{0,N}.$ Likewise the basis $\{\hat{a}_k\:|\:k\in\LC_n^{(1)}\}$ is orthonormal relative to  the scalar product ${\sproddd{\cdot}{\cdot}}_{0}.$ The definition of the isomorphism $\o_*$ now implies that  for  any $f,g \in \F$ we have 
\begin{equation}
 {\sprod{f}{g}}_{0,N} = {\sproddd{\o_*(f)}{\o_*(g)}}_0.  \label{eq:e3}
\end{equation}

We are ready now to formulate one of the  important  technical statements of this section.
\begin{lemma} \label{l:spF=spA}
The isomorphism $\o_*$ is an isometry. Which means that for any $f,g$ $\in$ $\F$ we have 
\begin{equation}
 {\sprod{f}{g}}_{\bb,N} = {\sproddd{\o_*(f)}{\o_*(g)}}_{\bb,N}.\label{eq:s=s}
\end{equation}
\end{lemma}
\begin{pf} Observe that ${\Delta}(z;\bb)$ is a symmetric Laurent polynomial. The following chain of identities gives the required result:
\begin{align}
& {\sproddd{\o_*(f)}{\o_*(g)}}_{\bb,N} =\text{(by \ref{eq:e1})}={\sproddd{\o_*(f)}{{\nabla}(\bb,N) \o_*(g)}}_0 = \text{(by Lemma \ref{l:Pf=})} = \nonumber \\ 
& = {\sproddd{\o_*(f)}{\o_*({\Delta}(z;\bb) g)}}_0 =\text{(by \ref{eq:e3})}={ \sprod{f}{{\Delta}(z;\bb) g}}_{0,N}=\text{(by \ref{eq:e2})}={\sprod{f}{g}}_{\bb,N}. \nonumber 
\end{align} 
\end{pf}

The following proposition is the main point of this section.

\begin{prop} \label{p:spF=spL}
The isomorphism $\Omega$ is an isometry. Which means that for any $f,g$ $\in$ $\F$ we have 
\begin{equation}
 {\sprod{f}{g}}_{\bb,N} = {\sprodd{\Omega(f)}{\Omega(g)}}_{\bb,N}.
\end{equation}
\end{prop}
\begin{pf}
Definitions of the scalar products ${\sprodd{\cdot}{\cdot}}_{\bb,N}$ and   ${\sproddd{\cdot}{\cdot}}_{\bb,N}$ along with the definition of the isomorphism $\o_{**}$ immediately lead to 
\begin{equation}
{\sprodd{\o_{**}(f)}{\o_{**}(g)}}_{\bb,N} = {\sproddd{f}{g}}_{\bb,N}, \label{eq:s=s2}
\end{equation}
for any $f,g \in A^{\pm}_n.$ The required statement now follows from Lemma \ref{l:spF=spA}.
\end{pf}

\subsection{Macdonald polynomials} 
Let $\l = (\l_1,\l_2,\dots \;)$ be a partition. As in Section \ref{sec:tame} we represent $\l$ by its diagram. For example below is the diagram for the  parition $\lambda = (6,4,4,3,1).$
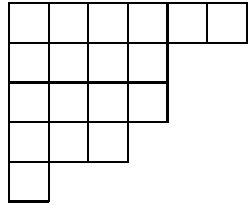
\begin{center}

\unitlength=0.75pt
\begin{picture}(120,110)
\put(0,100){\line(1,0){120}}
\put(0,80){\line(1,0){120}}
\put(0,60){\line(1,0){80}}
\put(0,40){\line(1,0){80}}
\put(0,20){\line(1,0){60}}
\put(0,0){\line(1,0){20}}
\put(0,0){\line(0,1){100}}
\put(20,0){\line(0,1){100}}
\put(40,100){\line(0,-1){80}}
\put(60,100){\line(0,-1){80}}
\put(80,100){\line(0,-1){60}}
\put(100,100){\line(0,-1){20}}
\put(120,100){\line(0,-1){20}}
\end{picture}  \\
\end{center}
For a square $s \in \lambda$  arm-length $ a_{\lambda}(s)$, leg-length $ l_{\lambda}(s)$, arm-colength $ a'(s)$ and leg-colength  $ l'(s)$ are defined as the number of squares in the diagram of $\lambda$ to the east, south, west and north from $s$ respectively. The  content of a square $s$ is $c(s) =  a'(s) - l'(s)$ and the  hook-length is $h_{\lambda}(s) =  a_{\lambda}(s) + l_{\lambda}(s) + 1.$ All partitions which we will encounter in this section have length less or equal to $n.$

For a partition $\lambda$ a  tableau $T$ of shape $\lambda$ \cite{MacBook} is a sequence of partitions:
\begin{equation}
\emptyset = \lambda^{(0)} \subset   \lambda^{(1)} \subset  \dots \subset  \lambda^{(r)} = \lambda   \label{eq:TT}
\end{equation}
such that each skew diagram $\theta^{(i)} = \lambda^{(i)} / \lambda^{(i-1)}$ $(1 \leq i \leq r)$ is a horizontal strip. The sequence $ ( |\theta^{(1)}|,\dots,|\theta^{(r)}|)$ is called the weight of $T.$

The symmetric monomial associated with partition $\l$ is an element of $\Lambda_n$ defined by  
\begin{equation} 
m_{\l} = \sum x_1^{\alpha_1} x_2^{\alpha_2}\cdots x_n^{\alpha_n}
\end{equation}
with the sum taken over all distinct rearrangements $(\alpha_1,\dots,\alpha_n)$ of $\l.$ Symmetric monomials  $m_{\l}$ where $\l$ runs through all partitions of length less or equal to $n$ form a basis of  $\Lambda_n.$ The transition matrix between this basis and the basis of Schur polynomials is upper unitriangular. That is to say $s_{\l}$ has the expansion of the form 
\begin{equation}
s_{\l} = m_{\l} + \sum_{\mu < \l} K_{\l \mu} m_{\mu} \label{eq:s=m+}
\end{equation}
with suitable coefficients $K_{\l \mu}$ called Kostka numbers.

Let $q$ and $t$ be  parameters. 
For each partition of length $\leq n$ the Macdonald polynomial $P_{\l}(q,t)$ is an element of $\Lambda_n \otimes_{\f}\f(q,t) $ defined by the following expansion in the basis of symmetric monomials $m_{\l}$  
\begin{equation}
P_{\lambda}(q,t) = m_{\lambda} + \sum_{\mu < \lambda} u_{\lambda\mu}(q,t) m_{\mu} \label{eq:MP=m+}
\end{equation}
where the summation is over partitions of length  $\leq n$ and the  coefficients $u_{\lambda\mu}(q,t)$  are  given by 
\begin{equation} 
u_{\lambda\mu}(q,t) = \sum_{T} \psi_{T}(q,t) 
\end{equation}
summed over all tableaux of shape $\lambda$ and weight $\mu.$ To describe the quantity $\psi_{T}(q,t),$ for partitions $\lambda$ and $\mu$ such that $\mu \subset \lambda$ and $\l/\mu$ is a horizontal strip, let $R_{\mu}$ denote the unique (possibly empty) row of $\mu$ which intersects $\l/\mu.$ Then   
\begin{equation}
\psi_{T}(q,t) = \prod_{i=1}^r \psi_{\lambda^{(i)}/\lambda^{(i-1)}}(q,t),
\end{equation}
where  
\begin{align}
& \psi_{\lambda/\mu}(q,t) = \prod_{s \in R_{\mu}} \frac{b_{\mu}(s;q,t)}{b_{\lambda}(s;q,t)}, \label{eq:psi(q,t)}\\ 
& \text{and} \quad b_{\lambda}(s;q,t) = \frac{1 - q^{a_{\lambda}(s)}t^{l_{\lambda}(s) + 1}}{1 - q^{a_{\lambda}(s)+ 1}t^{l_{\lambda}(s)}}.\label{eq:b(q,t)}
\end{align}

Symmetric monomials and Schur polynomials are obtained as specializations of $P_{\l}(q,t).$ We have 
\begin{equation}
             P_{\l}(q,1)  = m_{\l}  \qquad \text{and} \qquad   P_{\l}(q,q)  = s_{\l}. \nonumber
 \end{equation}
Another well-known specialization of $P_{\l}(q,t)$ is the Jack polynomial $J_{\l}^{(\alpha)}$ defined for a positive real $\alpha$ by  
\begin{equation}
J_{\l}^{(\alpha)} = \lim_{q \rightarrow 1} P_{\l}(q,q^{\frac{1}{\alpha}}).
\end{equation}

Let $k$ be a non-negative integer number. Introduce the symmetric Laurent polynomial  
\begin{equation}
\Box(q,q^k) = \prod_{1\leq i \neq j \leq n} \prod_{r=0}^{k-1} ( 1 - q^r x_i x_j^{-1}). \label{eq:MPwf}
\end{equation}
And define on  $\Lambda_n \otimes_{\f} \f(q) $ a scalar product by 
\begin{equation}
{\sprodd{f}{g}}_{q,q^k} = \frac{1}{n!}[\Box(q,q^k)f^* g]_1 . \label{eq:MPsp}
\end{equation}
Polynomials  $P_{\lambda}(q,q^k)$ are known to be pairwise orthogonal relative to this  scalar product.

\subsection{Jack($\gl_N$) polynomials} \label{sec:JP}
Here we will consider a certain specialization of Macdonald polynomials, other than the specializations we have mentioned above. Let $\o_N$ be an $N$th primitive root of unity. That is to say  $(\o_N)^N = 1$ and $(\o_N)^i \neq 1$ for all $i=1,\dots,N-1.$  
\begin{df}  \label{def:JP}
For any positive real $\gamma$ the limit of the Macdonald polynomial $P_{\lambda}(q,t)$ defined by 
\begin{equation}
P_{\lambda}^{(\gamma,N)} = \lim_{p\rightarrow 1} P_{\l}(\o_N p,\o_N p^{\gamma})\label{eq:limit}
\end{equation}
is called the {\em Jack($\gl_N$)} polynomial.
\end{df}
Note that if we take $N=1$ in this definition we get the usual  Jack polynomial $J_{\lambda}^{(\frac{1}{\gamma})}.$ 

\noindent Taking the limit (\ref{eq:limit}) in the expansion (\ref{eq:MP=m+}) we see that
\begin{equation}
P_{\lambda}^{(\gamma,N)} = m_{\lambda} + \sum_{\mu < \lambda} u_{\lambda\mu}^{(\gamma,N)} m_{\mu} \label{eq:JP=m+}
\end{equation}
where the sum runs over partitions of length $\leq n,$ and the coefficients $u_{\lambda\mu}^{(\gamma,N)}$ are given by 
\begin{equation} 
u_{\lambda\mu}^{(\gamma,N)} = \sum_{T} \psi_{T}^{(\gamma,N)} \label{eq:uLM}
\end{equation}
summed over all tableaux of shape $\lambda$ and weight $\mu.$  In the same notations as in  (\ref{eq:TT}) and (\ref{eq:psi(q,t)}) we now have 
\begin{equation}
\psi_{T}^{(\gamma,N)}  = \prod_{i=1}^r \psi_{\lambda^{(i)}/\lambda^{(i-1)}}^{(\gamma,N)},
\end{equation}
where  
\begin{align}
& \psi_{\lambda/\mu}^{(\gamma,N)} = \prod_{s \in R_{\mu}} \frac{b_{\mu}^{(\gamma,N)}(s)}{b_{\l}^{(\gamma,N)}(s)}, \label{eq:psi(gamma,N)} \\ 
& \text{and} \quad b_{\l}^{(\gamma,N)}(s) = \begin{cases} \frac{ a_{\lambda}(s) + \gamma l_{\lambda}(s) + \gamma }{ a_{\lambda}(s) + \gamma l_{\lambda}(s) + 1 }  & \text{ if $ h_{\lambda}(s) \equiv 0\bmod N,$} \\ 1 & \text{otherwise.} \end{cases}  \label{eq:b}
\end{align}
These formulas show that each $u_{\lambda\mu}^{(\gamma,N)}$ is a rational function of $\gamma$ with positive rational coefficients. In particular,  $u_{\lambda\mu}^{(\gamma,N)}$ is positive for all positive real $\gamma.$

Now let $\bb$ be a positive integer.
Consider the behaviour of the scalar product (\ref{eq:MPsp}) in the limit (\ref{eq:limit}) wherein we put $\gamma =N\bb+1.$  We have 
\begin{equation}
\lim_{p\rightarrow 1}\Box(\o_N p, \o_N p^{Nb+1} ) =  \prod_{1\leq i \neq j \leq n} \prod_{r=0}^{N\bb}(1 - \o_N^r x_i x_j^{-1}) = \prod_{1\leq i \neq j \leq n} (1 -  x_i^N x_j^{-N})^{\bb}  (1 -  x_i x_j^{-1}).
\end{equation}
The right-hand side of this equation is exactly the weight function $\Delta(\bb,N)$ of the scalar product ${\sprodd{\cdot}{\cdot}}_{\bb,N}$ which was introduced in (\ref{eq:spbN}). Thus for $\gamma = N\bb + 1 $ the scalar product for Macdonald polynomials ${\sprodd{\cdot}{\cdot}}_{q,q^{\gamma}}$ degenerates in the limit (\ref{eq:limit}) into the scalar product ${\sprodd{\cdot}{\cdot}}_{\bb,N}$ defined in (\ref{eq:spbN}). 

The orthogonality of Macdonald polynomials entails that for all non-negative integer $b$ the Jack($\gl_N$) polynomials  $P_{\lambda}^{(N\bb +1,N)}$ are pairwise orthogonal relative to the scalar product $(\ref{eq:spbN}).$

Since the basis of symmetric monomials and the basis of Schur polynomials are related by an upper unitriangular matrix (\ref{eq:s=m+}), the expansion (\ref{eq:JP=m+}) gives  
\begin{equation}
P_{\lambda}^{(\gamma,N)} = s_{\lambda} + \sum_{\mu < \lambda} v_{\lambda\mu}^{(\gamma,N)} s_{\mu} \label{eq:JP=s+}
\end{equation}
where the coefficients $v_{\lambda\mu}^{(\gamma,N)}$ are rational functions of $\gamma.$ Let $\bb \in \nat.$ The upper unitriangularity of either the expansion (\ref{eq:JP=s+}) or the expansion (\ref{eq:JP=m+}), and the condition 
\begin{equation}
{\sprodd{P_{\lambda}^{(N\bb+1,N)}}{P_{\mu}^{(N\bb+1,N)}}}_{\bb,N} = 0 \quad \text{ if $ \l \neq \mu$ }
\end{equation}
define the polynomials $P_{\lambda}^{(N\bb+1,N)}$ uniquely by Gram-Schmidt orthogonalization.

\subsection{Jack($\gl_N$) polynomials as eigenvectors of the maximal commutative subalgebra of the Yangian action} Recall the definition of the vacuum sequence $o$ given in Section \ref{sec:isomOmega}. Suppose  $k \in \LC_n^{(1)}.$ If  $ k \supseteq o,$ then $\ov{k} \subseteq \ov{o}.$ Moreover, if $n\equiv 0 \bmod N,$ then $\ov{k} \subseteq \ov{o}$ implies  $ k \supseteq o.$ Let $\F^{0}$ be the linear subspace of the space of states $\F$ defined by  \begin{equation} 
\F^{0} = \oplus_{k : k \supseteq o} \f \hat{u}_k.   \label{eq:F0n}
\end{equation}
If $n\equiv 0 \bmod N,$ then an equivalent definition of this subspace is 
\begin{equation} 
\F^{0} = \oplus_{k : \ov{k} \subseteq \ov{o} } \f \hat{u}_k. 
\end{equation}
The condition  $k \supseteq o$ is equivalent to the inequality $k_n \geq o_n.$ Therefore  $k \supseteq o$ and  $k \geq l$ imply   $l \supseteq o.$ Due to Theorem \ref{t:Xbasis}(i) we have $ X_{k}^{(\beta,N)} \in \F^{0}$ provided  $k \supseteq o.$ Moreover, again by Theorem \ref{t:Xbasis}(i), the set  $ \{X_{k}^{(\beta,N)} \: | \: k \in \LC_n^{(1)}, k \supseteq o\}$ is a basis of $\F^{0}.$ The vectors  $X_{k}^{(\beta,N)}$ such that $\ov{k} \subseteq \ov{o}$ form the same basis of $\F^{0}$ provided  $n\equiv 0 \bmod N.$  Since  $X_{k}^{(\beta,N)}$ with fixed $\ov{k}$ form a basis of an invariant subspace of the Yangian action $Y(\gl_N,\beta),$ we see that $\F^{0}$ is invariant with respect to this action when $n\equiv 0 \bmod N.$ 

Consider now the image of the subspace  $\F^{0}$ under the action of the isomorphism $\Omega$ (\ref{eq:Om}). The conditions $k\in \LC_n^{(1)}$ and  $k \supseteq o$ imply that $k - o$ is a partition, say $\l$, of length less or equal to $n.$ Therefore $\Omega$ maps $\F^{0}$  isomorphically into the space of symmetric polynomials $\Lambda_n.$ By definition of this isomorphism and in view of Theorem \ref{t:Xbasis}(i) we have  
\begin{equation}
\Omega\left(X_{o+\l}^{(\beta,N)}\right) = s_{\l} + \sum_{\mu < \l} x_{\l \mu}^{(\beta,N)}\: s_{\mu} 
\label{eq:X=s}\end{equation}
where $x_{\l \mu}^{(\beta,N)}$ are certain rational functions of $\beta.$
\begin{thm} \label{t:X=P}
For each partition $\l$ of length less or equal to $n$ we have 
\begin{equation}
\Omega\left(X_{o+\l}^{(\beta,N)}\right) = P_{\l}^{(N\beta + 1,N)}.
\end{equation}
\end{thm}
\begin{pf}
Theorem \ref{t:Xbasis}(iv) and Proposition \ref{p:spF=spL} give 
\begin{equation}
{\sprodd{\Omega\left(X_{o+\l}^{(\bb,N)}\right)}{\Omega\left(X_{o+\mu}^{(\bb,N)}\right)}}_{\bb,N} = 0 \quad \text{if $ \l \neq \mu$ }
\end{equation}
for all $\bb \in \nat.$ This and the upper unitriangularity (\ref{eq:X=s}) define the basis of symmetric polynomials $\Omega\left(X_{o+\l}^{(\bb,N)}\right)$ uniquely by Gram-Schmidt orthogonalization. The Jack($\gl_N$) polynomials $P_{\l}^{(N\bb + 1,N)}$ satisfy the same two conditions. Therefore  the statement of the theorem follows for all $\beta \in \nat.$ 

The coefficients $x_{\l \mu}^{(\beta,N)}$ of the expansion (\ref{eq:X=s}) are rational functions in $\beta,$ and so are the  coefficients $v_{\l \mu}^{(N\beta+1,N)}$ of the expansion (\ref{eq:JP=s+}). Since we have 
\begin{equation}
x_{\l \mu}^{(\beta,N)} = v_{\l \mu}^{(N\beta+1,N)} 
\end{equation}
for all non-negative integer $\beta,$ these coefficients are equal as rational functions. Hence  the statement of the theorem. 
\end{pf}

\section{Fock space representation of the affine Kac-Moody algebra  $\sll$} \label{bigsec:fockspace}
In this section we review the Fock space representation of the affine Kac-Moody algebra  $\sll.$ This review is to be considered as a preparation to Section \ref{sec:YangactF} where we  define a Yangian action on the Fock  space and describe how it is decomposed into irreducible Yangian submodules.

\subsection{Fock space representation of $\sll$}

\subsubsection{The affine Lie algebra $\sll$}
Let $N\geq 2$ and let ${\frak h}$ be a $N+1$-dimensional vector space over $\ra$ with basis $\{h_0,h_1,\dots,h_{N-1},D\}.$ Let $\{ \Lambda_0,\Lambda_1,\dots,\Lambda_{N-1},\delta\}$ be the corresponding dual basis of the dual space  ${\frak h}^*.$ That is 
\begin{equation}
\langle \Lambda_i , h_j \rangle  = \delta_{ij},\quad \langle \Lambda_i , D \rangle = 0, \quad \langle \delta , h_i \rangle = 0, \quad \langle \delta , D \rangle = 1. \nonumber 
\end{equation}
It will be convenient to extend the index set so that $\Lambda_i = \Lambda_{(i\pmod N)}$ for all integer $i.$ For all $i \in \zint$ we set $ \alpha_i = 2\Lambda_i - \Lambda_{i+1}-\Lambda_{i-1} + \delta_{i0}\delta,$ where $ \delta_{ij} = 1$ if $i-j \equiv 0\bmod N$ and $\delta_{ij} =0$ otherwise. 

The $N \times N$ matrix $ \| \langle \alpha_i, h_j \rangle \| $ is the generalized Cartan matrix of type $A_{N-1}^{(1)}.$ The associated affine Kac-Moody algebra is denoted by $\sll.$ It is defined as the Lie algebra generated by $D$ and $e_i,f_i,h_i$ with $0\leq i < N$ subject to the relations  
\begin{gather*}
[h_i,h_j]=0, \quad [h_i,D]=0, \\ \mbox{}
[h_i,e_j]=\langle \alpha_j, h_i \rangle e_j, \quad [D,e_j]=  \delta_{j0}e_j, \\ \mbox{} 
[h_i,f_j]=-\langle \alpha_j, h_i \rangle f_j, \quad [D,f_j]=  -\delta_{j0}f_j, \\ \mbox{} 
[e_i,f_j] = \delta_{ij}h_i, \\ \mbox{}
({\mathrm {ad}} e_i)^{1-\langle \alpha_j, h_i \rangle} e_j =0 \quad (i\neq j), \\ \mbox{} 
({\mathrm {ad}} f_i)^{1-\langle \alpha_j, h_i \rangle} f_j =0 \quad (i\neq j).
\end{gather*}
The centre of the algebra  $\sll$ is one-dimensional. It is generated by  the element $c = h_0+h_1+ \dots + h_{N-1}.$ The $\Lambda_i$ are known as the fundamental weights and the $\alpha_i$ as the simple roots of  $\sll.$ 

\subsubsection{Fock space representation of $\sll$} \label{sec:Fockrep} The space of states $\F$ admits an action of the algebra $\sll$ determined by the assignments    
\begin{gather*}
e_i \mapsto e_i^{(n)} = \sum_{j=1}^n z_j^{\delta_{i0}} \otimes E_{i, i+1}^{(j)}, \quad 
f_i \mapsto f_i^{(n)} = \sum_{j=1}^n z_j^{-\delta_{i0}} \otimes E_{i+1, i}^{(j)}, \\ 
h_i \mapsto h_i^{(n)} = \sum_{j=1}^n 1 \otimes E_{i,i}^{(j)} - 1 \otimes E_{i+1,i+1}^{(j)}, \quad  
D \mapsto D^{(n)} =  \sum_{j=1}^n D_j \otimes 1,
\end{gather*}
where the indices $i$ are considered modulo $N.$

Let $\zint^{\infty}$ be the set of all semi-infinite sequences $k = (k_1,k_2,k_3,\dots)$ with integer elements  $ k_i$ such that $ k_i = -i+1$ for all but finite number of $i\in \nat.$ The set $\zint^{\infty}$ contains the distinguished sequence $o = (0,-1,-2,\dots).$ We will call it the vacuum sequence.  

The Fock space $F$ is defined as a $\f$-linear  space  generated by the set of semi-infinite expressions    
\begin{equation}
\{ \hat{u}_k =  u_{k_1}\wedge u_{k_2} \wedge u_{k_3} \wedge \dots \; | \; k = (k_1,k_2,k_3,\dots) \in \zint^{\infty} \}, 
\end{equation}
modulo the relations $ u_{k_i}\wedge u_{k_{i+1}} = -  u_{k_i}\wedge u_{k_{i+1}}$ for each pair of indices $i$ and $i+1.$ An element $\hat{u}_k \in F$ will be  called a semi-infinite wedge. Let $\LC_{\infty}^{(1)}$ be the subset of  $\zint^{\infty}$ which consists of strictly decreasing sequences. Then the set  $ \{ \hat{u}_k \: | \:  k \in   \LC_{\infty}^{(1)} \}$  is a basis of  $F.$ An element of this basis is called a normally ordered wedge. 

Let us reserve the notation $ | m \rangle $ for the formal expression $ u_{m}\wedge u_{m-1} \wedge u_{m-2} \wedge \dots \,.$ There is an obvious homomorphism of linear spaces    $\F \rightarrow F$ determined by the assignment   
\begin{equation}
u_{k_1}\wedge u_{k_2} \wedge \dots \wedge u_{k_n} \mapsto u_{k_1}\wedge u_{k_2} \wedge \dots \wedge u_{k_n} \wedge | -n \rangle .
\end{equation}
Conversely, for each $v \in F$ there is a large enough $n$ and $v^{(n)} \in \F$ such that  
\begin{equation}
 v = v^{(n)} \wedge | -n \rangle .
\end{equation}

For all integer $m$ define the action of $e_i,f_i,h_i$ $(0\leq i < N-1)$ on  $ | m \rangle $ by 
\begin{align}
& e_i | m \rangle  = 0 , \\ 
& f_i | m \rangle = \begin{cases} u_{m+1}\wedge | m-1 \rangle & \text{ if $ i\equiv m \bmod N;$} \\
                                  0  &  \text{otherwise,} \end{cases} \\ 
& h_i | m \rangle = \begin{cases}  | m \rangle & \text{ if $ i\equiv m \bmod N;$} \\
                                  0  &  \text{otherwise.} \end{cases}
\end{align}
Let $x$ be any of the generators  $e_i,f_i,h_i$ $(0\leq i < N-1).$ The action of $x$ on  the Fock space is then defined by 
\begin{equation}
x v  =  (x^{(n)} v^{(n)})\wedge | -n \rangle +   v^{(n)} \wedge x | -n \rangle .
\end{equation}

The Fock space $F$ has two gradings. Let $\hat{u}_k $ be a normally ordered wedge. The homogeneous grading is defined by 
\begin{equation}
\deg_h(\hat{u}_k ) = | \ov{o} - \ov{k} | = \sum_{i\geq 1} \ov{o_i} - \ov{k_i}.
\end{equation}
And the principal grading is defined by 
\begin{equation}
\deg_p(\hat{u}_k ) = | k - o | = \sum_{i\geq 1}  {k_i} -{o_i}.
\end{equation}
The action of the degree generator on $F$ is then defined as $ D v = - \deg_h(v) v$ on any vector $v$ of homogeneous degree $\deg_h(v).$

\subsection{Constructions of the Fock space by inverse limits}  
\subsubsection{Principal construction} Let the  sequence $(k_1,k_2,\dots,k_n)$ be an element of the set $\LC_n^{(1)}$ (\ref{eq:LC}). In this section for any such sequence we will use the shorthand notation $k_{[n]}$ to indicate explicitly that it contains $n$ elements. If $k_{[m]} \in  \LC_m^{(1)}$ and $m\geq n,$ we will use  $k_{[n]}$ to denote the element of $\LC_n^{(1)}$ obtained by deleting last $m-n$ elements in  $k_{[m]}.$

Recall that the subspace $\F^{0} \subset \F$ (\ref{eq:F0n}) is defined by 
\begin{equation}
\F^{0} = \oplus_{k_{[n]} \supseteq o_{[n]}}\f \hat{u}_{k_{[n]}},
\end{equation}
where $ o_{[n]} = (0,-1,\dots,-n+1)$ is the vacuum sequence.

For any $k_{[n]} \in  \LC_n^{(1)}$ such that $k_{[n]} \supseteq o_{[n]},$ the difference $k_{[n]} - o_{[n]}$ is a partition of length less or equal to $n.$ We define the principal grading on  $\F^{0}$ by setting the {\em principal degree} of $\hat{u}_{k_{[n]}}$ to be equal  $|k_{[n]} - o_{[n]}|.$ For each $s=0,1,\dots, $ let $\F^{0,\langle s \rangle }$ be the homogeneous component of $\F^0$ of  principal degree  $s.$     

For $m \geq n$ define the projection 
\begin{equation}
\pi_{m,n}^{\langle s \rangle } : F_{N,m}^{0,\langle s \rangle } \rightarrow F_{N,n}^{0,\langle s \rangle } : \hat{u}_{k_{[m]}} \mapsto \begin{cases} \hat{u}_{k_{[n]}} & \text{if $ l(k_{[m]} - o_{[m]}) \leq n;$ }\\ 0 & \text{otherwise.}  \end{cases}\label{eq:princproj}
\end{equation}
This projection is a surjective homomorphism of linear spaces for all $m\geq n.$ Moreover if $m\geq n\geq s$ it is a bijection.

Let us  now form the inverse limit  
\begin{equation}
F_N^{0,\langle s \rangle } = \lim\begin{Sb} \longleftarrow \\ n \end{Sb} F_{N,n}^{0,\langle s \rangle }
\end{equation}
of the linear spaces $ F_{N,n}^{0,\langle s \rangle }$ relative to the homomorphisms $\pi_{m,n}^{\langle s \rangle }.$ By definition a vector in the linear space $ F_N^{0,\langle s \rangle }$ is a sequence $ ( f_n )_{n \geq 0} $ such that $f_n \in   F_{N,n}^{0,\langle s \rangle }$ and $\pi_{m,n}^{\langle s \rangle }(f_m) = f_n$ for all $m \geq n.$ Since  $\pi_{m,n}^{\langle s \rangle }$ is an isomorphism for  $m\geq n\geq s,$ the definition (\ref{eq:princproj}) implies 
\begin{equation}
 f_m = f_n\wedge u_{-n}\wedge u_{-n-1} \wedge \dots \wedge u_{-m+1}  \label{eq:princ-tail} 
\end{equation}provided $m\geq n\geq s.$ 
Let $F^{\langle s \rangle} $ be the component of the Fock space $F$ of  principal degree $s.$ From (\ref{eq:princ-tail}) it follows  that the map  
\begin{equation}
\Upsilon_p^{\langle s \rangle }: F_N^{0,\langle s \rangle } \rightarrow F^{\langle s \rangle}:( f_n )_{n \geq 0} \mapsto f_m \wedge | -m \rangle   \label{eq:princ-fock}
\end{equation}
does not depend on $m$ as soon as $m\geq s,$ and for any such $m$ is an isomorphism of linear spaces.

\subsubsection{Homogeneous construction} \label{sec:homconstr}
Now consider the linear space $F_{N,rN}^{0}$ for some $r\geq 0.$ On this space we introduce the homogeneous  grading by setting the {\em homogeneous degree}  of $\hat{u}_{k_{[n]}}$ to be equal  $|  \ov{o_{[rN]}} - \ov{k_{[rN]}}|.$ For each $d=0,1,\dots $ let $F_{N,rN}^{0,(d)}$ be the homogeneous component of  $F_{N,rN}^{0}$ of homogeneous  degree  $d.$

Similarly to the principal case, for $l \geq r$ we define the projection 
\begin{equation}
\pi_{l,r}: F_{N,l N}^{0} \rightarrow F_{N,rN}^{0 } : \hat{u}_{k_{[l N]}} \mapsto \begin{cases} \hat{u}_{k_{[rN]}} & \text{if $ l(k_{[lN]} - o_{[lN]}) \leq rN ;$ }\\ 0 & \text{otherwise.}  \end{cases}\label{eq:homproj}
\end{equation}
This projection obviously preserves the homogeneous degree. Let $\pi_{l,r}^{(d)}: F_{N,lN}^{0,(d)} \rightarrow F_{N,rN}^{0,(d)}$   be the  restriction of  $\pi_{l,r}$ onto the component $F_{N,lN}^{0,(d)}.$ The projection  $\pi_{l,r}^{(d)}$ is a surjective homomorphism of linear spaces for all $l\geq r.$ Moreover if $l\geq r\geq d$ it is a bijection because $ l(k_{[lN]} - o_{[lN]}) \leq dN $ if the homogeneous degree of $ \hat{u}_{k_{[l N]}} $ is less or equal to $d.$

Again we form the inverse limit  
\begin{equation}
F_N^{0,(d) } = \lim\begin{Sb} \longleftarrow \\ r \end{Sb} F_{N,rN}^{0,(d)}, \label{eq:hom}
\end{equation}
now of the linear spaces $ F_{N,rN}^{0,(d)}$ relative to the homomorphisms $\pi_{l,r}^{(d)}.$ By definition a vector in the linear space $ F_N^{(d)}$ is a sequence $ ( g_r )_{r \geq 0} $ such that $g_r \in   F_{N,rN}^{0,(d) }$ and $\pi_{l,r}^{(d) }(g_l) = g_r$ for all $l \geq r.$ Since  $\pi_{l,r}^{(d)}$ is an isomorphism for  $l\geq r\geq d,$ the definition (\ref{eq:homproj}) implies 
\begin{equation}
 g_l = g_r\wedge u_{-Nr}\wedge u_{-Nr-1} \wedge \dots \wedge u_{-Nl+1} \label{eq:hom-tail}
\end{equation} provided $l\geq r\geq d.$ 

Let $F^{(d)} $ be the component of the Fock space $F$ of  homogeneous  degree $d.$ From (\ref{eq:hom-tail}) it follows that the map  
\begin{equation}
\Upsilon_h^{(d)}: F_N^{0,(d) } \rightarrow F^{(d)}:( g_r )_{r \geq 0} \mapsto g_r \wedge | -Nr \rangle   
\label{eq:isomhomfock}\end{equation}
does not depend on $r$ if  $r\geq d,$ and for any such $r$ is an isomorphism of linear spaces.

\begin{prop} \label{p:hom-vs-princ}
Let $s\geq 0,$ and let $(f_n)_{n\geq 0}$ be a vector in $F_N^{0,\langle s \rangle }.$ Suppose that for all $r\geq 0$ we have $f_{rN} \in F_{N,rN}^{0,(d )}$ for a certain $d\geq 0.$ Then the sequence $(f_{rN} )_{r\geq 0}$ is a vector in $F_N^{0,(d) },$ and moreover we have 
\begin{equation}
\Upsilon_h^{(d)}\left( (f_{rN} )_{r\geq 0} \right) = \Upsilon_p^{(s)}\left( (f_n )_{n\geq 0} \right).\label{eq:hom-fock}
\end{equation}
\end{prop}
\begin{pf}
From  the definition of the inverse limit $F_N^{0,\langle s \rangle }$ it follows that  
\begin{equation}
\pi_{l,r}^{(d)} (f_{lN}) = f_{rN} 
\end{equation}
for all $l \geq r.$ Hence $(f_{rN} )_{r\geq 0}$ is a vector in $F_N^{0,(d) }.$  

Let us choose $r$ to be large enough, so that both of the inequalities $ r \geq d$ and $ rN \geq s $ hold. 
Then by (\ref{eq:princ-fock}) and (\ref{eq:isomhomfock}) we obtain 
\begin{equation}
\Upsilon_h^{(d)}\left( (f_{rN} )_{r\geq 0} \right)  = \Upsilon_p^{(s)}\left( (f_{n} )_{n\geq 0} \right) = f_{rN}\wedge |-rN \rangle   
\end{equation}
which completes the proof.
\end{pf}

\subsection{Symmetric functions}
Let $\Lambda_n$ be the linear space of symmetric polynomials in variables $x_1,\dots,x_n.$ This space has the basis which consists of all Schur polynomials $s_{\l}(x_1,\dots,x_n)$ labelled by partitions $\l$ of length less or equal to $n.$ The space  $\Lambda_n$ is graded. We have  
\begin{equation}
 \Lambda_n = \oplus_{s\geq 0}  \Lambda_n^s
\end{equation}
where $\Lambda_n^s$ consists of all homogeneous symmetric polynomials of degree $s.$ Schur polynomials $s_{\l}(x_1,\dots,x_n)$ such that $|\l| =s$ and $l(\l)\leq n$ form a basis of $\Lambda_n^s.$ Since length $l(\l)$ of any partition $\l$ is less or equal to its weight $|\l|,$ Schur polynomials labelled by  all partitions of $s$ form a basis of   $\Lambda_n^s$ provided $s \leq n.$ 

Let $m \geq n,$ and let $f(x_1,\dots,x_m)$ be an element of $\Lambda_m^s.$ Consider the projection 
\begin{equation}
\rho_{m,n}^s :  \Lambda_m^s \rightarrow \Lambda_n^s : f(x_1,\dots,x_m) \mapsto f(x_1,\dots,x_n,0,\dots,0). \label{eq:proj}
\end{equation} 
The effect of this projection on  Schur polynomials is well-known: it sends $s_{\l}(x_1,\dots,x_m)$ to $s_{\l}(x_1,\dots,x_n)$ if $l(\l) \leq n,$ and to $0$ if $l(\l) > n.$ Thus this projection is an isomorphism of linear spaces for $m\geq n \geq s.$  

Form the inverse limit 
\begin{equation}
\Lambda^s =  \lim\begin{Sb} \longleftarrow \\ n \end{Sb}\Lambda_n^s  \label{eq:Lambda^s}
\end{equation}
of linear spaces $\Lambda_n^s$ relative to the homomorphisms $\rho_{m,n}^s.$ The linear space 
\begin{equation}
\Lambda = \oplus_{s\geq 0} \Lambda^s 
\end{equation}
is called the space of symmetric functions.

Now recall, that we have the isomorphism (\ref{eq:Om}) of linear spaces $F_{N,n}^0$ and $\Lambda_n.$ In this section we denote this isomorphism by $\Omega_n.$ Thus for each $k_{[n]} \in \LC_n^{(1)}$ such that $k_{[n]}  \supseteq o_{[n]}$ we have 
\begin{equation}
\Omega_n( u_{k_{[n]}} ) = s_{k_{[n]} - o_{[n]}}(x_1,\dots,x_n).
\end{equation}
Obviously this isomorphism respects the principal grading of $F_{N,n}^0$ and the grading of $\Lambda_n.$ Moreover comparing definitions of the projections (\ref{eq:princproj}) and (\ref{eq:proj}) we see that for each $s \geq 0$ the map  
\begin{equation}
\Omega_{\infty}^s : F_{N}^{0,\langle s \rangle } \rightarrow  \Lambda^s : ( f_n )_{n\geq 0} \mapsto ( \Omega_n(f_n) )_{n\geq 0}\label{eq:isomLFprinc}
\end{equation}
is an isomorphism of linear spaces.

The Fock space $F$ is the  direct sum $ \oplus_{s\geq 0} F^{\langle s \rangle }.$ Recall  that  for each $s \geq 0$ we have the isomorphism    of linear spaces $ F^{\langle s \rangle } $ and $ F_N^{0,\langle s \rangle }$ by (\ref {eq:princ-fock}).
From (\ref{eq:isomLFprinc}) it now follows that $F$ with principal grading and $\Lambda$ are isomorphic as graded linear spaces. Let $\Omega_{\infty}: F \rightarrow \Lambda$ be the corresponding isomorphism map. For each partition $\l$  the element of $\Lambda$ defined by $ \Omega_{\infty}( \hat{u}_{\l +  o} ) $ is known as  Schur function $s_{\l}.$ Let $\Pi$ be the set of all partitions. Since the set of normally ordered wedges  $\{ \hat{u}_{\l +  o} \:|\: \l \in \Pi \}$ is a basis of $F,$ the set of all Schur functions $\{ s_{\l} \:|\: \l \in \Pi \}$ is a basis of $\Lambda.$

Proposition \ref{p:hom-vs-princ} entails 
\begin{prop} \label{p:hom-vs-symm}
Let $s\geq 0,$ and let  $ h = (h_n)_{n\geq 0}$ be a vector in $\Lambda^s.$ Suppose that for each $r\geq 0$ we have $\Omega_{rN}^{-1}(h_{rN}) \in F_{N,rN}^{0,(d)}$ for some $d \geq 0.$ 

Then the sequence $\tilde{h} = \left(\Omega_{rN}^{-1}(h_{rN})\right)_{r\geq 0}$ is a vector in $F_{N}^{0,(d)}.$ Moreover    
\begin{equation}
\Omega_{\infty}\left( \Upsilon_h^{(d)}( \tilde{h} ) \right) = h.
\end{equation}
\end{prop}

The  space of symmetric functions is isomorphic as a linear space to the  polynomial algebra $\f [p_1,p_2,p_3,\dots \; ].$  The generators $p_{i}\; (i=1,2,3,\dots)$ are known as the power-sums. For each partition $\l = (\l_1,\l_2,\dots ),$ the vector $p_{\l} \in  \f [p_1,p_2,p_3,\dots \; ]$ is defined by 
\begin{equation}
p_{\l} = p_{\l_1}p_{\l_2}p_{\l_3} \cdots \;.
\end{equation}
The set $ \{ p_{\l} \:|\: \l \in \Pi \}$ is a basis of $\Lambda = \f [p_1,p_2,p_3,\dots \; ].$ For each partition $\l = (\l_1,\l_2,\dots )$ let $m_i(\l)$ be the multiplicity of $i$ in $\l.$ Set  
\begin{equation}
z_{\l} = \prod_{i\geq 1}i^{m_i(\l)}m_i(\l)!
\end{equation}
and define a scalar product on $\Lambda$ by 
\begin{equation}
\sprodd{ p_{\l} }{ p_{\mu} } = z_{\l} \delta_{\l \mu}. \label{eq:standardsp}
\end{equation}
The basis of Schur functions is orthonormal relative to this scalar product:
\begin{equation}
\sprodd{ s_{\l} }{ s_{\mu} } =  \delta_{\l \mu}.
\end{equation}
Consider the automorphism 
\begin{equation}
\o: \Lambda \rightarrow \Lambda :  p_{i} \mapsto (-1)^{i-1}p_{i}. 
\end{equation}
This automorphism is clearly an involution. Its action on the basis of Schur functions is given by 
\begin{equation}
\o(s_{\l}) = s_{\l'}
\end{equation}
where $\l'$ stands for the partition conjugated to $\l.$

\subsection{Jack($\gl_N$) symmetric functions}  \label{sec:JSF}
Let $\l$ be a partition of $s \geq 0,$ i.e. $|\l| =s.$ Let 
\begin{equation}
P_{\l}(q,t)(x_1,\dots,x_n) \in \Lambda_n^s \otimes_{\f} \f(q,t)
\end{equation}
 be the Macdonald polynomial labelled by $\l.$ If $l(\l) > n$ we set $P_{\l}(q,t)(x_1,\dots,x_n) =0.$ It is well-known that the sequence  
\begin{equation}
P_{\l}(q,t) = \left( P_{\l}(q,t)(x_1,\dots,x_n) \right)_{n\geq 0}  \label{eq:MSF}
\end{equation}
is an element of the inverse limit $\Lambda_{\f(q)}^s =\Lambda^s \otimes_{\f} \f(q,t)$ (see \ref{eq:Lambda^s}). In other words  $P_{\l}(q,t)$ is a symmetric function of degree $s.$ This symmetric function is called the Macdonald symmetric function.

Introduce on $\Lambda_{\f(q,t)} = \oplus_{s\geq 0} \Lambda_{\f(q,t)}^s$ a scalar product by  
\begin{equation}
{\sprodd{ p_{\l} }{ p_{\mu} }}_{q,t} = z_{\l} \delta_{\l \mu} \prod_{i\geq 1}\frac{ 1 - q^{\l_i}}{1 - t^{\l_i}}.\label{eq:MSFsp}
\end{equation}

For each  $\l$ which is a  partition of $s,$ let 
\begin{equation}
m_{\l}(x_1,\dots,x_n) \in \Lambda_n^s 
\end{equation}
 be the corresponding symmetric monomial. We set $m_{\l}(x_1,\dots,x_n) =0$ if $l(\l) > n.$ One easily checks that the sequence  
\begin{equation}
m_{\l} = \left( m_{\l}(x_1,\dots,x_n) \right)_{n\geq 0}
\end{equation}
is an element of the inverse limit $\Lambda^s.$ It is called the monomial symmetric function.

\begin{prop}[\cite{MacBook}] \label{p:MSFchar}
The set of Macdonald symmetric functions $\{ P_{\l}(q,t)\: |\: \l \in \Pi \}$ is the unique basis of  $\Lambda_{\f(q,t)}$ which satisfies the following two properties:

\noindent {\em (a)} The transition matrix that expresses these symmetric functions  in terms of  monomial symmetric functions  is upper unitriangular. That is the expansion of $P_{\lambda}(q,t)$ has the form 
\begin{equation}
P_{\lambda}(q,t) = m_{\lambda} + \sum_{\mu < \lambda} u_{\lambda\mu}(q,t) m_{\mu}. \label{eq:MSF=m+}
\end{equation}
with certain coefficients $u_{\lambda\mu}(q,t)$ in $\f(q,t).$

\noindent {\em (b)} Symmetric functions   $P_{\lambda}(q,t)$ are pairwise orthogonal relative to the scalar product $(\ref{eq:MSFsp}).$ 
\end{prop}

Let $\l$ be a partition of $s.$  And let   
\begin{equation}
P_{\l}^{(\gamma,N)}(x_1,\dots,x_n) \in \Lambda_n^s 
\end{equation}
 be the Jack($\gl_N$) polynomial defined in Section \ref{sec:JP}. If $l(\l) > n$ we set $P_{\l}^{(\gamma,N)}(x_1,\dots,x_n)=0.$ Taking the limit (\ref{eq:limit}) in (\ref{eq:MSF}) we see that 
\begin{equation}
P_{\l}^{(\gamma,N)} = \left( P_{\l}^{(\gamma,N)}(x_1,\dots,x_n) \right)_{n\geq 0} \label{eq:JSF}
\end{equation}
is an element of the inverse limit $\Lambda^s.$ Thus $P_{\l}^{(\gamma,N)}$ is a symmetric function of degree $s.$ We will call it the Jack($\gl_N$) symmetric function.

Now let us introduce on $\Lambda$ another scalar product. Define 
\begin{equation}
{\sprodd{ p_{\l} }{ p_{\mu} }}_{\gamma,N} = z_{\l} \delta_{\l \mu}  \gamma^{-l_{N}(\l)} \label{eq:JSFsp}
\end{equation}
where for each partition $\l$ we set $l_{N}(\l) = \#\{\l_i >0 \: | \: \l_i \equiv 0\bmod N\}.$ This scalar product is  the degeneration  of the scalar product (\ref{eq:MSFsp}) in the limit (\ref{eq:limit}).

From the characterization of Macdonald symmetric functions given by Proposition \ref{p:MSFchar} we immediately obtain the following characterization of Jack($\gl_N$) symmetric functions:

\begin{prop} \label{p:JSFchar}
The set of {\em Jack($\gl_N$)} symmetric functions $\{ P_{\l}^{(\gamma ,N)}\: |\: \l \in \Pi \}$ is the unique basis of  $\Lambda$ which satisfies the following two properties:

\noindent {\em (a)} The transition matrix that expresses these symmetric functions  in terms of  monomial symmetric functions  is upper unitriangular. That is the expansion of $ P_{\l}^{(\gamma ,N)}$ has the form 
\begin{equation}
 P_{\l}^{(\gamma ,N)} = m_{\lambda} + \sum_{\mu < \lambda} u_{\lambda\mu}^{(\gamma ,N)}\; m_{\mu}. \label{eq:JSF=m+}
\end{equation}
with certain coefficients $u_{\lambda\mu}^{(\gamma ,N)}$ in $\f.$

\noindent {\em (b)} Symmetric functions   $P_{\l}^{(\gamma ,N)}$ are pairwise orthogonal relative to the scalar product $(\ref{eq:JSFsp}).$ 
\end{prop}

Let us now mention two other properties of Jack($\gl_N$) symmetric functions. The first is the explicit formulas for their normalizations with respect to the scalar product $(\ref{eq:JSFsp}).$ These formulas are obtained as limiting cases of the well-known formulas for normalizations of Macdonald symmetric functions \cite{MacBook}. For any partition $\l$ and a square $s \in \l$ let $b_{\l}^{(\gamma,N)}(s)$ be the expression introduced in (\ref{eq:b}). Set 
\begin{equation}  
b_{\l}^{(\gamma,N)} = \prod_{s \in \l} b_{\l}^{(\gamma,N)}(s).
\end{equation}
Then the normalization of the  Jack($\gl_N$) symmetric function  $P_{\l}^{(\gamma ,N)}$ is given by 
\begin{equation}
{\sprodd{ P_{\l}^{(\gamma ,N)}}{ P_{\l}^{(\gamma ,N)}}}_{\gamma,N} = 1/b_{\l}^{(\gamma,N)}.
\end{equation}
Another property is the so-called duality. Introduce the automorphism 
\begin{equation}
\o_{\gamma,N}: \Lambda \rightarrow \Lambda :  p_{i} \mapsto (-1)^{i-1}\gamma^{-\delta(i\equiv 0\bmod N)} p_{i}. 
\end{equation}
The duality for  Jack($\gl_N$) symmetric functions is expressed by the relation
\begin{equation}
\o_{\gamma,N}\left(P_{\l}^{(\gamma ,N)}\right) = b_{\l'}^{(\gamma^{-1},N)} P_{\l'}^{(\gamma^{-1},N)} .
\end{equation}

Now we will describe an important relationship between the eigenbases of the commutative algebra $A(\gl_N;\beta)$ described by Theorem \ref{t:Xbasis} and Jack($\gl_N$) symmetric functions. Let $k = (k_1,k_2,\dots ) \in \LC_{\infty}^{(1)}.$ Then $k - o$ is a partition, say $\l.$ Let $d= | \ov{o} -  \ov{k}|.$ In other words $d$ is the homogeneous degree of the normally ordered wedge $\hat{u}_k \in F.$ Observe that we have $l(k - o)\leq dN.$ For each $n\geq 0$ we set $k_{[n]} = (k_1,\dots,k_n).$ Then $k_{[n]}- o_{[n]} = \l$ as soon as $n \geq l(k - o).$

Let $n \geq l(k - o).$ Then the vector $X_{k_{[n]}}^{(\beta,N)}$ is an element of $F_{N,n}^{0,(d)}.$ Let $n < l(k - o),$ in this case we set $X_{k_{[n]}}^{(\beta,N)} =0.$

By Theorem \ref{t:X=P} for all $n\geq 0$ we have 
\begin{equation}
\Omega_{n}\left(X_{k_{[n]}}^{(\beta,N)}\right) = P_{\l}^{(N\beta + 1,N)}(x_1,\dots,x_n). \label{eq:lc1}
\end{equation}

Recall that the sequence 
\begin{equation}
\left(P_{\l}^{(N\beta + 1,N)}(x_1,\dots,x_n)\right)_{n\geq 0}
\end{equation}
 is an element of $\Lambda^{|\l|},$ that is a symmetric function of degree $|\l|.$ Of course it is just the  Jack($\gl_N$) symmetric function $P_{\l}^{(N\beta + 1,N)}.$ Taking into account (\ref{eq:lc1}) and  Proposition \ref{p:hom-vs-symm} we  conclude that the sequence 
\begin{equation}
\left(X_{k_{[rN]}}^{(\beta,N)}\right)_{r\geq 0} 
\end{equation}
is an element of the linear space $F_{N}^{0,(d)}.$ 

We have the isomorphism $\Upsilon_h^{(d)}$  (see \ref{eq:hom-fock}) of linear spaces $F_{N}^{0,(d)}$ and  $F^{(d)}.$  Let us define the vector $X_{k}^{(\beta,N)} \in F^{(d)}$ by 
\begin{equation}
X_{k}^{(\beta,N)} = \Upsilon_h^{(d)}\left( \left(X_{k_{[rN]}}^{(\beta,N)}\right)_{r\geq 0} \right).
\end{equation}
Proposition \ref{p:hom-vs-symm} and (\ref{eq:lc1})  lead to the semi-infinite analogue of Theorem \ref{t:X=P}:
\begin{prop}\label{p:X=JSF}
Let $\Omega_{\infty}: F \rightarrow \Lambda$ be the isomorphism of the Fock space and the space of symmetric functions defined by the assignment $\Omega_{\infty}: \hat{u}_k \mapsto s_{\l}$  for each $k \in \LC_{\infty}^{(1)} $ and $ \l = k - o.$

Then we have 
\begin{equation}
\Omega_{\infty}\left( X_{o+\l}^{(\beta,N)}\right) = P_{\l}^{(N\beta + 1,N)}
\end{equation}
for all partitions $\l.$

\end{prop}

\section{Yangian action on the Fock space representation  of $\sll$}\label{sec:YangactF}
The homogeneous construction of the Fock space $F$ given in Section \ref{sec:homconstr} allows to define an action of the Yangian $Y(\gl_N)$ on $F.$

The defining relations (\ref{eq:Ydef}) of the algebra   $Y(\gl_N)$ imply that this algebra is generated by the coefficients of the quantum determinant $A_N(u)$ (cf. Proposition \ref{p:quantumdet}) and the elements $T_{ab}^{(1)}, T_{ab}^{(2)}$ $(a,b=1,\dots,N).$ Recall that for any series $f(u)\in \f [[u^{-1}]]$ of the form\begin{equation}
f(u) = 1 + u^{-1} f^{(1)} + u^{-2} f^{(2)} +  \dots 
\end{equation}
the map $ \o_f : T_{ab}(u) \mapsto f(u)T_{ab}(u) $ extends to an automorphism of $Y(\gl_N).$ The effect of this map on the quantum determinant and the generators $T_{ab}^{(1)}, T_{ab}^{(2)}$ is given by 
\begin{align}
& \o_f : A_N(u) \mapsto f(u)f(u-1)\dots f(u-N+1)  A_N(u), \\
& \o_f : T_{ab}^{(1)} \mapsto  T_{ab}^{(1)} + \delta_{ab} f^{(1)}, \label{eq:fT1}\\  
& \o_f :  T_{ab}^{(2)} \mapsto  T_{ab}^{(2)} + f^{(1)} T_{ab}^{(1)}+ \delta_{ab} f^{(2)}.\label{eq:fT2}
\end{align}
In Section \ref{sec:YangactFn} we introduced the Yangian action $Y(\gl_N;\beta)$ on the linear space $F_{N,n}.$ In this section we will denote the generators of this action by  $T_{ab}^{(s)}(\beta;n)$ $(s=1,2,\dots )$ and the corresponding generating series by   $T_{ab}(u;\beta;n).$ The explicit forms of   $T_{ab}^{(1)}(\beta;n)$ and   $T_{ab}^{(2)}(\beta;n)$ are   
\begin{equation}
T_{ab}^{(1)}(\beta;n) = \sum_{i=1}^n E_{ba}^{(i)}; \quad  
T_{ab}^{(2)}(\beta;n) = - \sum_{i=1}^n d_i(\beta)^{(n)} E_{ba}^{(i)}  + \sum_{1\leq i < j \leq n} \sum_{c=1}^N E_{ca}^{(i)} E_{bc}^{(j)}, 
\end{equation}
where $d_i(\beta)^{(n)}$ $(i=1,\dots,n)$ are the Cherednik-Dunkl operators (\ref{eq:CD}).

\begin{prop}\label{p:inv} Let $r\geq 0.$ For all $s=1,2,\dots $ and $a,b=1,\dots,N$ the operators $T_{ab}^{(s)}(\beta;rN)$ leave invariant the following two subspaces of $F_{N,rN}:$

{\em (a)} The linear space $F_{N,rN}^0$ defined in {\em (\ref{eq:F0n})}. 

{\em (b)} For each $d\geq 0$ the  linear space $F_{N,rN}^{0,(d)} \subset F_{N,rN}^0$ defined in Section {\em \ref{sec:homconstr}}.
\end{prop}
\begin{pf} From the definition (\ref{eq:F0n}) it follows  that  
\begin{equation}
F_{N,rN}^0 = F_{N,rN} \cap \left( (z_1\cdots z_{rN} )^r  \f [z_1^{-1},\dots,z_{rN}^{-1}] \otimes V^{\otimes rN }\right).
\end{equation}
The Cherednik-Dunkl operators leave $(z_1\cdots z_{rN} )^r \f [z_1^{-1},\dots,z_{rN}^{-1}] $ invariant, hence the statement (a).

The Cherednik-Dunkl operators commute with the degree operator $D^{(rN)} = z_1\frac{\p}{\p z_1} + \dots + z_{rN}\frac{\p}{\p z_{rN}}.$ This implies (b). 
\end{pf}

\subsection{Intertwining relations} \label{sec:inter}
For each $r\geq 0$ define the series $ 1 + \sum_{s\geq 1} u^{-s} f^{(s)}(r)$ as the expansion in $u^{-1}$ of the rational function  
\begin{equation} 
f(u;r) = \prod_{s=1}^r \frac{u + (\beta^{-1} + N) s - 1}{u + (\beta^{-1} + N) s}.
\end{equation}
The first two coefficients of this series are 
\begin{equation}
f^{(1)}(r) = -r; \qquad  f^{(2)}(r) = r(r-1)/2 + (\beta^{-1} + N)r(r+1)/2. \label{eq:f1f2}
\end{equation}
Consider  the following renormalized generating series 
\begin{equation}
\ov{T}_{ab}(u;\beta;rN) = f(u;r) T_{ab}(u;\beta;rN).
\end{equation}

\begin{prop} \label{p:inter}
For each $r\geq 0$ let $\pi_{r+1,r}: F_{N,rN+N}^0 \rightarrow  F_{N,rN}^0$ be the projection defined in Section {\em \ref{sec:homconstr}}.

For all $a,b =1,\dots,N$ we have the intertwining relations 
\begin{equation}
\pi_{r+1,r} \ov{T}_{ab}(u;\beta;rN+N) = \ov{T}_{ab}(u;\beta;rN)\pi_{r+1,r}.
\end{equation}

\end{prop}

\begin{pf}
Since the algebra $Y(\gl_N)$ is generated by the elements $T_{ab}^{(1)}, T_{ab}^{(2)}$ $(a,b =1,\dots,N)$ and the coefficients of the quantum determinant $A_N(u),$ it is sufficient to prove the intertwining relations 
\begin{align}
&\pi_{r+1,r} \ov{T}_{ab}^{(1)}(\beta;rN+N) = \ov{T}_{ab}^{(1)}(\beta;rN)\pi_{r+1,r}, \label{eq:i1} \\
&\pi_{r+1,r} \ov{T}_{ab}^{(2)}(\beta;rN+N) = \ov{T}_{ab}^{(2)}(\beta;rN)\pi_{r+1,r}, \label{eq:i2} \\
&\pi_{r+1,r} \ov{A}_N(u;\beta;rN+N) = \ov{A}_N(u;\beta;rN)\pi_{r+1,r}.\label{eq:i3}
\end{align}
We will prove the first two relations. The relation (\ref{eq:i3}) is proved in a similar way. 

Let $k =(k_1,\dots,k_{rN+N})$ be an element of the set $\LC_{rN+N}^{(1)}$ such that $k \supseteq o$ where $o = (0,-1,\dots,-rN-N+1)$ is the vacuum sequence.  By the definition (see \ref{eq:F0n}) of the linear space $F_{N,rN+N}^{0}$ the set of normally ordered wedges $\{ \hat{u}_k \: |\: k \in  \LC_{rN+N}^{(1)}, k \supseteq o \}$ is a basis of $F_{N,rN+N}^{0}.$ We will prove (\ref{eq:i1}) and (\ref{eq:i2}) by comparing actions of both sides on elements of this basis.

First we prove the relation (\ref{eq:i1}). Consider the expression 
\begin{equation}
{T}_{ab}^{(1)}(\beta;rN+N) \hat{u}_k  = \sum_{i=1}^{rN+N} u_{k_1}\wedge \cdots  \wedge E_{ba} u_{k_i}\wedge  \cdots  \wedge u_{k_{rN+N}}.\label{eq:i4}
\end{equation}

Suppose $\hat{u}_k  \in \Ker \pi_{r+1,r}.$ By  definition (\ref{eq:homproj}) this means that $l ( k - o) > rN$ and, consequently, $ \ov{k}_{rN+1} < \ov{o}_{rN+1}.$ The right-hand side of (\ref{eq:i4}) is a linear combination of normally ordered wedges $\hat{u}_m$ such that $\ov{m} = \ov{k}.$ Hence, in particular, $ \ov{m}_{rN+1} <  \ov{o}_{rN+1}.$ The last inequality implies $l ( m - o) > rN.$ Therefore $\hat{u}_m \in \Ker \pi_{r+1,r},$ and we have 
\begin{equation}
\pi_{r+1,r} {T}_{ab}^{(1)}(\beta;rN+N)\hat{u}_k = 0 = {T}_{ab}^{(1)}(\beta;rN)\pi_{r+1,r}\hat{u}_k. \label{eq:i5} 
\end{equation}

Suppose now that  $\hat{u}_k \in F_{N,rN+N}^{0},$ $ \hat{u}_k \not\in \Ker \pi_{r+1,r}.$ Then for all $i=rN+1,\dots,rN+N$ the definition (\ref{eq:homproj}) implies $ k_i = k_i^{(0)} = -i+1.$ The right-hand side of (\ref{eq:i4}) takes the form 
\begin{equation}
\left( {T}_{ab}^{(1)}(\beta;rN) \hat{u}_{k'} \right)\wedge u_{-rN}\wedge \cdots \wedge u_{-rN-N+1} + \delta_{ab} \hat{u}_k
\end{equation}
where $k' = (k_1,\dots,k_{rN}).$ Therefore we have 
\begin{equation}
\pi_{r+1,r} {T}_{ab}^{(1)}(\beta;rN+N)\hat{u}_k =  \left( {T}_{ab}^{(1)}(\beta;rN) + \delta_{ab}1\right)\pi_{r+1,r} \hat{u}_k.
\end{equation}
Which is equivalent to 
\begin{equation}
\pi_{r+1,r} \ov{T}_{ab}^{(1)}(\beta;rN+N)\hat{u}_k =   \ov{T}_{ab}^{(1)}(\beta;rN)\pi_{r+1,r} \hat{u}_k
\end{equation}
because of (\ref{eq:fT1}) and (\ref{eq:f1f2}). Thus (\ref{eq:i1}) is proven. 

We come now to the proof of the intertwining relation (\ref{eq:i2}). 

\begin{lemma}  \label{l:i1}
Let $A_n$ be the operator {\em (\ref{eq:antisymm})} of the total antisymmetrization in the linear space $\cz \otimes V\pn.$ For all $a,b = 1,\dots,N$ we have on the space  $\cz \otimes V\pn$ the following  operator identity:
\begin{equation} 
{T}_{ab}^{(2)}(\beta;n) A_n = A_n \check{T}_{ab}^{(2)}(\beta;n)   
\end{equation} 
where 
\begin{equation} 
\check{T}_{ab}^{(2)}(\beta;n) = - \sum_{i=1}^n d_i(\beta)^{(n)} E_{ba}^{(i)}  + \sum_{1\leq j < i \leq n} \sum_{c=1}^N E_{ca}^{(i)} E_{bc}^{(j)}.  \nonumber 
\end{equation} 
\end{lemma}
\begin{pf}
Taking into account the  relation $(K_{ij}+P_{ij})A_n = 0$ we obtain 
\begin{equation}
{T}_{ab}^{(2)}(\beta;n) A_n  = -\left( \sum_{i=1}^n \ov{d}_i (\beta)^{(n)} E_{ba}^{(i)} \right) A_n 
\end{equation}
where for each $i=1,\dots,n :$  
\begin{equation}
\ov{d}_i (\beta)^{(n)} = {d}_i (\beta)^{(n)} + \sum_{j < i} K_{ij}. \nonumber 
\end{equation}
The operators $\ov{d}_i (\beta)^{(n)}$ are covariant with respect to the permutation operators $K_{ij}.$
That is 
\begin{equation}
 K_{ij} \ov{d}_j (\beta)^{(n)}  = \ov{d}_i (\beta)^{(n)}K_{ij}, \quad  K_{ij} \ov{d}_k (\beta)^{(n)} = \ov{d}_k (\beta)^{(n)} K_{ij}  \quad \text{if $k \neq i,j.$} \nonumber 
\end{equation}
This implies that $\sum_{i=1}^n \ov{d}_i (\beta)^{(n)} E_{ba}^{(i)}$ commutes with $A_n.$ Application of  $A_n (K_{ij}+P_{ij}) = 0$  proves the lemma.
\end{pf}

For each $r\geq 0$ define $L_r \subset \f [z_1^{\pm 1},\dots,z_{rN}^{\pm 1}] $ as the linear span of monomials $z^m$ such that $m_i \leq r$ for all $i=1,\dots,rN$ and $ \#\{ m_i \: |\: m_i = r\} < N.$   

\begin{lemma} \label{l:i2}
Let $r\geq 0$ and suppose  $z^m \in L_{r+1}.$ Then for any $ v \in  V^{\otimes (rN+N)}$ we have \begin{equation} \pi_{r+1,r} A_{rN+N} ( z^m \otimes v ) =0. \end{equation}
\end{lemma}

\begin{pf} The expression $A_{rN+N} ( z^m \otimes v )$ is either zero, or is equal to a linear combination of normally ordered wedges $\hat{u}_k$ such that $\ov{k}$ is a permutation of $m.$ Since for a normally ordered wedge $\hat{u}_k$  we have $ \ov{k}_1 \leq \ov{k}_2 \leq \dots \leq \ov{k}_{rN+N},$ it follows that $\ov{k}_{rN+1} < r + 1.$ On the other hand  $\ov{o}_{rN+1} = r+1$ by definition. Thus the length of the partition $k - o$ is greater than $rN.$ Definition (\ref{eq:homproj}) of the projection $\pi_{r+1,r}$ gives now the required statement.
 \end{pf}

Let $\hat{u}_k $ be a normally ordered wedge from $ F_{N.rN+N}^0.$ Consider the expression    
\begin{gather} 
{T}_{ab}^{(2)}(\beta;rN+N)\hat{u}_k = A_{rN+N}\left( \check{T}_{ab}^{(2)}(\beta;rN+N)  z^{\ov{k}} \otimes v(\un{k})\right) = \label{eq:i6} \\ = A_{rN+N}\left( - \sum_{i=1}^{rN+N} d_i(\beta)^{(rN+N)} z^{\ov{k}} \otimes E_{ba}^{(i)} v(\un{k})  + z^{\ov{k}} \otimes  \sum_{1\leq j < i \leq rN+N} \sum_{c=1}^N E_{ca}^{(i)} E_{bc}^{(j)}v(\un{k})\right).  \nonumber 
\end{gather}

Suppose to start with that $\hat{u}_k \in \Ker \pi_{r+1,r},$ that is  $l(k - o) > rN.$ Then $\ov{k}_{rN+1} < \ov{o}_{rN+1}$ and hence $ z^{\ov{k}} \in L_{r+1}.$  By Lemma \ref{l:dm=dn} we have $ d_i(\beta)^{(rN+N)} z^{\ov{k}} \in L_{r+1}$  for all $i=1,\dots,N.$ By Lemma \ref{l:i2} we now have        
\begin{equation}
{T}_{ab}^{(2)}(\beta;rN+N)\hat{u}_k \in \Ker \pi_{r+1,r}.
\end{equation}
Which gives  
\begin{equation}
\pi_{r+1,r} {T}_{ab}^{(2)}(\beta;rN+N)\hat{u}_k  = 0 = {T}_{ab}^{(2)}(\beta;rN)\pi_{r+1,r}\hat{u}_k .
\end{equation}

Now suppose $\hat{u}_k \not\in \Ker \pi_{r+1,r}.$ Then for all $i=rN+1,\dots,rN+N$ we necessarily have $k_i =  o_i =   -i+1,$ and hence  $ \ov{k}_i =  \ov{o}_i = r+1.$ Since $\ov{k}$ can have no more than $N$ equal elements, we also have $\ov{k}_i < r+1$ if $ i \leq rN.$

By Lemma \ref{l:dm=dn}(eq.\ref{eq:dm=dn1}) we obtain 
\begin{equation}
d_i(\beta)^{(rN+N)}z^{\ov{k}} \equiv d_i(\beta)^{(rN)}z^{\ov{k}} \bmod L_{r+1} \quad (i=1,\dots,rN). \nonumber
\end{equation}
And by  Lemma \ref{l:dm=dn}(eq.\ref{eq:dm=dn2}) we obtain 
\begin{equation}
d_i(\beta)^{(rN+N)}z^{\ov{k}} \equiv (\beta^{-1}(r+1)-i+2rN + N) z^{\ov{k}} \bmod L_{r+1} \quad (i=rN+1,\dots,rN+N).\nonumber
\end{equation}
Application of Lemma \ref{l:i2} and Lemma \ref{l:i1} allows to transform the right-hand side of (\ref{eq:i6}) and get 
\begin{gather}
T_{ab}^{(2)}(\beta;rN+N)\hat{u}_k \equiv \\ 
\equiv  
\left( T_{ab}^{(2)}(\beta;rN) \hat{u}_{k'} \right)  
\wedge u_{-Nr}\wedge \cdots \wedge u_{-Nr-N+1} +  \nonumber \\ +  \left(  T_{ab}^{(1)}(\beta;rN) \hat{u}_{k'} \right) \wedge u_{-Nr}\wedge \cdots \wedge u_{-Nr-N+1} - \nonumber \\  
- \delta_{ab}\left((\beta^{-1}+N)(r+1)  -1\right) \hat{u}_k   \quad \bmod \Ker \pi_{r+1,r} \nonumber 
\end{gather}
where $k'= (k_1,\dots,k_{rN}).$ This entails 
\begin{equation}
\pi_{r+1,r} T_{ab}^{(2)}(\beta;rN+N)\hat{u}_k = \left\{ T_{ab}^{(2)} + T_{ab}^{(1)} -\delta_{ab}\left((\beta^{-1}+N)(r+1)  -1\right)\right\} \pi_{r+1,r} \hat{u}_k. \nonumber 
\end{equation}
Which is equivalent to 
\begin{equation}
\pi_{r+1,r} \ov{T}_{ab}^{(2)}(\beta;rN+N)\hat{u}_k =   \ov{T}_{ab}^{(2)}(\beta;rN)\pi_{r+1,r} \hat{u}_k
\end{equation}
because of (\ref{eq:fT2}) and (\ref{eq:f1f2}). Thus (\ref{eq:i2}) is proven. 

\end{pf}

\subsection{Yangian action on the Fock space} \label{sec:yactfock}
Now we are ready to define an action of the algebra $Y(\gl_N)$ on the Fock space module $F$ of the affine Kac-Moody algebra $\sll.$ Propositions \ref{p:inv} and \ref{p:inter} allow us to define a Yangian action on the inverse limit $F_N^{0,(d)}$ (see (\ref{eq:hom})) for each $d \geq 0.$ Let $ g = (g_r)_{r\geq 0}$ be a vector in   $F_N^{0,(d)}.$ For each generator $T_{ab}^{(s)}$ of the Yangian algebra $Y(\gl_N)$ we define the corresponding action on   $F_N^{0,(d)}$  by the map 
\begin{equation}
T_{ab}^{(s)} : \left( g_r \right)_{r\geq 0} \mapsto  \left( \ov{T}_{ab}^{(s)}(\beta;rN)g_r \right)_{r\geq 0}. \label{eq:defYactF}
\end{equation}
Proposition \ref{p:inter}  guarantees that this action is well-defined.

The inverse limit  $F_N^{0,(d)}$ is isomorphic as a linear space to the component $F^{(d)}$ of the Fock space $F$ by (\ref{eq:isomhomfock}). By this isomorphism  the Fock space admits a homogeneous degree preserving $Y(\gl_N)$-action. Let us denote this action by $Y(\gl_N;\beta),$ and by $T_{ab}^{(s)}(\beta)$ $(a,b=1,\dots,N),$ $(s=1,2,\dots)$ the generators of this action. 
 
\subsection{Decomposition of the Fock space into irreducible $Y(\gl_N)$-submodules} \label{sec:F-Y}From now on  we will  identify  the Fock space module $F$ of the algebra $\sll$ and the space of symmetric functions $\Lambda$ by the isomorphism of linear spaces   
\begin{equation}    
\Omega_{\infty} : F \rightarrow \Lambda : \hat{u}_k \mapsto s_{k - o} \quad ( k \in \LC_{\infty}^{(1)} ).
\end{equation}

\mbox{} In Section \ref{sec:yactfock} we introduced the $Y(\gl_N)$-action  $Y(\gl_N;\beta)$ on the  Fock space $F.$ Let $A(\gl_N;\beta)$ be the corresponding action of the maximal commutative subalgebra $A(\gl_N) \subset Y(\gl_N).$ And let the $A_1(u;\beta),\dots,A_N(u;\beta) \in \End{(F)}[[u^{-1}]]$ denote the generating series for elements of  $A(\gl_N;\beta).$ 

The Fock space has a  basis formed by the Jack($\gl_N$) symmetric functions $P_{\l}^{(N\beta+1,N)}$ (see Section \ref{sec:JSF}) where $\l$ runs through the set of all partitions ${\Pi}.$

\begin{prop} \label{p:JSF=eigenbase}
The set $\{ P_{\l}^{(N\beta+1,N)} \: | \: \l \in \Pi \}$ is the unique, up to normalization of eigenvectors, eigenbasis of the commutative algebra $A(\gl_N;\beta)$ in the Fock space $F.$ For each $m=1,\dots,N$ we have  
\begin{gather}
 A_m(u;\beta) P_{\l}^{(N\beta+1,N)} = A_m(u;\beta;\l) P_{\l}^{(N\beta+1,N)}  \\ 
 \text{where} \quad  A_m(u;\beta;\l) = \prod_{i\geq 1} \frac{u+\beta^{-1}\ov{(o + \l)}_i+i-1 + \delta(\un{(o + \l)}_i \leq m) }{u+\beta^{-1}\ov{o}_i+i-1 + \delta(\un{o}_i \leq m) }\frac{u+\beta^{-1}\ov{o}_i+i-1  }{u+\beta^{-1}\ov{(o + \l)}_i+i-1 } \nonumber  \end{gather}
and  $o = (0,-1,-2,\dots )$ is  the vacuum sequence.
\end{prop}
\begin{pf}
Theorem \ref{t:Xbasis}, Proposition \ref{p:X=JSF} and the definition of the action $Y(\gl_N;\beta)$ given in Section \ref{sec:yactfock} give the required statement.
\end{pf}

\begin{df}[Colouring] A partition $\l$ is said to be coloured if, for all $i=1,2,\dots $ and $j=1,\dots,\l_i,$ the $(i,j)$th square of $\l$ is filled with value $ j-i \pmod N .$ 
\end{df}

\noindent For example, colouring the parition $(6,4,4,3,1)$  for $N=3$  we obtain 
\begin{center}
\unitlength=0.75pt
\begin{picture}(120,110) 
\put(0,100){\line(1,0){120}}
\put(0,80){\line(1,0){120}}
\put(0,60){\line(1,0){80}}
\put(0,40){\line(1,0){80}}
\put(0,20){\line(1,0){60}}
\put(0,0){\line(1,0){20}}
\put(0,0){\line(0,1){100}}
\put(20,0){\line(0,1){100}}
\put(40,100){\line(0,-1){80}}
\put(60,100){\line(0,-1){80}}
\put(80,100){\line(0,-1){60}}
\put(100,100){\line(0,-1){20}}
\put(120,100){\line(0,-1){20}}
\put(0,80){\makebox(20,20){0}}\put(20,80){\makebox(20,20){1}}\put(40,80){\makebox(20,20){2}}
\put(60,80){\makebox(20,20){0}}\put(80,80){\makebox(20,20){1}}\put(100,80){\makebox(20,20){2}}
\put(0,60){\makebox(20,20){2}}\put(20,60){\makebox(20,20){0}}\put(40,60){\makebox(20,20){1}}\put(60,60){\makebox(20,20){2}}
\put(0,40){\makebox(20,20){1}}\put(20,40){\makebox(20,20){2}}\put(40,40){\makebox(20,20){0}}\put(60,40){\makebox(20,20){1}}
\put(0,20){\makebox(20,20){0}}\put(20,20){\makebox(20,20){1}}\put(40,20){\makebox(20,20){2}}
\put(0,0){\makebox(20,20){2}}
\end{picture}  
\end{center}
In what follows we will consider each partition to be coloured.

\begin{df} For each partition $\l$ and $i=0,1,\dots,N-1$ let $c_i$ be the multiplicity of the colour $i$ in $\l.$ Define 
\begin{equation}
\wt(\l) = \Lambda_0 - \sum_{i=0}^{N-1} c_i \alpha_i . 
\end{equation}
\end{df}
Propositions \ref{p:wtn} and \ref{p:X=JSF} now imply  
\begin{prop} \label{p:weight} For each partition $\l$ the {\em Jack($\gl_N$)} symmetric function $P_{\l}^{(N\beta + 1,N)}$ is a weight vector of $\sll$ of the weight $\wt(\l).$ 
\end{prop}

Let $o = (0,-1,-2,\dots )$ be the vacuum sequence. The corresponding sequence $\ov{o}$ (see Section \ref{sec:preliminary}) is $ ( (1)^N (2)^N (3)^N \dots ),$ and the corresponding sequence $\un{o}$ is $((N,N-1,\dots,1)^{\infty}).$ 
\begin{df} \label{df:W} A semi-infinite sequence of integers $m=(m_1,m_2,m_3,\dots)$ is called a {\em coordinate configuration} iff the following three conditions are satisfied

{\em (i)} $m_{i+1} \geq m_{i}$ for all $i\in \nat.$

{\em (ii)} Multiplicity of each element in $m$ does not exceed $N.$ That is $m_{i+N} > m_i$ for all $i\in \nat.$

{\em (iii)} $m_i = \ov{o}_i$ for all but finite number of $i\in \nat.$
\end{df}
The set of all coordinate configurations will be  denoted by $W.$ If $m \in W,$ we may represent $m$ as $((r_1)^{p_1}(r_2)^{p_2}(r_3)^{p_3}\dots)$ where $r_1 < r_2 < r_3 < \dots ,$ and $p_i$ denotes the multiplicity of $r_i$ in $m.$ Due to the condition (ii) in Definition \ref{df:W} we have $ p_i \leq N,$ and due to the asymptotic condition (iii) we have $p_i = N$ for all large enough $i.$ We set $l(m)$ to be equal the maximal $i$ such that $p_i < N.$

\begin{df}
A semi-infinite sequence $( \ep_i )_{i\geq 1}$ such that $ \ep_i \in \{1,\dots,N\}$ $(i\in \nat)$ is called a {\em spin configuration } iff we have the equalities $ \ep_i = \un{o}_i $ for all but finite number of $i\in \nat.$ The sequence $ (\un{o}_i )_{i\geq 1}  $ is called the {\em vacuum spin configuration.} The set of all spin configuration is denoted by $\Sigma.$    
\end{df}

It is easy to see that for each partition $\l$ the sequence $m(\l) = \ov{o+\l}$ is a coordinate configuration, and the 
sequence $\ep(\l) = \un{o+\l}$ is a spin configuration, moreover, the map 
\begin{equation}
(m,\ep): \Pi \rightarrow W\times \Sigma : \l \mapsto (m(\l),\ep(\l)) 
\end{equation}
defines a one-to-one correspondence between the set of all partitions, $\Pi,$ and the set  $ W\times \Sigma.$ For each partition $\l$  we will call $m(\l)$ (resp. $\ep(\l)$) the coordinate (resp. spin) configuration of $\l.$ 

\begin{prop} \label{p:decF} As a $Y(\gl_N)$-module the Fock space representation of $\sll$ decomposes as 
\begin{equation}
F = \bigoplus_{m \in W} F(m;\beta)\label{eq:decFw1}
\end{equation}
where for  each $m = ((r_1)^{p_1} (r_2)^{p_2} (r_3)^{p_3} \dots ) \in W$ the linear space 
\begin{equation}
 F(m;\beta) = \bigoplus_{ \{\l  | m(\l) = m \}} \f P_{\l}^{(N\beta + 1,N)} \label{eq:decFw2}
\end{equation}
 is invariant and irreducible with respect to the Yangian action $Y(\gl_N;\beta).$ Moreover up to some automorphism of the  form $\o_f,$  $ F(m;\beta)$ is isomorphic as a $Y(\gl_N)$-module to the tensor product 
\begin{equation}
V_{(1^{p_1})}( a_1)\otimes V_{(1^{p_2})}( a_2 ) \otimes \cdots \otimes V_{(1^{p_{l(m)}})}( a_{l(\mm)})\label{eq:decFw3}
\end{equation}
where $ a_s = \beta^{-1}r_s -1 + p_1 + \dots + p_s$ $(s=1,\dots,l(m)).$
\end{prop}
\begin{pf}
Propositions \ref{p:decn} and \ref{p:X=JSF} give the required statement in view of the definition of the $Y(\gl_N)$-action  $Y(\gl_N;\beta)$ given in Section \ref{sec:yactfock}.
\end{pf}

For a given $m \in W$ the set of partitions $\l$ such that $m(\l) = m$ contains the unique partition $\l^{\max}$ of minimal degree $|\l|.$ By Proposition \ref{p:weight} the $\ssl_N$ weight of this partition is maximal among all $\ssl_N$ weights of the irreducible submodule $F(m;\beta).$ Thus the Jack($\gl_N$) symmetric function $P_{\l^{\max}}^{(N\beta + 1,N)}$ is the singular vector of the irreducible Yangian submodule $F(m;\beta).$

\begin{eg} \label{ex:FS}Let $N=3.$ The basis of the irreducible $Y(\gl_3)$-submodule which corresponds to the coordinate configuration  $m=((0)^2 (1) (2)^3 (3)^3 \dots)$  is given by the Jack($\gl_3$) symmetric functions  $P_{\l}^{(3\beta+1,3)}$  labelled by the nine partitions displayed below. At the right of each partition $\l$  are shown the first three entries of the  corresponding spin configuration $\ep(\l)$ with entries displayed from the top to the bottom. The rest of the entries of each $\ep(\l)$ coincide with the vacuum spin configuration.
\begin{center}
\unitlength=10pt
\vspace{0.35cm}
\begin{picture}(7,4)(0,-4) 
\put(0,0){\line(1,0){2}}
\put(0,-1){\line(1,0){2}}
\put(0,-2){\line(1,0){2}}
\put(0,0){\line(0,-1){2}}
\put(1,0){\line(0,-1){2}}
\put(2,0){\line(0,-1){2}}
\put(0,-1){\makebox(1,1){0}}\put(1,-1){\makebox(1,1){1}}
\put(0,-2){\makebox(1,1){2}}\put(1,-2){\makebox(1,1){0}}
\put(3,-1){\makebox(1,1){2}}\put(3,-2){\makebox(1,1){1}}\put(3,-3){\makebox(1,1){1}}
\end{picture} 
{\begin{picture}(7,4)(0,-4) 
\put(0,0){\line(1,0){3}}
\put(0,-1){\line(1,0){3}}
\put(0,-2){\line(1,0){2}}
\put(0,0){\line(0,-1){2}}
\put(1,0){\line(0,-1){2}}
\put(2,0){\line(0,-1){2}}
\put(3,0){\line(0,-1){1}}
\put(0,-1){\makebox(1,1){0}}\put(1,-1){\makebox(1,1){1}}\put(2,-1){\makebox(1,1){2}}
\put(0,-2){\makebox(1,1){2}}\put(1,-2){\makebox(1,1){0}}
\put(4,-1){\makebox(1,1){3}}\put(4,-2){\makebox(1,1){1}}\put(4,-3){\makebox(1,1){1}}
\end{picture}}
{\begin{picture}(7,4)(0,-4) 
\put(0,0){\line(1,0){2}}
\put(0,-1){\line(1,0){2}}
\put(0,-2){\line(1,0){2}}
\put(0,-3){\line(1,0){1}}
\put(0,0){\line(0,-1){3}}
\put(1,0){\line(0,-1){3}}
\put(2,0){\line(0,-1){2}}
\put(0,-1){\makebox(1,1){0}}\put(1,-1){\makebox(1,1){1}}
\put(0,-2){\makebox(1,1){2}}\put(1,-2){\makebox(1,1){0}}
\put(0,-3){\makebox(1,1){1}}
\put(3,-1){\makebox(1,1){2}}\put(3,-2){\makebox(1,1){1}}\put(3,-3){\makebox(1,1){2}}
\end{picture}}  
{\begin{picture}(7,4)(0,-4) 
\put(0,0){\line(1,0){3}}
\put(0,-1){\line(1,0){3}}
\put(0,-2){\line(1,0){3}}
\put(0,0){\line(0,-1){2}}
\put(1,0){\line(0,-1){2}}
\put(2,0){\line(0,-1){2}}
\put(3,0){\line(0,-1){2}}
\put(0,-1){\makebox(1,1){0}}\put(1,-1){\makebox(1,1){1}}\put(2,-1){\makebox(1,1){2}}
\put(0,-2){\makebox(1,1){2}}\put(1,-2){\makebox(1,1){0}}\put(2,-2){\makebox(1,1){1}}
\put(4,-1){\makebox(1,1){3}}\put(4,-2){\makebox(1,1){2}}\put(4,-3){\makebox(1,1){1}}
\end{picture}}
{\begin{picture}(7,4)(0,-4) 
\put(0,0){\line(1,0){3}}
\put(0,-1){\line(1,0){3}}
\put(0,-2){\line(1,0){2}}\put(0,-3){\line(1,0){1}}
\put(0,0){\line(0,-1){3}}
\put(1,0){\line(0,-1){3}}
\put(2,0){\line(0,-1){2}}
\put(3,0){\line(0,-1){1}}
\put(0,-1){\makebox(1,1){0}}\put(1,-1){\makebox(1,1){1}}\put(2,-1){\makebox(1,1){2}}
\put(0,-2){\makebox(1,1){2}}\put(1,-2){\makebox(1,1){0}}\put(0,-3){\makebox(1,1){1}}
\put(4,-1){\makebox(1,1){3}}\put(4,-2){\makebox(1,1){1}}\put(4,-3){\makebox(1,1){2}}
\end{picture}} \\ 
{\begin{picture}(7,4)(0,-4) 
\put(0,0){\line(1,0){2}}
\put(0,-1){\line(1,0){2}}
\put(0,-2){\line(1,0){2}}\put(0,-3){\line(1,0){2}}
\put(0,0){\line(0,-1){3}}
\put(1,0){\line(0,-1){3}}
\put(2,0){\line(0,-1){3}}
\put(0,-1){\makebox(1,1){0}}\put(1,-1){\makebox(1,1){1}}
\put(0,-2){\makebox(1,1){2}}\put(1,-2){\makebox(1,1){0}}
\put(0,-3){\makebox(1,1){1}}\put(1,-3){\makebox(1,1){2}}
\put(3,-1){\makebox(1,1){2}}\put(3,-2){\makebox(1,1){1}}\put(3,-3){\makebox(1,1){3}}
\end{picture}}  
{\begin{picture}(7,4)(0,-4) 
\put(0,0){\line(1,0){3}}
\put(0,-1){\line(1,0){3}}
\put(0,-2){\line(1,0){3}}
\put(0,-3){\line(1,0){1}}
\put(0,0){\line(0,-1){3}}
\put(1,0){\line(0,-1){3}}
\put(2,0){\line(0,-1){2}}
\put(3,0){\line(0,-1){2}}
\put(0,-1){\makebox(1,1){0}}\put(1,-1){\makebox(1,1){1}}\put(2,-1){\makebox(1,1){2}}
\put(0,-2){\makebox(1,1){2}}\put(1,-2){\makebox(1,1){0}}\put(2,-2){\makebox(1,1){1}}
\put(0,-3){\makebox(1,1){1}}
\put(4,-1){\makebox(1,1){3}}\put(4,-2){\makebox(1,1){2}}\put(4,-3){\makebox(1,1){2}}
\end{picture}}
{\begin{picture}(7,4)(0,-4) 
\put(0,0){\line(1,0){3}}
\put(0,-1){\line(1,0){3}}
\put(0,-2){\line(1,0){2}}
\put(0,-3){\line(1,0){2}}
\put(0,0){\line(0,-1){3}}
\put(1,0){\line(0,-1){3}}
\put(2,0){\line(0,-1){3}}
\put(3,0){\line(0,-1){1}}
\put(0,-1){\makebox(1,1){0}}\put(1,-1){\makebox(1,1){1}}\put(2,-1){\makebox(1,1){2}}
\put(0,-2){\makebox(1,1){2}}\put(1,-2){\makebox(1,1){0}}
\put(0,-3){\makebox(1,1){1}}\put(1,-3){\makebox(1,1){2}}
\put(4,-1){\makebox(1,1){3}}\put(4,-2){\makebox(1,1){1}}\put(4,-3){\makebox(1,1){3}}
\end{picture}}  
{\begin{picture}(7,4)(0,-4) 
\put(0,0){\line(1,0){3}}
\put(0,-1){\line(1,0){3}}
\put(0,-2){\line(1,0){3}}
\put(0,-3){\line(1,0){2}}
\put(0,0){\line(0,-1){3}}
\put(1,0){\line(0,-1){3}}
\put(2,0){\line(0,-1){3}}
\put(3,0){\line(0,-1){2}}
\put(0,-1){\makebox(1,1){0}}\put(1,-1){\makebox(1,1){1}}\put(2,-1){\makebox(1,1){2}}
\put(0,-2){\makebox(1,1){2}}\put(1,-2){\makebox(1,1){0}}\put(2,-2){\makebox(1,1){1}}
\put(0,-3){\makebox(1,1){1}}\put(1,-3){\makebox(1,1){2}}
\put(4,-1){\makebox(1,1){3}}\put(4,-2){\makebox(1,1){2}}\put(4,-3){\makebox(1,1){3}}
\end{picture}}

\end{center}
\unitlength=1pt
Partition $(2,2)$ corresponds to the singular vector of this submodule.
\end{eg}

\section{Yangian action on the basic representation of $\sll$}
\subsection{The basic representation of $\sll$} 
Under the action of $\sll$ the Fock space $F$ decomposes as \cite{JM}:
\begin{equation}
F=V_N \oplus S_N
\end{equation}
where $ V_N = U(\sll) s_0 = \f[\, p_i\:|\: i\in \nat, i\not\equiv 0\bmod N],$ and $S_N$ is the ideal generated by $p_{kN}$ $(k \in \nat ).$ Moreover one has the isomorphisms of $\sll$-modules
\begin{equation}
 V_N \cong F/S_N \cong V(\Lambda_0)
\end{equation}
where $V(\Lambda_0)$ is the basic representation of $\sll.$

Our first goal in this section is to define an action of the algebra $Y(\gl_N)$ on the linear space $V_N$ starting from the $Y(\gl_N)$-action $Y(\gl_N;\beta)$ on the Fock space.

\subsection{The projection map and a Yangian action on the basic representation of $\sll$} \label{sec:Y-B} The projection map $\pi$ from $F$ onto the linear space $V_N$ is defined in the basis of the power-sums $\{ p_{\l} | \l \in \Pi \}$ as follows: 
\begin{equation}
\pi : F \rightarrow V_N : p_{\l} \mapsto \begin{cases} 0 & \text{ if $ \exists \l_i > 0 $ such that $ \l_i \equiv 0\bmod N, $ } \\ p_{\l} & \text{ otherwise.} \end{cases} \label{eq:pro}
\end{equation}

Consider the Yangian action $Y(\gl_N;\beta)$ on the Fock space $F.$ Thus far we have assumed that the parameter $\beta$ is a positive real number. If we are concerned only with the definition of the action $Y(\gl_N;\beta),$ and not with the decomposition of the Fock space into irreducible Yangian submodules, this assumption is not essential. Indeed, all the results of Sections \ref{sec:inter} and \ref{sec:yactfock} are easily seen to be valid for all $\beta$ in $\f \setminus \{0\}.$ In particular taking the limit $\beta \rightarrow \infty$ in $Y(\gl_N;\beta)$ gives rise to a well-defined $Y(\gl_N)$-action on the Fock space which we will denote by $Y(\gl_N;\infty).$ 

The following proposition is obtained as corollary to  Proposition 14 in the paper \cite{STU}.
\begin{prop} \label{p:Ycommp}The $Y(\gl_N)$-action $Y(\gl_N;\infty)$ commutes with the multiplication operator by $p_{jN}$ for each $j\in \nat.$
\end{prop}
This proposition allows to define a Yangian action on the linear space $V_N$ as a projection of the action $Y(\gl_N;\infty).$ Precisely we have 
\begin{prop} \label{p:YactB} For each generator $T_{ab}^{(s)}$ of the algebra $Y(\gl_N)$ define $\pi(T_{ab}^{(s)}) \in \End(V_N)$ by 
\begin{equation}
\pi(T_{ab}^{(s)}) \pi(v) = \pi\left( T_{ab}^{(s)}(\infty) v\right)
\end{equation}
for all $v$ in the Fock space $F.$ 

Then the assignment 
\begin{equation}
T_{ab}^{(s)} \mapsto \pi(T_{ab}^{(s)})
\end{equation}
gives rise to a $Y(\gl_N)$-action on the linear space $V_N.$
\end{prop}
\begin{pf} This follows from Proposition \ref{p:Ycommp} and the definition of the projection $\pi.$ 
\end{pf}

We will denote by $\pi(Y(\gl_N))$ and $\pi(A(\gl_N))$ the $Y(\gl_N)$-action and the action of the maximal commutative subalgebra $A(\gl_N)  \subset Y(\gl_N)$ defined by the above proposition. 

\subsection{Projecting the Jack($\gl_N$) symmetric functions} An examination of explicit expressions for the coefficients (\ref{eq:uLM}) shows that these coefficients have no poles at $\gamma=\infty$ and $\gamma=0.$ Therefore for each partition $\l$ we have well-defined  Jack($\gl_N$) symmetric functions $P_{\l}^{(\infty,N)}$ and $P_{\l}^{(0,N)}.$ We will now consider the effect of the projection $\pi$ on these symmetric functions. 

\begin{df}
A partition $\l=(\l_1,\l_2,\dots )$ is said to be {\em $N$-regular} iff $\l_i - \l_{i+1} < N$ for all $i=1,2,\dots .$ The set of all $N$-regular partitions is denoted by $\Pi_N.$ 
\end{df}
\begin{rmk} Our definition of $N$-regular partitions differs from the standard definition by conjugation of partitions.
\end{rmk}

Recall from Section \ref{sec:JSF} that there is an automorphism of the Fock space defined by 
\begin{equation}
\o_{\gamma,N}: F \rightarrow F :  p_{k} \mapsto (-1)^{k-1} \gamma^{-\delta(k\equiv 0\bmod N)} p_{k} \label{eq:omg}
\end{equation}
and that we have the duality of Jack($\gl_N$) symmetric functions expressed by 
\begin{equation}
\o_{\gamma,N}\left(P_{\l}^{(\gamma ,N)}\right) = b_{\l'}^{(\frac{1}{\gamma},N)} P_{\l'}^{(\frac{1}{\gamma},N)}  \label{eq:dual2}
\end{equation}
where for each partition $\l$ 
\begin{equation}
b_{\l}^{(\gamma,N)} = \prod_{\{s \in \l | h_{\l}(s)\equiv 0\bmod N\}}\frac{ a_{\lambda}(s) + \gamma l_{\lambda}(s) + \gamma }{ a_{\lambda}(s) + \gamma l_{\lambda}(s) + 1 }. \label{eq:bbg}
\end{equation}

Let $\o$ be the standard involution of the Fock space defined by 
\begin{equation}
\o: F \rightarrow F :  p_{k} \mapsto (-1)^{k-1} p_{k}. 
\end{equation}
This involution clearly leaves the subspace $V_N$ invariant. Moreover the definition (\ref{eq:pro}) of the projection $\pi$ and the definition (\ref{eq:omg}) of $\o_{\gamma,N}$ give the equalities 
\begin{equation}
\o \o_{\infty,N}(f) = \o_{\infty,N}\o (f)= \pi(f)  \label{eq:omg=pi}
\end{equation}
for any $f \in F.$ 

\begin{prop} \label{p:projP}For any partition $\l$ we have 
\begin{align}
&\pi(P_{\l}^{(\infty ,N)}) = \begin{cases}  b_{\l'} \o(P_{\l'}^{(0,N)}) & \text{{\em if $\l\in \Pi_N;$}}  \\ 0  & \text{{\em otherwise}}  \end{cases}\label{eq:x1}\\ 
\intertext{where for each $\l \in \Pi_N$ }
& b_{\l'} = \prod_{\{s \in \l | h_{\l}(s)\equiv 0\bmod N\}}\frac{ l_{\lambda}(s)}{ l_{\lambda}(s) + 1 } > 0.
\end{align}  
\end{prop}
\begin{pf}
From (\ref{eq:bbg}) it follows that $b_{\l'}^{(0,N)}=0$ if $\l \not\in \Pi_N$ and $b_{\l'}^{(0,N)} =  b_{\l'} >0$ if $\l \in \Pi_N.$ Setting $\gamma =\infty$ in (\ref{eq:dual2}) and taking into account (\ref{eq:omg=pi}) we obtain the required  statement.
\end{pf}

\begin{cor}\label{p:projPcor}\mbox{} 

{\em (i)} For each $\l \in \Pi_N$ we have $P_{\l'}^{(0,N)} \in V_N.$ 

{\em (ii)} The sets $\{ \pi(P_{\l}^{(\infty ,N)})\:|\: \l \in \Pi_N \}$ and $\{ P_{\l'}^{(0,N)}\:|\: \l \in \Pi_N \}$ are bases of $V_N.$
\end{cor}
\begin{pf}
(i) Follows from (\ref{eq:x1}) since $\o$ leaves $V_N$ invariant.

The symmetric functions $P_{\l}^{(0,N)}$ are linearly independent. This implies (ii).
\end{pf}
In the sequel for each $\l \in \Pi_N$ we will use the notation $P_{\l}^{(N)}$ defined by  
\begin{equation}
P_{\l}^{(N)} = \pi(P_{\l}^{(\infty,N)}) =  b_{\l'} \o(P_{\l'}^{(0,N)}).
\end{equation}
Thus $\{ P_{\l}^{(N)}\: | \: \l \in \Pi_N \}$ is a basis of $V_N.$

\subsection{Yangian decomposition of the basic representation of $\sll$}
Each $N$-regular partition $\l$ is completely determined by its spin configuration $\ep(\l) = \un{o+\l}.$ Moreover, it is straightforward to verify that the map  
\begin{equation}
\ep : \Pi_N \rightarrow \Sigma : \l \mapsto \ep(\l)
\end{equation}
gives a one-to-one correspondence of sets. Hence one can label $N$-regular partitions by the corresponding spin configurations which is sometimes convenient to do. 

\begin{prop} \label{p:piAeigenbasis}
The set $\{ P_{\l}^{(N)}\: | \: \l \in \Pi_N \}$ is an eigenbasis of the commutative algebra $\pi(A(\gl_N)).$ For each $m=1,\dots,N$ we have 
\begin{align}
& \pi(A_m(u)) P_{\l}^{(N)} = A_m(u;\l) P_{\l}^{(N)}  \\ 
& \text{where} \quad  A_m(u;\l) = \prod_{i\geq 1} \frac{u +i-1 + \delta(\ep(\l)_i \leq m) }{u+i-1 + \delta(\ep(0)_i \leq m) }. \nonumber  
\end{align}
\end{prop}
\begin{pf} Set $\beta=\infty$ in the statement of Proposition \ref{p:JSF=eigenbase} and apply the projection $\pi.$ 
\end{pf}

When $\beta =\infty ,$ the decomposition of the Fock space into invariant subspaces relative to the action $Y(\gl_N;\infty)$ is still described by Proposition \ref{p:decF} with the only difference from the case of positive real $\beta$ that now  these invariant subspaces are not, in general, irreducible $Y(\gl_N)$-modules. Setting $\beta=\infty$ in (\ref{eq:decFw1}) and    (\ref{eq:decFw2}) we have 
\begin{equation}
F = \bigoplus_{m \in W} F(m;\infty),   \label{eq:decFwinf1}
\end{equation}
where $F(m;\infty)$ is a  $Y(\gl_N)$-submodule given by  
\begin{equation}
F(m;\infty) = \bigoplus_{\{ \l \in \Pi |  m(\l) = m \} } \f  P_{\l}^{(\infty,N)}.
\end{equation}

From Proposition  \ref{p:projP} and Corollary \ref{p:projPcor} we know the effect of the projection $\pi$ (\ref{eq:pro}) on the elements of the basis $\{  P_{\l}^{(\infty,N)}\: | \: \l \in \Pi \}.$ This allows to derive the decomposition of the basic representation of $\sll$ simply  by applying $\pi$ to the both sides of (\ref{eq:decFwinf1}). 

Let $W_N'$ be the set of all coordinate configurations which correspond to $N$-regular partitions. That is to say $ W_N' = \{ m(\l) \: |\: \l \in \Pi_N \}.$ If an element $m$ appears in $W_N'$ with multiplicity, discard all the copies of $m$ but one. In this way we obtain the set $W_N$ all of whose elements are distinct, and for each element $m$ of $W_N$ there is at least one  $N$-regular partition $\l$ such that $m=m(\l).$ It is easy to describe the set $W_N.$ It consists of all semi-infinite sequences of non-negative integers $m = (m_1,m_2,\dots)$ such that $m_i=\ov{o}_i$ for all but finite number of $i\in \nat$ and for each $i\in \nat$ we have 

(i) $ m_{i+1}-m_i  \in \{0,1\}, $ 

(ii) $ m_{i+N} > m_i.$ 

Since each $N$-regular partition $\l$ is uniquely determined by its spin configuration $\ep(\l),$ the coordinate configuration $m(\l)$ is also determined by $\ep(\l).$ The formula which, together with the asymptotic condition, gives $m(\l)$ in terms of $\ep(\l)$ reads as follows:
\begin{equation}
m(\l)_i = \begin{cases}m(\l)_{i+1} -  1  & \text{ if $\ep(\l)_{i+1} \geq  \ep(\l)_{i},$} \\
m(\l)_{i+1} 1  & \text{ if $\ep(\l)_{i+1} <  \ep(\l)_{i}.$} \end{cases}
\end{equation}
Now we describe the decomposition of the basic representation $V_N$ relative to the $Y(\gl_N)$-action $\pi( Y(\gl_N)).$ 
\begin{prop} \label{p:decB1}
As a $Y(\gl_N)$-module the basic representation of $\sll$ decomposes as 
\begin{equation}
V_N = \bigoplus_{m \in W_N} V_N(m)  \label{eq:decB11}
\end{equation}
where $V_N(m)$ is a  $Y(\gl_N)$-submodule given by 
\begin{equation}
V_N(m) =  \bigoplus_{\{ \l \in \Pi_N | m(\l) = m \}} \f P_{\l}^{(N)}. \label{eq:decB12}
\end{equation}
\end{prop}
\begin{pf} Set $\beta =\infty$ in  (\ref{eq:decFw1}) and (\ref{eq:decFw2}) and apply the projection $\pi$ (\ref{eq:pro}). The required statement now follows from the  definition of the Yangian action $\pi(Y(\gl_N))$ given in Section  \ref{sec:Y-B}.
\end{pf}

\begin{eg} \label{ex:BS} Let $N=3,$ we continue from Example \ref{ex:FS}. The coordinate configuration $m=((0)^2 (1)(2)^3 (3)^3 \dots )$ is an element of $W_3.$ The corresponding invariant subspace of the Yangian action $\pi(Y(\gl_3))$ has the basis formed by symmetric functions $P_{\l}^{(3)}$ labelled by the  $3$-regular partitions shown below. As in Example \ref{ex:FS} we display the corresponding spin configuration alongside each partition.
\unitlength=10pt
\begin{center}
\vspace{0.35cm}
\begin{picture}(7,4)(0,-4) 
\put(0,0){\line(1,0){2}}
\put(0,-1){\line(1,0){2}}
\put(0,-2){\line(1,0){2}}
\put(0,0){\line(0,-1){2}}
\put(1,0){\line(0,-1){2}}
\put(2,0){\line(0,-1){2}}
\put(0,-1){\makebox(1,1){0}}\put(1,-1){\makebox(1,1){1}}
\put(0,-2){\makebox(1,1){2}}\put(1,-2){\makebox(1,1){0}}
\put(3,-1){\makebox(1,1){2}}\put(3,-2){\makebox(1,1){1}}\put(3,-3){\makebox(1,1){1}}
\end{picture} 
{\begin{picture}(7,4)(0,-4) 
\put(0,0){\line(1,0){3}}
\put(0,-1){\line(1,0){3}}
\put(0,-2){\line(1,0){2}}
\put(0,0){\line(0,-1){2}}
\put(1,0){\line(0,-1){2}}
\put(2,0){\line(0,-1){2}}
\put(3,0){\line(0,-1){1}}
\put(0,-1){\makebox(1,1){0}}\put(1,-1){\makebox(1,1){1}}\put(2,-1){\makebox(1,1){2}}
\put(0,-2){\makebox(1,1){2}}\put(1,-2){\makebox(1,1){0}}
\put(4,-1){\makebox(1,1){3}}\put(4,-2){\makebox(1,1){1}}\put(4,-3){\makebox(1,1){1}}
\end{picture}}
{\begin{picture}(7,4)(0,-4) 
\put(0,0){\line(1,0){2}}
\put(0,-1){\line(1,0){2}}
\put(0,-2){\line(1,0){2}}
\put(0,-3){\line(1,0){1}}
\put(0,0){\line(0,-1){3}}
\put(1,0){\line(0,-1){3}}
\put(2,0){\line(0,-1){2}}
\put(0,-1){\makebox(1,1){0}}\put(1,-1){\makebox(1,1){1}}
\put(0,-2){\makebox(1,1){2}}\put(1,-2){\makebox(1,1){0}}
\put(0,-3){\makebox(1,1){1}}
\put(3,-1){\makebox(1,1){2}}\put(3,-2){\makebox(1,1){1}}\put(3,-3){\makebox(1,1){2}}
\end{picture}}  
{\begin{picture}(7,4)(0,-4) 
\end{picture}}
{\begin{picture}(7,4)(0,-4) 
\put(0,0){\line(1,0){3}}
\put(0,-1){\line(1,0){3}}
\put(0,-2){\line(1,0){2}}\put(0,-3){\line(1,0){1}}
\put(0,0){\line(0,-1){3}}
\put(1,0){\line(0,-1){3}}
\put(2,0){\line(0,-1){2}}
\put(3,0){\line(0,-1){1}}
\put(0,-1){\makebox(1,1){0}}\put(1,-1){\makebox(1,1){1}}\put(2,-1){\makebox(1,1){2}}
\put(0,-2){\makebox(1,1){2}}\put(1,-2){\makebox(1,1){0}}\put(0,-3){\makebox(1,1){1}}
\put(4,-1){\makebox(1,1){3}}\put(4,-2){\makebox(1,1){1}}\put(4,-3){\makebox(1,1){2}}
\end{picture}} \\ 
{\begin{picture}(7,4)(0,-4) 
\put(0,0){\line(1,0){2}}
\put(0,-1){\line(1,0){2}}
\put(0,-2){\line(1,0){2}}\put(0,-3){\line(1,0){2}}
\put(0,0){\line(0,-1){3}}
\put(1,0){\line(0,-1){3}}
\put(2,0){\line(0,-1){3}}
\put(0,-1){\makebox(1,1){0}}\put(1,-1){\makebox(1,1){1}}
\put(0,-2){\makebox(1,1){2}}\put(1,-2){\makebox(1,1){0}}
\put(0,-3){\makebox(1,1){1}}\put(1,-3){\makebox(1,1){2}}
\put(3,-1){\makebox(1,1){2}}\put(3,-2){\makebox(1,1){1}}\put(3,-3){\makebox(1,1){3}}
\end{picture}}  
{\begin{picture}(7,4)(0,-4) 
\put(0,0){\line(1,0){3}}
\put(0,-1){\line(1,0){3}}
\put(0,-2){\line(1,0){3}}
\put(0,-3){\line(1,0){1}}
\put(0,0){\line(0,-1){3}}
\put(1,0){\line(0,-1){3}}
\put(2,0){\line(0,-1){2}}
\put(3,0){\line(0,-1){2}}
\put(0,-1){\makebox(1,1){0}}\put(1,-1){\makebox(1,1){1}}\put(2,-1){\makebox(1,1){2}}
\put(0,-2){\makebox(1,1){2}}\put(1,-2){\makebox(1,1){0}}\put(2,-2){\makebox(1,1){1}}
\put(0,-3){\makebox(1,1){1}}
\put(4,-1){\makebox(1,1){3}}\put(4,-2){\makebox(1,1){2}}\put(4,-3){\makebox(1,1){2}}
\end{picture}}
{\begin{picture}(7,4)(0,-4) 
\put(0,0){\line(1,0){3}}
\put(0,-1){\line(1,0){3}}
\put(0,-2){\line(1,0){2}}
\put(0,-3){\line(1,0){2}}
\put(0,0){\line(0,-1){3}}
\put(1,0){\line(0,-1){3}}
\put(2,0){\line(0,-1){3}}
\put(3,0){\line(0,-1){1}}
\put(0,-1){\makebox(1,1){0}}\put(1,-1){\makebox(1,1){1}}\put(2,-1){\makebox(1,1){2}}
\put(0,-2){\makebox(1,1){2}}\put(1,-2){\makebox(1,1){0}}
\put(0,-3){\makebox(1,1){1}}\put(1,-3){\makebox(1,1){2}}
\put(4,-1){\makebox(1,1){3}}\put(4,-2){\makebox(1,1){1}}\put(4,-3){\makebox(1,1){3}}
\end{picture}}  
{\begin{picture}(7,4)(0,-4) 
\put(0,0){\line(1,0){3}}
\put(0,-1){\line(1,0){3}}
\put(0,-2){\line(1,0){3}}
\put(0,-3){\line(1,0){2}}
\put(0,0){\line(0,-1){3}}
\put(1,0){\line(0,-1){3}}
\put(2,0){\line(0,-1){3}}
\put(3,0){\line(0,-1){2}}
\put(0,-1){\makebox(1,1){0}}\put(1,-1){\makebox(1,1){1}}\put(2,-1){\makebox(1,1){2}}
\put(0,-2){\makebox(1,1){2}}\put(1,-2){\makebox(1,1){0}}\put(2,-2){\makebox(1,1){1}}
\put(0,-3){\makebox(1,1){1}}\put(1,-3){\makebox(1,1){2}}
\put(4,-1){\makebox(1,1){3}}\put(4,-2){\makebox(1,1){2}}\put(4,-3){\makebox(1,1){3}}
\end{picture}}


\end{center}
\unitlength=1pt
Partition $(3,3)$ which was present in the decomposition of the submodule $F(m;\infty)$ is not $3$-regular and therefore the corresponding Jack($\gl_N$) symmetric function $P_{(3,3)}^{(\infty,3)} $ vanishes under the action of the projection $\pi.$  
\end{eg}

\subsection{Irreducibility of the Yangian decomposition of the basic representation of $\sll$} Our goal now is to prove irreducibility of each of the  Yangian submodules $V_N(m)$ $(m\in W_N)$ and to compute corresponding Drinfeld polynomials. To do this it is convenient to parameterize the set of all $N$-regular partitions, or, equivalently, the set of all spin configurations by semi-standard tableaux of ribbon (Young) diagrams. This parameterization was introduced and its relation to the Yangian decomposition was conjectured in the paper \cite{KKN}.

\subsubsection{Ribbon diagrams} Let $\theta$ be a skew diagram. Two squares in $\theta$ are {\em adjacent } if they share a common side. A skew diagram $\theta$ is {\em connected} if for any pair of squares $s$ and $s'$ in $\theta$ there is a series of squares $s_1=s,s_2,\dots,s_n=s'$ in $\theta$ such that $s_i$ and $s_{i+1}$ are adjacent. A skew diagram is called a {\em ribbon} if it is connected and contains no $2\times 2$ blocks of squares. The {\em length} of a ribbon is the total number of squares it contains. We will let  $[p_1,\dots,p_l]$ denote the ribbon of $l$ columns such that the height of $i$th column from the right is $p_i.$ A ribbon is said to be of {\em rank} $N$ if heights of all its columns do not exceed $N.$ 

\begin{eg} \label{ex:ribbon} Here is the ribbon $[3,2,1,1,1,2,1,1,1,2].$ It has the length $15$ and is of the rank $3.$
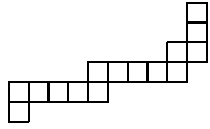
\begin{center}
\unitlength=0.75pt
\begin{picture}(100,70) 
\put(90,60){\line(1,0){10}}
\put(90,50){\line(1,0){10}}
\put(80,40){\line(1,0){20}}
\put(40,30){\line(1,0){60}}
\put(0,20){\line(1,0){90}}
\put(0,10){\line(1,0){50}}
\put(0,0){\line(1,0){10}}
\put(0,0){\line(0,1){20}}
\put(10,0){\line(0,1){20}}
\put(20,10){\line(0,1){10}}
\put(30,10){\line(0,1){10}}
\put(40,10){\line(0,1){20}}
\put(50,10){\line(0,1){20}}
\put(60,20){\line(0,1){10}}
\put(70,20){\line(0,1){10}}
\put(80,20){\line(0,1){20}}
\put(90,20){\line(0,1){40}}
\put(100,30){\line(0,1){30}}

\end{picture}  
\end{center}

\end{eg}
Let $R_N$ be the set of all rank $N$ ribbons $\theta$ which satisfy the following two conditions  

(i) The leftmost column of $\theta$ is of the height less than $N.$ 

(ii) The length of $\theta$ is divisible by $N.$

\noindent Thus the  ribbon of Example \ref{ex:ribbon} belongs to each of the sets $R_3,R_5,R_{15}.$

Let $m =(m_i)_{i\geq 1}$ be a coordinate configuration from the set $W_N.$ Recall that by the definition of $W_N$ the difference  $m_{i+1}-m_{i}$ is either $0$ or $1$ for all $i.$ Moreover if we represent $m$ as $((r_1)^{p_1}(r_2)^{p_2}(r_3)^{p_3} \dots)$ then $r_{i+1}= r_i+1$ $(i\in \nat),$ and the multiplicities $p_i$ do not exceed $N.$ For a given coordinate configuration $m \in W_N$ we associate a ribbon $\theta(m)$ by the following procedure 

1. Write the first square. 

2. For all $i=1,2,\dots,p_1+p_2+\dots+p_{l(m)}-1$ attach the $i+1$th square at the bottom (resp. at the left) of the  $i$th square if $m_{i+1}=m_{i}$ (resp. $m_{i+1}=m_{i}+1$).

\noindent It is easy to see that $\theta(m) = [p_1,\dots,p_{l(m)}].$ Hence $\theta(m)$ is of rank $N,$ its leftmost column has the length $p_{l(m)}< N$ by the definition of $l(m)$ (see Section \ref{sec:F-Y}) and the length of $\theta(m)$ is divisible by $N$ due to the asymptotic condition on the coordinate configuration $m.$    

\begin{prop} The correspondence of sets 
\begin{equation}
\theta : W_N \rightarrow R_N : m \mapsto \theta(m)
\end{equation}
is one-to-one.
\end{prop}
\begin{pf} 
The injectivity of $\theta$ follows since each coordinate configuration 
\begin{equation}
m = (m_i)_{i\geq 1}=((r_1)^{p_1}(r_2)^{p_2}(r_3)^{p_3} \dots) \in W_N
\end{equation}
is uniquely determined by its multiplicities $p_1,\dots,p_{l(m)}$ in view of the condition $r_{i+1} = r_i + 1$ $(i\in \nat),$ and the asymptotic condition  $m_i = \ov{o}_i$ $ (i \gg 1).$  

Let $[p_1,\dots,p_l]$ be a ribbon from the set $R_N,$ in particular we have $p_1+\dots+p_l =s N$ for a certain $s\geq 0.$ The sequence  
\begin{equation}
m = ((s-l+1)^{p_1}(s-l+2)^{p_2}\dots (s)^{p_l} (s+1)^N (s+2)^N \dots )
\end{equation}
is an element of the set $W_N$ such that $\theta(m) = [p_1,\dots,p_l].$ This shows that $\theta$ is surjective.  
\end{pf}

\begin{eg} The coordinate configuration $m=((0)^2 (1) (2)^3 (3)^3 \dots )$ is an element of the set $W_3$ with  $l(m)=2.$ The ribbon which corresponds to this $m$ is   

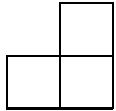
\begin{center}
\begin{picture}(40,60)(0,-10) 
\put(20,40){\line(1,0){20}}
\put(0,20){\line(1,0){40}}
\put(0,0){\line(1,0){40}}
\put(0,0){\line(0,1){20}}
\put(20,0){\line(0,1){40}}
\put(40,0){\line(0,1){40}}
\end{picture}  
\end{center}
\end{eg}

\subsubsection{Semi-standard tableaux}In each square of a given ribbon $\theta$ inscribe one of the numbers $1,2,\dots,N.$ We will call such an arrangement of numbers a {\em semi-standard tableau T of shape} $\theta$ if the numbers decrease downward along each column and do not increase from left to right along each row. 

\begin{rmk}
Note that the definition of a semi-standard tableau which we use is different from the usual definition (see e.g. \cite{KKN},\cite{MacBook}) where 
the inscribed numbers strictly increase downward along  columns and do not increase from right to left along rows.   
\end{rmk}

Let $n$ be the length of the ribbon $\theta,$ then $T$ is uniquely represented as the sequence $(\ep_1,\ep_2,\dots,\ep_n)$ where $\ep_i$ is the number inscribed in the $i$th square along $\theta$ counting from the right to the left  and from the top to the bottom. Now we have a natural map $\phi_m$ from the set of semi-standard tableaux of shape $\theta(m)$ to the set of spin configurations $\Sigma$ defined by  
\begin{equation}
\phi_m : (\ep_1,\ep_2,\dots,\ep_n) \mapsto (\ep_1,\ep_2,\dots,\ep_n, (N,N-1,\dots,1)^{\infty}).
\end{equation}
For each element $m$ of the set $W_N$ we define $\Sigma(m) \subset \Sigma$ by  
\begin{equation}
\Sigma(m) = \bigsqcup_{\{ \l \in \Pi_N | m(\l) = m \}} \{ \ep(\l) \}.
\end{equation}
The following proposition is proved in \cite{KKN}.
\begin{prop}
For each $m \in W_N$ the map $\phi_m$ gives a one-to-one correspondence between the set of semi-standard tableaux of shape $\theta(m)$ and spin configurations in $\Sigma(m).$
\end{prop}

\noindent Since spin configurations and $N$-regular partitions are in one-to-one correspondence, we now may label the basis $\{ P_{\l}^{(N)} \: | \: \l \in \Pi_N, m(\l) = m\} $  of the $Y(\gl_N)$-submodule $V_N(m)$ by semi-standard tableaux of shape $\theta(m).$  

Consider a ribbon $\theta(m)$ for some $m\in W_N.$ Construct a semi-standard tableau  by filling all columns of $\theta(m)$ with numbers $1,2,\dots,N$ so that a column of height $n$ is filled with  $1,2,\dots,n.$ We will call this tableau the {\em maximal tableau of shape $\theta(m).$} This tableau corresponds to the maximal partition ( i.e. the partition of minimal degree $|\l|$) of the set $\{  \l \in \Pi_N \:|\: m(\l) = m\} $ and hence the $\ssl_N$ weight of the basis vector labelled by this tableau is maximal in the Yangian submodule $V_N(m).$

\begin{eg} Let $N=3.$ Continuing from Example \ref{ex:BS} we find that the basis of the $Y(\gl_N)$-submodule $V_3(m)$ with  $m=((0)^2 (1) (2)^3 (3)^3 \dots )$ is labelled by the following eight semi-standard tableaux  
\unitlength=10pt
\begin{center}


\begin{picture}(6,6)(0,-1) 
\put(2,4){\line(1,0){2}}
\put(0,2){\line(1,0){4}}
\put(0,0){\line(1,0){4}}
\put(0,0){\line(0,1){2}}
\put(2,0){\line(0,1){4}}
\put(4,0){\line(0,1){4}}
\put(0,0){\makebox(2,2){1}}
\put(2,0){\makebox(2,2){1}}
\put(2,2){\makebox(2,2){2}}
\end{picture}  
\begin{picture}(8,4)(0,-4) 
\put(0,0){\line(1,0){2}}
\put(0,-1){\line(1,0){2}}
\put(0,-2){\line(1,0){2}}
\put(0,0){\line(0,-1){2}}
\put(1,0){\line(0,-1){2}}
\put(2,0){\line(0,-1){2}}
\put(0,-1){\makebox(1,1){0}}\put(1,-1){\makebox(1,1){1}}
\put(0,-2){\makebox(1,1){2}}\put(1,-2){\makebox(1,1){0}}
\end{picture} 
\begin{picture}(6,6)(0,-1) 
\put(2,4){\line(1,0){2}}
\put(0,2){\line(1,0){4}}
\put(0,0){\line(1,0){4}}
\put(0,0){\line(0,1){2}}
\put(2,0){\line(0,1){4}}
\put(4,0){\line(0,1){4}}
\put(0,0){\makebox(2,2){1}}
\put(2,0){\makebox(2,2){1}}
\put(2,2){\makebox(2,2){3}}
\end{picture}  
{\begin{picture}(8,4)(0,-4) 
\put(0,0){\line(1,0){3}}
\put(0,-1){\line(1,0){3}}
\put(0,-2){\line(1,0){2}}
\put(0,0){\line(0,-1){2}}
\put(1,0){\line(0,-1){2}}
\put(2,0){\line(0,-1){2}}
\put(3,0){\line(0,-1){1}}
\put(0,-1){\makebox(1,1){0}}\put(1,-1){\makebox(1,1){1}}\put(2,-1){\makebox(1,1){2}}
\put(0,-2){\makebox(1,1){2}}\put(1,-2){\makebox(1,1){0}}
\end{picture}} \\ 
\begin{picture}(6,6)(0,-1) 
\put(2,4){\line(1,0){2}}
\put(0,2){\line(1,0){4}}
\put(0,0){\line(1,0){4}}
\put(0,0){\line(0,1){2}}
\put(2,0){\line(0,1){4}}
\put(4,0){\line(0,1){4}}
\put(0,0){\makebox(2,2){2}}
\put(2,0){\makebox(2,2){1}}
\put(2,2){\makebox(2,2){2}}
\end{picture}  
{\begin{picture}(8,4)(0,-4) 
\put(0,0){\line(1,0){2}}
\put(0,-1){\line(1,0){2}}
\put(0,-2){\line(1,0){2}}
\put(0,-3){\line(1,0){1}}
\put(0,0){\line(0,-1){3}}
\put(1,0){\line(0,-1){3}}
\put(2,0){\line(0,-1){2}}
\put(0,-1){\makebox(1,1){0}}\put(1,-1){\makebox(1,1){1}}
\put(0,-2){\makebox(1,1){2}}\put(1,-2){\makebox(1,1){0}}
\put(0,-3){\makebox(1,1){1}}
\end{picture}}  
\begin{picture}(6,6)(0,-1) 
\put(2,4){\line(1,0){2}}
\put(0,2){\line(1,0){4}}
\put(0,0){\line(1,0){4}}
\put(0,0){\line(0,1){2}}
\put(2,0){\line(0,1){4}}
\put(4,0){\line(0,1){4}}
\put(0,0){\makebox(2,2){2}}
\put(2,0){\makebox(2,2){1}}
\put(2,2){\makebox(2,2){3}}
\end{picture}  
{\begin{picture}(8,4)(0,-4) 
\put(0,0){\line(1,0){3}}
\put(0,-1){\line(1,0){3}}
\put(0,-2){\line(1,0){2}}\put(0,-3){\line(1,0){1}}
\put(0,0){\line(0,-1){3}}
\put(1,0){\line(0,-1){3}}
\put(2,0){\line(0,-1){2}}
\put(3,0){\line(0,-1){1}}
\put(0,-1){\makebox(1,1){0}}\put(1,-1){\makebox(1,1){1}}\put(2,-1){\makebox(1,1){2}}
\put(0,-2){\makebox(1,1){2}}\put(1,-2){\makebox(1,1){0}}\put(0,-3){\makebox(1,1){1}}
\end{picture}} \\ 
\begin{picture}(6,6)(0,-1) 
\put(2,4){\line(1,0){2}}
\put(0,2){\line(1,0){4}}
\put(0,0){\line(1,0){4}}
\put(0,0){\line(0,1){2}}
\put(2,0){\line(0,1){4}}
\put(4,0){\line(0,1){4}}
\put(0,0){\makebox(2,2){3}}
\put(2,0){\makebox(2,2){1}}
\put(2,2){\makebox(2,2){2}}
\end{picture}  
{\begin{picture}(8,4)(0,-4) 
\put(0,0){\line(1,0){2}}
\put(0,-1){\line(1,0){2}}
\put(0,-2){\line(1,0){2}}\put(0,-3){\line(1,0){2}}
\put(0,0){\line(0,-1){3}}
\put(1,0){\line(0,-1){3}}
\put(2,0){\line(0,-1){3}}
\put(0,-1){\makebox(1,1){0}}\put(1,-1){\makebox(1,1){1}}
\put(0,-2){\makebox(1,1){2}}\put(1,-2){\makebox(1,1){0}}
\put(0,-3){\makebox(1,1){1}}\put(1,-3){\makebox(1,1){2}}
\end{picture}}   
\begin{picture}(6,6)(0,-1) 
\put(2,4){\line(1,0){2}}
\put(0,2){\line(1,0){4}}
\put(0,0){\line(1,0){4}}
\put(0,0){\line(0,1){2}}
\put(2,0){\line(0,1){4}}
\put(4,0){\line(0,1){4}}
\put(0,0){\makebox(2,2){2}}
\put(2,0){\makebox(2,2){2}}
\put(2,2){\makebox(2,2){3}}
\end{picture}  
{\begin{picture}(8,4)(0,-4) 
\put(0,0){\line(1,0){3}}
\put(0,-1){\line(1,0){3}}
\put(0,-2){\line(1,0){3}}
\put(0,-3){\line(1,0){1}}
\put(0,0){\line(0,-1){3}}
\put(1,0){\line(0,-1){3}}
\put(2,0){\line(0,-1){2}}
\put(3,0){\line(0,-1){2}}
\put(0,-1){\makebox(1,1){0}}\put(1,-1){\makebox(1,1){1}}\put(2,-1){\makebox(1,1){2}}
\put(0,-2){\makebox(1,1){2}}\put(1,-2){\makebox(1,1){0}}\put(2,-2){\makebox(1,1){1}}
\put(0,-3){\makebox(1,1){1}}
\end{picture}} \\
\begin{picture}(6,6)(0,-1) 
\put(2,4){\line(1,0){2}}
\put(0,2){\line(1,0){4}}
\put(0,0){\line(1,0){4}}
\put(0,0){\line(0,1){2}}
\put(2,0){\line(0,1){4}}
\put(4,0){\line(0,1){4}}
\put(0,0){\makebox(2,2){3}}
\put(2,0){\makebox(2,2){1}}
\put(2,2){\makebox(2,2){3}}
\end{picture}  
{\begin{picture}(8,4)(0,-4) 
\put(0,0){\line(1,0){3}}
\put(0,-1){\line(1,0){3}}
\put(0,-2){\line(1,0){2}}
\put(0,-3){\line(1,0){2}}
\put(0,0){\line(0,-1){3}}
\put(1,0){\line(0,-1){3}}
\put(2,0){\line(0,-1){3}}
\put(3,0){\line(0,-1){1}}
\put(0,-1){\makebox(1,1){0}}\put(1,-1){\makebox(1,1){1}}\put(2,-1){\makebox(1,1){2}}
\put(0,-2){\makebox(1,1){2}}\put(1,-2){\makebox(1,1){0}}
\put(0,-3){\makebox(1,1){1}}\put(1,-3){\makebox(1,1){2}}
\end{picture}}  
\begin{picture}(6,6)(0,-1) 
\put(2,4){\line(1,0){2}}
\put(0,2){\line(1,0){4}}
\put(0,0){\line(1,0){4}}
\put(0,0){\line(0,1){2}}
\put(2,0){\line(0,1){4}}
\put(4,0){\line(0,1){4}}
\put(0,0){\makebox(2,2){3}}
\put(2,0){\makebox(2,2){2}}
\put(2,2){\makebox(2,2){3}}
\end{picture}  
{\begin{picture}(8,4)(0,-4) 
\put(0,0){\line(1,0){3}}
\put(0,-1){\line(1,0){3}}
\put(0,-2){\line(1,0){3}}
\put(0,-3){\line(1,0){2}}
\put(0,0){\line(0,-1){3}}
\put(1,0){\line(0,-1){3}}
\put(2,0){\line(0,-1){3}}
\put(3,0){\line(0,-1){2}}
\put(0,-1){\makebox(1,1){0}}\put(1,-1){\makebox(1,1){1}}\put(2,-1){\makebox(1,1){2}}
\put(0,-2){\makebox(1,1){2}}\put(1,-2){\makebox(1,1){0}}\put(2,-2){\makebox(1,1){1}}
\put(0,-3){\makebox(1,1){1}}\put(1,-3){\makebox(1,1){2}}
\end{picture}}


\end{center}
\unitlength=1pt
where the corresponding (colored) $3$-regular partition is  shown on the right of each tableaux.
The tableau $(2,1,1)$ is the maximal tableau of shape $\theta(m).$ \end{eg}

\subsubsection{Irreducibility and Drinfeld polynomials} Let $m\in W_N$ and let $\theta(m)$ be the corresponding ribbon from the set $R_N.$  With this ribbon we associate the unique pair of partitions $\nu=(\nu_1,\nu_2,\dots)$ and $\mu=(\mu_1,\mu_2,\dots)$ such  that $\theta(m) = \nu/\mu$ and the length of $\mu$ (resp. $\mu'$) is less than the length of $\nu$ (resp. $\nu'$). Let $V_{\nu,\mu}$ be the tame $Y(\gl_N)$-module associated with the skew diagram $\nu/\mu$ in the manner described in Section \ref{sec:tame}. We will denote by  $^{\s}V_{\nu,\mu}$ the tame Yangian module obtained from  $V_{\nu,\mu}$  by the pullback through the automorphism $\s_N$ defined  in (\ref{eq:sigma}). Lemma 6.2 contained in \cite{KKN} implies that the dimension of $V_{\nu,\mu}$ is equal to the number of semi-standard tableaux of shape $\theta(m) = \nu/\mu.$ Hence dimensions of the Yangian modules $V_N(m)$ and $^{\s}V_{\nu,\mu}$ are equal. For a complex number $h$ we  denote  by $^{\s}V_{\nu,\mu}(h)!
!
$ the tame Yangian module obtain
ed from $^{\s}V_{\nu,\mu}$ through the automorphism $\xi(h)$ (\ref{eq:xi}). 

\begin{prop}
The $Y(\gl_N)$-modules $V_N(m)$ and $^{\s}V_{\nu,\mu}(\nu_1-1)$ are isomorphic up to an automorphism of the algebra $Y(\gl_N)$ of the form {\em (\ref{eq:omf}).}
\end{prop}
\begin{pf}
Let $\l^{\max}$ be the unique partition from the set $\{ \l \in \Pi_N \:|\: m(\l)=m \}$ such that its degree $|\l^{\max}| $ is minimal among degrees of all partitions contained in this set. Then the $\ssl_N$ weight of the basis vector $P_{\l^{\max}}^{(N)}$ is maximal among the  $\ssl_N$ weights of the submodule   $V_N(m).$ This implies that for all $N \geq a > b \geq 1$ and $s=1,2,\dots $ we have 
\begin{equation}
 \pi(T_{ab}^{(s)})P_{\l^{\max}}^{(N)} = 0.
\end{equation}
By Proposition \ref{p:piAeigenbasis},  $P_{\l^{\max}}^{(N)}$ is an eigenvector of the maximal commutative subalgebra $\pi(A(\gl_N)).$ Using   the eigenvalue formula given in  Proposition \ref{p:piAeigenbasis} we obtain  for all $m=1,\dots,N-1$ 
\begin{equation}
\frac{\pi(A_{m+1}(u))\pi(A_{m-1}(u-1))}{\pi(A_{m}(u))\pi(A_{m}(u-1))}P_{\l^{\max}}^{(N)} = \frac{P_m(u-1)}{P_m(u)}P_{\l^{\max}}^{(N)}
\end{equation}
where $P_m(u)$ are polynomials given by  
\begin{equation}
P_m(u) =\prod_{c}(u + \nu_1 - c - m ); \quad m=1,\dots,N-1, \label{eq:loc}
\end{equation}
with the product taken over contents  of bottom squares of all  columns of height $m$ in the ribbon $\theta(m) = \nu/\mu.$  

Let $V_N(m)'$ be the submodule of $V_N(m)$ generated by the vector  $ P_{\l^{\max}}^{(N)}.$ Suppose $V_N(m)''$ is the maximal proper submodule of $V_N(m)'.$ 

Without loss of generality we may assume that  $P_{\l^{\max}}^{(N)} \not\in V_N(m)''.$ Then comparing polynomials (\ref{eq:loc}) with the Drinfeld polynomials of the module  $^{\s}V_{\nu,\mu}$ given by Proposition \ref{p:tameDP} and Proposition \ref{p:sigmaDP} we see that  $V_N(m)'/V_N(m)''$ and $ ^{\s}V_{\nu,\mu}(\nu_1-1)$ are isomorphic up to an automorphism of the form (\ref{eq:omf}).  The required statement now follows since dimensions of  $^{\s}V_{\nu,\mu}(\nu_1-1)$ and $V_N(m)$ are equal.  
\end{pf}

\end{document}